\documentclass[twoside,australian,english,3p,sort&compress, review]{elsarticle}
\usepackage[T1]{fontenc}
\usepackage[utf8]{inputenc}
\usepackage{color}
\usepackage{babel}
\usepackage{float}
\usepackage{units}
\usepackage{amsmath}
\usepackage{graphicx}
\usepackage{esint}
\usepackage[unicode=true,
 bookmarks=true,bookmarksnumbered=false,bookmarksopen=false,
 breaklinks=false,pdfborder={0 0 1},backref=false,colorlinks=true]
 {hyperref}

\makeatletter

\newcommand{\lyxdot}{.}

\journal{Computers \& Structures}
\usepackage{lineno}
\usepackage{multirow}
\usepackage{upgreek}
\usepackage{paralist}
\usepackage[section]{placeins}
\usepackage{scalerel}

\usepackage{xfrac}
\usepackage{xcolor}
\usepackage[bottom]{footmisc}

\@ifundefined{showcaptionsetup}{}{%
 \PassOptionsToPackage{caption=false}{subfig}}
\usepackage{subfig}
\makeatother

\begin{document}

\begin{frontmatter}{}

\title{High-order implicit time integration scheme with controllable numerical
dissipation based on mixed-order Padé expansions}

\author[unsw]{Chongmin Song\corref{cor1}}

\ead{c.song@unsw.edu.au}

\author[unsw]{Xiaoran~Zhang}

\author[unsw,tud]{Sascha~Eisenträger}

\author[unsw]{Ankit Ankit}

\cortext[cor1]{Corresponding author}

\address[unsw]{Centre for Infrastructure Engineering and Safety, School of Civil
and Environmental Engineering, University of New South Wales, Sydney,
NSW 2052, Australia.}

\address[tud]{Department of Civil and Environmental Engineering, Institute for Mechanics,
Computational Mechanics Group, Technical University of Darmstadt,
Darmstadt 64287, Germany }
\begin{abstract}
A single-step high-order implicit time integration scheme with controllable
numerical dissipation at high frequencies is presented for the transient
analysis of structural dynamic problems. The amount of numerical dissipation
is controlled by a user-specified value of the spectral radius $\rho_{\infty}$
in the high frequency limit. Using this user-specified parameter as
a weight factor, a Padé expansion of the matrix exponential solution
of the equation of motion is constructed by mixing the diagonal and
sub-diagonal expansions. An efficient time-stepping scheme is designed
where systems of equations, similar in complexity to the standard
Newmark method, are solved recursively. It is shown that the proposed
high-order scheme achieves high-frequency dissipation, while minimizing
low-frequency dissipation and period errors. The effectiveness of
the provided dissipation control and the efficiency of the scheme
are demonstrated by numerical examples. A simple guideline for the
choice of the controlling parameter and time step size is provided.
The source codes written in \texttt{MATLAB} and \texttt{FORTRAN} are
available for download at: \href{https://github.com/ChongminSong/HighOrderTimeIntegration}{https://github.com/ChongminSong/HighOrderTimeIntegration}.
\end{abstract}

\end{frontmatter}{}


\section{Introduction\label{sec:Introduction}}

In several engineering and science disciplines, it is needed to evaluate
the solution of the response history of a system to a dynamic action
\citep{Hernandez2020,Kim2021,Gao2021}. Various direct time integration
methods have been developed to discretize the time-continuous dynamic
equations \citep{Bathe2005,Zienkiewvicz:2005}. A review on the recent
progress has been reported by the authors in Ref.~\citep{Song2022}.
In an explicit method, the response at the current time step is formulated
from the information at previous time steps. Explicit methods can
be designed without the solution of simultaneous equations, but are
generally conditionally stable. In an implicit method, the response
at the current time step is formulated by combining the information
at the current and previous time steps. Commonly used implicit methods
are unconditionally stable, but require the solution of simultaneous
equations. In the present paper, only implicit methods are considered.
The error induced by the temporal discretization of a specific implicit
method is related to the time step size and its order of accuracy.

In the numerical analysis of a dynamic system expressed by governing
partial differential equations, for example wave propagation in a
continuous medium, the problem domain is discretized spatially, resulting
in a system of semi-discrete equations of motion. Numerical methods
such as the finite element method \citep{BookBathe2014,Zienkiewvicz:2005},
the spectral element method \citep{ArticleDuczek2019b}, isogeometric
analysis \citep{Cottrell2006}, overlapping finite elements \citep{Kim2021},
and the scaled boundary finite element method \citep{Song2009,Song2018}
can be employed for the spatial discretization. The spatial discretization
error in a semi-discrete system is related to the number of nodes
per wavelength and the polynomial degree of the shape functions of
the elements. The error increases as the wavelength becomes shorter,
i.e., as the frequency becomes higher. In many engineering applications,
the mesh is chosen to accurately represent the vibration modes below
a maximum frequency of interest. These modes are referred to as lower
modes. The spatial discretization error will become unacceptable for
higher modes that are spatially unresolved. When higher modes are
excited, the solution will be polluted by spurious, i.e., non-physical
oscillations. Therefore, a large number of time integration schemes
possessing numerical damping, for example \citep{Houbolt1950,Newmark1959,Wilson1972,Hilber1977,Chung1993,Bathe2007,Noh2018,Malakiyeh2021,Kim2017a},
have been developed aiming to achieve high-frequency dissipation while
minimizing low-frequency dissipation. The most commonly used implicit
methods in commercial finite element software and scientific applications
include the Houbolt method~\citep{Houbolt1950}, the Wilson-$\theta$
method~\citep{Wilson1972}, the Newmark method~\citep{Newmark1959},
the HHT-$\alpha$ method~\citep{Hilber1977}, the generalized-$\alpha$
method~\citep{Chung1993}, and the Bathe method~\citep{Bathe2005}.
The computer implementations of these methods are straightforward.
The system of simultaneous equations to be solved at a time step is
similar to that in statics. 

Various high-order time integration methods have been proposed recently~
\citep{Kim2017a,Kim2017b,Kim2018,Kim2020,Soares2020,Behnoudfar2020,Behnoudfar2021}.
The Bathe method has been extended to third-order and fourth-order
accuracy in \citep{Kwon2021,Choi2022}. In most methods, the response
history is approximated with polynomial functions in time. Numerical
methods such as collocation, differential quadrature, and weighted-residuals
are applied to derive the time-stepping formulations. Another technique
for constructing highly accurate time integration methods is to apply
the matrix theory for the solution of systems of ordinary differential
equations. Various time-stepping algorithms are designed by approximating
the matrix exponential function in the solution with functions that
are suitable for numerical computation, such as Taylor series, Chebyshev
expansions, and Padé expansions~\citep{Reusch1988,Wang2007,Barucq2018,Gao2021}.
Generally speaking, high-order direct time integration methods are
computationally more expensive than the second-order methods mentioned
above for advancing one time step. However, much larger time step
sizes can be used to obtain solutions of similar accuracy, which may
lead to a more efficient overall performance.

A computationally effective high-order time integration method is
developed in Ref.~\citep{Song2022} by applying the technique of
partial fractions to Padé expansions. Based on the diagonal Padé expansion
of order $(M,M)$ an \emph{A-}stable time-stepping scheme of order
$2M$ is constructed. A salient feature of the approach is that the
simultaneous equations to be solved are similar to that in the standard
Newmark method. When $M$ is odd, one system of equations with a real
matrix and $(M-1)/2$ systems of equations with complex matrices are
solved recursively to advance one time step. When $M$ is even, there
are $M/2$ systems of equations with complex matrices. It is observed
from a \texttt{FORTRAN} implementation using Intel MKL \texttt{PARDISO}
direct solver that a scheme of order $2M$ is approximately only $M$-times
more costly than the Newmark method, which results in significant
gains in terms of efficiency reaching the same accuracy level. However,
the time-stepping scheme based on diagonal Padé expansions does not
possess numerical dissipation.

This paper extends the high-order time integration scheme proposed
in Ref.~\citep{Song2022} to include controllable numerical dissipation
by mixing the diagonal Padé expansion of order $(M,M)$ with the sub-diagonal
Padé expansion of order $(M-1,M).$ The spectral radius in the high-frequency
limit $\rho_{\infty}$ is to be specified by the user as a parameter
in the range between $\rho_{\infty}=0$ (\emph{L-}stable) and $\rho_{\infty}=1$
(\emph{A-}stable) to control the amount of numerical dissipation.
The order of accuracy of the scheme based on the mixed-order Padé
expansion is $2M-1$, when numerical dissipation is specified.

The subsequent development of the paper is organized as follows: Section~\ref{sec:Summary-time-stepping-schemes}
summarizes the matrix exponential solution of the equation of motion.
Section~\ref{sec:Mixed-order-Pad=0000E9-expansion} explains the
construction of the mixed-order Padé expansion with a user-specified
parameter. In Section~\ref{sec:Time-stepping-schemes}, a computationally
effective time-stepping scheme is described. The numerical dissipation
and dispersion of the proposed scheme are analyzed in Section~\ref{sec:Analysis-of-time-stepping}.
Numerical examples are presented in Section~\ref{sec:Numerical-examples}
to demonstrate the effectiveness of numerical dissipation and the
efficiency of the proposed scheme. The selection of the user-specified
parameter $\rho_{\infty}$ for controlling the amount of numerical
dissipation and time step size $\Delta t$ are also investigated.
Finally, conclusions are drawn in Section~\ref{sec:Conclusions}.

\section{Summary of time-stepping using matrix exponential\label{sec:Summary-time-stepping-schemes}}

\subsection{Development of time integration method}

In this section, the theory used in developing the single-step high-order
implicit time integration scheme featuring controllable numerical
dissipation is presented. The method is based on rewriting the equation
of motion as a system of first-order ordinary differential equations
(ODEs) in state-space and approximating the matrix exponential function
in the exact solution by Padé expansions.

The equation of motion in structural dynamic problems can be expressed
as a system of second-order ODEs and is written as
\begin{equation}
\mathbf{M\ddot{u}}(t)+\mathbf{C\dot{u}}(t)+\mathbf{Ku}(t)=\mathbf{f}(t)\label{eq: equation of motion}
\end{equation}
with the initial conditions
\begin{equation}
\mathbf{u}(t=0)=\mathbf{u}_{0}\,,\label{eq:initial displacement}
\end{equation}
\begin{equation}
\dot{\mathbf{u}}(t=0)=\dot{\mathbf{u}}_{0}\,,\label{eq:initial acceleration}
\end{equation}
where $\mathbf{M},\mathbf{\,C}$, and $\mathbf{K}$ denote the mass,
damping, and stiffness matrices, respectively. $\mathbf{f}$ is the
external excitation force vector, and $\mathbf{u}\mathbf{\mathrm{,\,}\boldsymbol{\mathbf{\dot{u}}}}$
and $\mathbf{\ddot{u}}$ represent displacement, velocity, and acceleration
vectors, respectively.

In a time-stepping scheme, the overall time duration is divided into
a finite number of time intervals. Without loss of generality, the
scheme is described for a time step $n$ over the interval $t_{n-1}\leq t\leq t_{n}$
($n=1,2,\ldots,n_{S}$). The size of the time step is denoted as $\Delta t=t_{n}-t_{n-1}.$
By introducing a dimensionless time variable $s,$ for each time step,
the time within the time step $n$ can be determined by 
\begin{equation}
t(s)=t_{n-1}+s\Delta t,\quad\;0\leq s\leq1\,,\label{eq:dimensionless time}
\end{equation}
with $t(s=0)=t_{n-1}$ at the beginning of the time step and $t(s=1)=t_{n}$
at the end of the time step.

Denoting the derivative with respect to the dimensionless time $s$
by a circle ($\circ$) above the symbol, the velocity and acceleration
vectors within the time step $t_{n-1}\leq t\leq t_{n}$ can be written
as
\begin{equation}
\mathbf{\dot{u}}=\frac{1}{\Delta t}\frac{\mathrm{d\boldsymbol{u}}}{\mathrm{d}s}=\frac{1}{\Delta t}\stackrel{{\scriptstyle {\scriptscriptstyle \circ}}}{\mathbf{u}}\,,\label{eq:velocity in dimensionless time}
\end{equation}
\begin{equation}
\mathbf{\ddot{u}}=\frac{1}{\Delta t^{2}}\frac{\mathrm{d^{2}\mathbf{u}}}{\mathrm{d}s^{2}}=\frac{1}{\Delta t^{2}}\stackrel{{\scriptstyle {\scriptscriptstyle \circ\circ}}}{\mathbf{u}}\,.\label{eq:acceleration in dimensinoless time}
\end{equation}
Therefore, the equation of motion can be expressed in the dimensionless
time as
\begin{equation}
\mathbf{M}\stackrel{{\scriptstyle {\scriptscriptstyle \circ\circ}}}{\mathbf{u}}+\Delta t\mathbf{C}\stackrel{{\scriptstyle {\scriptscriptstyle \circ}}}{\mathbf{u}}+\Delta t^{2}\mathbf{K}\mathbf{u}=\Delta t^{2}\mathbf{f}\,.\label{eq:eqn of motion in dimensionless}
\end{equation}

Introducing a state-space vector $\boldsymbol{z}$ defined as
\begin{equation}
\mathbf{z}=\left\{ \begin{array}{c}
\stackrel{{\scriptstyle {\scriptscriptstyle \circ}}}{\mathbf{u}}\\
\mathbf{u}
\end{array}\right\} \,,\label{eq:state space vector}
\end{equation}
Eq.~\eqref{eq:eqn of motion in dimensionless} can be transformed
into a system of first-order ODEs
\begin{equation}
\stackrel{{\scriptstyle {\scriptscriptstyle \circ}}}{\mathbf{z}}\equiv\frac{\mathrm{d\mathbf{z}}}{\mathrm{d}s}=\mathbf{Az+F}\,,\label{eq:1st order ODE}
\end{equation}
where $\mathbf{A}$ is the constant coefficient matrix defined as
\begin{equation}
\mathbf{A}=\left[\begin{array}{cc}
-\Delta t\mathbf{M^{-1}C} & -\Delta t^{2}\mathbf{M^{-1}K}\\
\mathbf{I} & \mathbf{0}
\end{array}\right]\,,\label{eq:matrix A}
\end{equation}
 and $\mathbf{F}$ is the non-homogeneous term
\begin{equation}
\mathbf{F}=\left\{ \begin{array}{c}
\Delta t^{2}\mathbf{M^{-1}f}\\
\mathbf{0}
\end{array}\right\} \,.\label{eq:force vector}
\end{equation}

The general solution of Eq.~\eqref{eq:1st order ODE} at time $t_{n}$
is obtained with the matrix exponential function as

\begin{equation}
\mathbf{z}_{n}=e^{\mathbf{A}s}\mathbf{z}_{n-1}+e^{\mathbf{A}s}\intop_{0}^{s}e^{-\mathbf{A}\tau}\mathbf{F}(\tau)\mathsf{d}\tau\,.\label{eq:general solution}
\end{equation}
The force vector $\mathbf{f}$, and thus $\mathbf{F}$, is expressed
as a polynomial expansion at the middle of the time step $n$
\begin{equation}
\mathbf{F}_{n}(s)=\sum\limits _{k=0}^{p_{\mathrm{f}}}\tilde{\mathbf{F}}_{\mathrm{m}n}^{(k)}(s-0.5)^{k}=\tilde{\mathbf{F}}_{\mathrm{m}n}^{(0)}+\tilde{\mathbf{F}}_{\mathrm{m}n}^{(1)}(s-0.5)+\tilde{\mathbf{F}}_{\mathrm{m}n}^{(2)}(s-0.5)^{2}+\ldots+\tilde{\mathbf{F}}_{\mathrm{m}n}^{(p_{\mathrm{f}})}(s-0.5)^{p_{\mathrm{f}}}\,.\label{eq:ForceExpansion}
\end{equation}
The solution in Eq.\ \eqref{eq:general solution} at the end of the
time step ($s=1$) is simplified to
\begin{equation}
\mathbf{z}_{n}=e^{\mathbf{A}}\mathbf{z}_{n-1}+\sum\limits _{k=0}^{p_{\mathrm{f}}}\mathbf{B}_{k}\tilde{\mathbf{F}}_{\mathrm{m}n}^{(k)}\,,\label{eq:sol s=00003D1}
\end{equation}
where $\mathbf{B}_{k}$ is integrated by parts and can be determined
recursively 
\begin{equation}
\mathbf{B}_{k}=e^{\mathbf{A}}\int\limits _{0}^{1}(\tau-0.5)^{k}e^{-\mathbf{A}\tau}\mathrm{d}\tau=\mathbf{A}^{-1}\left(k\mathbf{B}_{k-1}+\left(-\cfrac{1}{2}\right)^{k}(e^{\mathbf{A}}-(-1)^{k}\mathbf{I})\right)\,,\quad\;\forall k=0,1,2,\ldots,p_{\mathrm{f}}\label{eq:intergrationExpmI}
\end{equation}
with the starting value at $k=0$ 
\begin{equation}
\mathbf{B}_{0}=e^{\mathbf{A}}\int\limits _{0}^{1}e^{-\mathbf{A}\tau}\mathrm{d}\tau=\mathbf{A}^{-1}\left(e^{\mathbf{A}}-\mathbf{I}\right)\,.\label{eq:integrationExpm0}
\end{equation}

\section{Mixed-order Padé expansion of matrix exponential\label{sec:Mixed-order-Pad=0000E9-expansion}}

\noindent Computing the matrix exponential $e^{\mathbf{A}}$ in Eq.~\eqref{eq:sol s=00003D1}
generally results in a full matrix. For practical engineering problems,
the operation is expensive in terms of computational time and memory
\citep{Golub1996}. To derive an efficient time-stepping scheme, the
direct computation of the matrix exponential is avoided by employing
a rational approximation $\mathbf{R}=\mathbf{R}(\mathbf{A})$, i.e.,
a ratio of two polynomials
\begin{equation}
e^{\mathbf{A}}\approx\mathbf{R}=\mathbf{\dfrac{\mathbf{P}}{Q}}\,,\label{eq:ExpmRational}
\end{equation}
where $\mathbf{P}=\mathbf{P}(\mathbf{A})$ and $\mathbf{Q}=\mathbf{Q}(\mathbf{A})$
are polynomials of matrix $\mathbf{A}$, expressed as\begin{subequations}\label{eq:polyPQ}
\begin{align}
\mathbf{P} & =\sum\limits _{i=0}^{N_{p}}p_{i}\mathbf{A}^{i}=p_{0}\mathbf{I}+p_{1}\mathbf{A}+\ldots+p_{N_{p}}\mathbf{A}^{N_{p}}\,,\label{eq:polyP}\\
\mathbf{Q} & =\sum\limits _{i=0}^{N_{q}}q_{i}\mathbf{A}^{i}=q_{0}\mathbf{I}+q_{1}\mathbf{A}+\ldots+q_{N_{q}}\mathbf{A}^{N_{q}}\,.\label{eq:polyQ}
\end{align}
\end{subequations}Here, $N_{p}$ and $N_{q}$ denote the degrees
of $\mathbf{P}$ and $\mathbf{Q}$, respectively, and the scalar coefficients
$p_{i}$ and $q_{i}$ are real. To ensure that $e^{\mathbf{0}}=\mathbf{I}$
holds, $p_{0}=q_{0}$ applies. Note that the matrix product $\mathbf{Q}^{-1}\mathbf{P}$
is commutative (i.e., $\mathbf{P}\mathbf{Q}^{-1}=\mathbf{Q}^{-1}\mathbf{\mathbf{P}}$).

The Padé expansion of a function is obtained by determining the coefficients
of the two polynomials to obtain the highest order of accuracy. For
the matrix exponential function $e^{\mathbf{\mathbf{A}}}$, the Padé
expansion of order $(L,M)$ is written as 
\begin{equation}
e^{\mathbf{A}}\approx e_{L/M}^{\mathbf{A}}=\cfrac{\boldsymbol{\mathbf{P}}_{L/M}(\mathbf{A})}{\mathbf{Q}_{L/M}(\mathbf{A})}\,,\label{eq:Pade expansion}
\end{equation}
where $\mathbf{P}_{L/M}(\mathbf{A})$ and $\mathbf{Q}_{L/M}(\mathbf{A})$
are polynomials of order $L$ and $M$, respectively, given by \begin{subequations}\label{eq:PQ_L/M}
\begin{equation}
\mathbf{P}_{L/M}(\mathbf{A})=\sum\limits _{i=0}^{L}\cfrac{(M+L-i)!}{i!(L-i)!}\mathbf{A}^{i}\,,\label{eq:P_L/M}
\end{equation}
\begin{equation}
\mathbf{Q}_{L/M}(\mathbf{A})=\frac{M!}{L!}\sum\limits _{i=0}^{M}\cfrac{(M+L-i)!}{i!(M-i)!}(-\mathbf{A})^{i}\,.
\end{equation}
\end{subequations} The truncation error is of order $O(\mathbf{A}^{L+M+1})$.

As it will be shown in Section~\ref{sec:Analysis-of-time-stepping},
the diagonal Padé expansions, in which the polynomial orders of the
numerator and denominator are identical ($L=M$), are \emph{A-}stable
and do not exhibit any numerical dissipation, while the sub-diagonal
Padé expansions ($L<M$) are \emph{L-}stable, leading to strong numerical
dissipation.

To develop a time integration scheme with controllable numerical dissipation,
a mixed order Padé expansion is constructed. The spectral radius in
the high-frequency limit $\rho_{\infty}$ is selected as the controlling
parameter. Equation~\eqref{eq:Pade expansion} is rewritten for the
diagonal Padé expansion of order $(M,M)$ and sub-diagonal expansion
of order $(L<M,M)$ as \begin{subequations}
\begin{align}
\mathbf{Q}_{M/M}(\mathbf{A})e^{\mathbf{A}} & =\mathbf{P}_{M/M}(\mathbf{A})+O(\mathbf{A}^{2M+1})\,,\label{eq:Pade-MM-error}\\
\mathbf{Q}_{L/M}(\mathbf{A})e^{\mathbf{A}} & =\mathbf{P}_{L/M}(\mathbf{A})+O(\mathbf{A}^{L+M+1})\,.\label{eq:Pade-LM-error}
\end{align}
\end{subequations}The degrees of the denominators of the two Padé
expansions are chosen to be the same, i.e., $M$, in the present work.
Summing Eqs.~\eqref{eq:Pade-MM-error} and \eqref{eq:Pade-LM-error}
by applying $\rho_{\infty}$ and $(1-\rho_{\infty})$, respectively,
as the weights leads to
\begin{equation}
\rho_{\infty}\mathbf{Q}_{M/M}(\mathbf{A})e^{\mathbf{A}}+(1-\rho_{\infty})\mathbf{Q}_{L/M}(\mathbf{A})e^{\mathbf{A}}=\rho_{\infty}\mathbf{P}_{M/M}(\mathbf{A})+(1-\rho_{\infty})\mathbf{P}_{L/M}(\mathbf{A})+O(\mathbf{A}^{L+M+1})\,.\label{eq:Pade-mixed-error}
\end{equation}
Defining two polynomial functions \begin{subequations}\label{eq:Pade-mixed-PQ}
\begin{align}
\mathbf{P} & =\rho_{\infty}\mathbf{P}_{M/M}(\mathbf{A})+(1-\rho_{\infty})\mathbf{P}_{L/M}(\mathbf{A})\,,\label{eq:Pade-mixed-P}\\
\mathbf{Q} & =\rho_{\infty}\mathbf{Q}_{M/M}(\mathbf{A})+(1-\rho_{\infty})\mathbf{Q}_{L/M}(\mathbf{A})\,,\label{eq:Pade-mixed-Q}
\end{align}
\end{subequations}a rational approximation of matrix exponential
$e^{\mathbf{A}}$ – see Eq.~\eqref{eq:ExpmRational} with Eq. \eqref{eq:polyPQ}
– is obtained from Eq.~\eqref{eq:Pade-mixed-error} as
\begin{equation}
e^{\mathbf{A}}=\mathbf{\dfrac{\mathbf{P}}{Q}}+O(\mathbf{A}^{L+M+1})\,.\label{eq:Pade-mixed}
\end{equation}
The degree $N_{q}$ of the denominator $\mathbf{Q}$ is equal to $M$.
The degree $N_{p}$ of the numerator $\mathbf{P}$ is equal to $L$
when $\rho_{\infty}=0$, and equal to $M$ otherwise. The rational
function
\begin{equation}
\mathbf{R}=\mathbf{\dfrac{\mathbf{P}}{Q}}\label{eq:Pade-mixed-def}
\end{equation}
with Eq.~\eqref{eq:Pade-mixed-PQ} is referred to as a mixed-order
Padé expansion. When $\rho_{\infty}=1$ is chosen, the diagonal Padé
expansion with the error order $O(\mathbf{A}^{2M+1})$ is recovered.
For any other value $0\leq\rho_{\infty}<1,$ the order of error is
at $O(\mathbf{A}^{L+M+1})$, dominated by the sub-diagonal expansion,
which leads to strong numerical dissipation.

\section{Time-stepping scheme\label{sec:Time-stepping-schemes}}

The mixed-order Padé expansions in Section~\ref{sec:Mixed-order-Pad=0000E9-expansion}
are employed to construct a high-order time-stepping scheme with controllable
numerical dissipation by extending the algorithm presented in Ref.~\citep{Song2022}
and therefore, only key equations are summarized below.

Using Eq.~\eqref{eq:Pade-mixed} and pre-multiplying with $\mathbf{Q}$,
Eq.~\eqref{eq:sol s=00003D1} is reformulated as 
\begin{equation}
\mathbf{Q}\mathbf{z}_{n}=\mathbf{P}\mathbf{z}_{n-1}+\sum\limits _{k=0}^{p_{\mathrm{f}}}\mathbf{C}_{k}\tilde{\mathbf{F}}_{\mathrm{m}n}^{(k)}\,,\label{eq:PadeStepping}
\end{equation}
where the matrices $\mathbf{C}_{k}$ follow from Eqs.~\eqref{eq:intergrationExpmI}
and \eqref{eq:Pade-mixed} as 
\begin{equation}
\mathbf{C}_{k}=\mathbf{Q}\mathbf{B}_{k}=\mathbf{A}^{-1}\left(k\mathbf{C}_{k-1}+\left(-\cfrac{1}{2}\right)^{k}(\mathbf{P}-(-1)^{k}\mathbf{Q})\right)\,,\quad\;\forall k=1,2,\ldots,p_{\mathrm{f}}\,,\label{eq:PadeStepping_C}
\end{equation}
which can be evaluated recursively starting from 
\begin{equation}
\mathbf{C}_{0}=\mathbf{Q}\mathbf{B}_{0}=\mathbf{A}^{-1}\left(\mathbf{P}-\mathbf{Q}\right)\,.\label{eq:PadeStepping_C0}
\end{equation}
where $\mathbf{B}_{0}$ in Eq.~\eqref{eq:integrationExpm0} has been
substituted into.

To develop an efficient algorithm, the polynomial $\mathbf{Q}$ is
factorized as 
\begin{equation}
\mathbf{Q}=\left(r_{1}\mathbf{I}-\mathbf{A}\right)\left(r_{2}\mathbf{I}-\mathbf{A}\right)\ldots\left(r_{M}\mathbf{I}-\mathbf{A}\right)\,.\label{eq:Qfactor}
\end{equation}
The roots $r$ are either real or pairs of complex conjugates. Using
Eq.~\eqref{eq:Qfactor}, Eq.~\eqref{eq:PadeStepping} is rewritten
as 
\begin{equation}
\left(r_{1}\mathbf{I}-\mathbf{A}\right)\left(r_{2}\mathbf{I}-\mathbf{A}\right)\ldots\left(r_{M}\mathbf{I}-\mathbf{A}\right)\mathbf{z}_{n}=\mathbf{b}_{n}\,,\label{eq:PadeSteppingEq}
\end{equation}
where the right-hand side is expressed as 
\begin{equation}
\mathbf{b}_{n}=\mathbf{P}\mathbf{z}_{n-1}+\sum\limits _{k=0}^{p_{\mathrm{f}}}\mathbf{C}_{k}\tilde{\mathbf{F}}_{\mathrm{m}n}^{(k)}\,.\label{eq:PadeStepping_b}
\end{equation}
Equation~\eqref{eq:PadeSteppingEq} is reformulated as series of
equations linear in matrix $\mathbf{A}$ by introducing the auxiliary
variables $\mathbf{z}^{(k)}$ ($k\in\{1,2,\ldots,M-1\}$) 
\begin{equation}
\begin{split}\left(r_{1}\mathbf{I}-\mathbf{A}\right)\mathbf{z}^{(1)} & =\mathbf{b}_{n}\,,\\
\left(r_{2}\mathbf{I}-\mathbf{A}\right)\mathbf{z}^{(2)} & =\mathbf{z}^{(1)}\,,\\
\cdots\\
\left(r_{M}\mathbf{I}-\mathbf{A}\right)\mathbf{z}_{n} & =\mathbf{z}^{(M-1)}\,.
\end{split}
\label{eq:SteppingSuccessive}
\end{equation}
The equations are solved successively by considering one root at a
time.

When a root is real, the corresponding line in Eq.~\eqref{eq:SteppingSuccessive}
is denoted as 
\begin{equation}
\left(r\mathbf{I}-\mathbf{A}\right)\mathbf{x}=\mathbf{g}\,,\label{eq:realRootEq}
\end{equation}
where $r$ is the real root. The unknown vector $\mathbf{x}$ is determined
in relation to a given right-hand side denoted by $\mathbf{g}$. Partitioning
$\mathbf{x}$ and $\mathbf{g}$ into two sub-vectors of equal size
\begin{equation}
\mathbf{x}=\begin{Bmatrix}\mathbf{x}_{1}\\
\mathbf{x}_{2}
\end{Bmatrix}\qquad\text{and}\qquad\mathbf{g}=\begin{Bmatrix}\mathbf{g}_{1}\\
\mathbf{g}_{2}
\end{Bmatrix}\,,\label{eq:realRootPartition}
\end{equation}
and using Eq.~\eqref{eq:matrix A}, Eq.~\eqref{eq:realRootEq} is
rewritten as
\begin{equation}
\left(r^{2}\mathbf{M}+r\Delta t\mathbf{C}+\Delta t^{2}\mathbf{K}\right)\mathbf{x}_{1}=r\mathbf{M}\mathbf{g}_{1}-\Delta t^{2}\mathbf{K}\mathbf{g}_{2}\label{eq:realRootEq5}
\end{equation}
for the solution of $\mathbf{x}_{1}$, and
\begin{equation}
\mathbf{x}_{2}=\cfrac{1}{r}(\mathbf{x}_{1}+\mathbf{g}_{2})\label{eq:realRootEq4}
\end{equation}
for determining $\mathbf{x}_{2}$.

For a pair of complex conjugate roots $r$ and $\bar{r}$ (the overbar
indicates a complex conjugate), the two equations are considered together
and expressed for the unknown vector $\mathbf{x}$ as 
\begin{equation}
\left(r\mathbf{I}-\mathbf{A}\right)\left(\bar{r}\mathbf{I}-\mathbf{A}\right)\mathbf{x}=\mathbf{g}\,.\label{eq:ComplexRootEq}
\end{equation}
Another auxiliary vector $\mathbf{y}$ is introduced 
\begin{equation}
\left(r\mathbf{I}-\mathbf{A}\right)\mathbf{y}=\mathbf{g}\,.\label{eq:ComplexRootEqy}
\end{equation}
Partitioning the vector $\mathbf{y}$ in the same way 
\begin{equation}
\mathbf{y}=\begin{Bmatrix}\mathbf{y}_{1}\\
\mathbf{y}_{2}
\end{Bmatrix}\,,\label{eq:ComplexRoot_sln_y}
\end{equation}
the sub-vector $\mathbf{y}_{1}$ is obtained from 
\begin{equation}
\left(r^{2}\mathbf{M}+r\Delta t\mathbf{C}+\Delta t^{2}\mathbf{K}\right)\mathbf{y}_{1}=r\mathbf{M}\mathbf{g}_{1}-\Delta t^{2}\mathbf{K}\mathbf{g}_{2}\,,\label{eq:ComplexRoot_sln_y1}
\end{equation}
while the second sub-vector $\mathbf{y}_{2}$ follows as
\begin{equation}
\mathbf{y}_{2}=\cfrac{1}{r}\left(\mathbf{y}_{1}+\mathbf{g}_{2}\right)\,.\label{eq:ComplexRoot_sln_y2}
\end{equation}
The solution of Eq.~\eqref{eq:ComplexRootEq} for the unknown vector
$\mathbf{x}$ is expressed as

\begin{equation}
\mathbf{x}=\cfrac{-1}{2\mathrm{Im}(r)\mathrm{i}}\left(\mathbf{y}-\bar{\mathbf{y}}\right)=-\cfrac{\mathrm{Im}(\mathbf{y})}{\mathrm{Im}(r)}\,.\label{eq:ComplexRoot-slnx}
\end{equation}

It is worthwhile to note that the systems of algebraic equations~\eqref{eq:realRootEq5}
and \eqref{eq:ComplexRoot_sln_y1} are in the same form as that found
in the Newmark time-stepping scheme.

\section{Numerical properties of time-stepping scheme for undamped systems\label{sec:Analysis-of-time-stepping}}

The numerical properties, such as numerical dissipation and dispersion,
of the proposed scheme are analyzed in this section. The discussion
is limited to the cases of mixing the diagonal ($M,M$) and the first
sub-diagonal ($M-1,M$) Padé expansions. These cases allow controllable
numerical dissipation at the cost of decreasing the order of accuracy
by one, i.e., from $2M$ to ($2M-1$).

It is also possible to mix the second sub-diagonal ($M-2,M$) expansions
to introduce even stronger numerical dissipation at the cost of reducing
the accuracy to ($2M-2$). Mixing more than two Padé expansions can
be considered as well. Although the numeri

cal properties of such mixed orders are not discussed here, the \texttt{MATHEMATICA}
functions provided in this section can be used for this purpose, if
needed.

\subsection{Discrete-time solution and analysis\label{subsec:Analysis-of-discrete-time}}

The free vibration case of an undamped system, for which $\mathbf{C}=\mathbf{0}$
and $\mathbf{f}=\mathbf{0}$ in Eq.~\eqref{eq: equation of motion}
apply, is considered
\begin{equation}
\mathbf{M\ddot{u}}(t)+\mathbf{Ku}(t)=\mathbf{0}\label{eq:free vibration of undamped system}
\end{equation}
in the evaluation of the dissipative and dispersive characteristics
of the time-stepping scheme. Equation~\eqref{eq:free vibration of undamped system}
can be reduced to a

series of independent single-degree-of-freedom systems using the eigenvalue
problem
\begin{equation}
\mathbf{K}\boldsymbol{\boldsymbol{\phi}}=\omega^{2}\mathbf{M}\boldsymbol{\boldsymbol{\phi}}\,,\label{eq:Analysis-eigenKM}
\end{equation}
where $\omega$ denotes the natural frequency and $\boldsymbol{\boldsymbol{\phi}}$
is the eigenvector. Hence, it is sufficient to consider a single-degree-of-freedom
equation with the natural frequency $\omega$ and the period of vibration
\begin{equation}
T=\dfrac{2\pi}{\omega}\,.\label{eq:Analysis-Tn}
\end{equation}

Following the procedures described in Section~\ref{sec:Summary-time-stepping-schemes}
(Eqs.~\eqref{eq: equation of motion} to \eqref{eq:force vector}),
Eq.~\eqref{eq:free vibration of undamped system} is expressed as
a system of first-order ODEs
\begin{equation}
\stackrel{{\scriptstyle {\scriptscriptstyle \circ}}}{\mathbf{z}}=\mathbf{Az}\label{eq:Analysis-ODE}
\end{equation}
with the coefficient matrix
\begin{equation}
\mathbf{A}=\left[\begin{array}{cc}
\mathbf{0} & -\Delta t^{2}\mathbf{M^{-1}K}\\
\mathbf{I} & \mathbf{0}
\end{array}\right]\,.\label{eq:Analysis-A}
\end{equation}

The eigenvalue problem of matrix \textbf{A} in Eq.~\eqref{eq:matrix A}
is expressed as
\begin{equation}
\mathbf{A}\boldsymbol{\boldsymbol{\psi}}=\lambda\boldsymbol{\boldsymbol{\psi}}\,.\label{eq:Analysis-eigenA}
\end{equation}
Considering Eq.~\eqref{eq:Analysis-eigenKM}, the eigenvalues $\lambda$
and eigenvectors $\boldsymbol{\boldsymbol{\psi}}$ are identified
as\begin{subequations}\label{eq:Analysis-A-eigen-sln}
\begin{align}
\lambda & =\pm\mathrm{i}\omega\Delta t\,,\label{eq:eq:Analysis-A-eigenvalue}\\
\boldsymbol{\boldsymbol{\psi}} & =\left\{ \begin{array}{c}
\pm(\mathrm{i}\omega\Delta t)\boldsymbol{\boldsymbol{\phi}}\\
\boldsymbol{\boldsymbol{\phi}}
\end{array}\right\} \,.\label{eq:Analysis-A-eigenvector}
\end{align}
\end{subequations}All the eigenvalues are purely imaginary and form
pairs of complex conjugates ($\lambda$, $\bar{\lambda}=-\lambda$).
The eigenvectors are also pairs of complex conjugates ($\boldsymbol{\boldsymbol{\psi}}$,
$\bar{\boldsymbol{\boldsymbol{\psi}}}$).

Equation~\eqref{eq:Analysis-ODE} is decoupled by introducing the
transformation
\begin{equation}
\mathbf{z}=\boldsymbol{\Psi}\mathbf{w}\,,\label{eq:Analysis-w-def}
\end{equation}
where the eigenvector matrix $\boldsymbol{\Psi}$ consists of all
the eigenvectors as individual columns, and $\mathbf{w}$ are the
generalized coordinates. Using Eq.~\eqref{eq:Analysis-w-def} and
the technique of eigen-decomposition, Eq.~\eqref{eq:Analysis-ODE}
is decoupled into a series of independent ODEs expressed as
\begin{equation}
\stackrel{{\scriptstyle {\scriptscriptstyle \circ}}}{w}=\lambda w\,.\label{eq:Analytical-w-ODE}
\end{equation}
For every natural frequency $\omega$ in Eq.~\eqref{eq:Analysis-eigenKM},
there are a pair of complex conjugate eigenvalues (Eq.~\eqref{eq:Analysis-A-eigen-sln}).
The vibration of a single-degree-of-freedom system is thus described
by \begin{subequations}\label{eq:Analytical-w-ODE-omega}
\begin{align}
\stackrel{{\scriptstyle {\scriptscriptstyle \circ}}}{w}= & +(\mathrm{i}\omega\Delta t)w\,,\label{eq:Analytical-w-ODE-omega-1}\\
\stackrel{{\scriptstyle {\scriptscriptstyle \circ}}}{\bar{w}}= & -(\mathrm{i}\omega\Delta t)\bar{w}\,.\label{eq:Analytical-w-ODE-omega-2}
\end{align}
\end{subequations}For the evaluation of the dissipative and dispersive
characteristics, it is sufficient to consider one of the two general
coordinates (or eigenvalues) since they are complex conjugates. Equation~\eqref{eq:Analytical-w-ODE-omega-1}
is chosen in the following with
\begin{equation}
\lambda=\mathrm{i}\omega\Delta t\,.\label{eq:Analytical-lambda}
\end{equation}

The continuous-time solution of Eq.~\eqref{eq:Analytical-w-ODE-omega-1}
is written as 
\begin{equation}
w=ce^{\lambda s}=ce^{\mathrm{i}\omega\Delta ts}\,,\label{eq:Analytical-w-sln-general}
\end{equation}
where $c$ is an integration constant. Considering one time step $0\leq s\leq1$
(Eq.~\eqref{eq:dimensionless time}) and using the response $w_{n-1}$
at the start of the time step ($s=0$) as the initial condition, the
response at the end of time step ($s=1$) is obtained as 
\begin{equation}
w_{n}=w_{n-1}e^{\mathrm{i}\omega\Delta t}=w_{n-1}e^{\mathrm{i}2\pi\frac{\Delta t}{T}}\,,\label{eq:Analytical-stepping-exact}
\end{equation}
where Eq.~\eqref{eq:Analysis-Tn} has been substituted into. Applying
Eq.~\eqref{eq:Analytical-w-sln-general} from $t=0$ with the initial
condition $w(0)=w_{0}$ and considering Eq.~\eqref{eq:dimensionless time}
result in
\begin{equation}
w=w_{0}e^{\mathrm{i}\omega t}=w_{0}e^{\mathrm{i}2\pi\frac{t}{T}}\,.\label{eq:Analytical-w-sln}
\end{equation}
Note that the amplitude of vibration remains constant. The phase of
the harmonic vibration is described by $e^{\mathrm{i}2\pi t/T}$. 

The time-stepping formulation in Eq.~\eqref{eq:PadeStepping} is
expressed for the free vibration case in Eq.~\eqref{eq:Analysis-ODE}
as 
\begin{equation}
\mathbf{Q}\mathbf{z}_{n}=\mathbf{P}\mathbf{z}_{n-1}\,.\label{eq:Analytical-stepping}
\end{equation}
With Eq.~\eqref{eq:Analysis-eigenA}, the eigenvalue of an integer
power $i$ of $\mathbf{A}$ is found to be equal to the $i$-th power
of its eigenvalue
\begin{equation}
\mathbf{A}^{i}\boldsymbol{\boldsymbol{\psi}}=\mathbf{A}^{i-1}\mathbf{A}\boldsymbol{\boldsymbol{\psi}}=\lambda\mathbf{A}^{i-1}\boldsymbol{\boldsymbol{\psi}}=\ldots=\lambda^{i}\boldsymbol{\boldsymbol{\psi}}\,.\label{eq:eigen-powerA-1}
\end{equation}
Therefore, the eigen-solution of the polynomial function of matrix
$\mathbf{A}$ in Eq.~\eqref{eq:polyPQ} is expressed as\begin{subequations}\label{eq:Analytical-eigen-PQ}
\begin{align}
\mathbf{P}\boldsymbol{\boldsymbol{\psi}} & =P(\lambda)\boldsymbol{\boldsymbol{\psi}}\,,\label{eq:Analytical-eigen-P}\\
\mathbf{Q}\boldsymbol{\boldsymbol{\psi}} & =Q(\lambda)\boldsymbol{\boldsymbol{\psi}}\,,\label{eq:Analytical-eigen-Q}
\end{align}
\end{subequations}where $P(\lambda)$ and $Q(\lambda)$ are the same
polynomial functions with the eigenvalue $\lambda$ replacing matrix
$\mathbf{A}$.

Introducing the generalized coordinates (Eq.~\eqref{eq:Analysis-w-def})
and using Eq.~\eqref{eq:Analytical-eigen-PQ}, the time-stepping
formulation in Eq.~\eqref{eq:Analytical-stepping} is decomposed
as
\begin{equation}
Q(\lambda)\tilde{w}_{n}=P(\lambda)\tilde{w}_{n-1}\,,\label{eq:eq:Analytical-stepping-decoupled}
\end{equation}
with $\tilde{w}$ denoting the approximate discrete-time solution.
As in the continuous-time solution, it is sufficient to consider the
general coordinate corresponding to $\lambda=\mathrm{i}\omega\Delta t$
(Eq.~\eqref{eq:Analytical-lambda}). In the analysis of time-stepping
schemes, it is customary to use the ratio of the size of time step
$\Delta t$ to the period $T$ and consequently, the eigenvalue $\lambda$
is written as 
\begin{equation}
\lambda=\mathrm{i}2\pi\dfrac{\Delta t}{T}\,.\label{eq:Analytical-lambda-T}
\end{equation}
The time-stepping scheme in Eq.~\eqref{eq:eq:Analytical-stepping-decoupled}
is expressed as
\begin{equation}
\tilde{w}_{n}=R(\lambda)\tilde{w}_{n-1}\,,\label{eq:Analytical-stepping-w-eq}
\end{equation}
where the amplification factor is expressed as
\begin{equation}
R(\lambda)=\dfrac{P(\lambda)}{Q(\lambda)}\label{eq:Analytical-amplificationFactor}
\end{equation}
which is the rational approximation of $e^{\lambda}$ (Eq.~\eqref{eq:ExpmRational}).
For the sake of a simple notation, the argument $\lambda$ will be
omitted hereafter.

Introducing the spectral radius describing the amplitude
\begin{equation}
\rho=|R|\label{eq:Analytical-rho}
\end{equation}
and the phase angle 
\begin{align}
\bar{\Omega} & =\arg\left(R\right)\,,\label{eq:Analytical-arg}
\end{align}
the amplification factor is expressed in the polar form as 
\begin{equation}
R=\rho e^{\mathrm{i}\bar{\Omega}}\,.\label{eq:Analytical-R-polar}
\end{equation}
Defining
\begin{equation}
\bar{\omega}=\dfrac{\bar{\Omega}}{\Delta t}\,,\label{eq:Analytical-stepping-omega}
\end{equation}
\begin{equation}
\bar{T}=\dfrac{2\pi}{\bar{\omega}}=\dfrac{2\pi\Delta t}{\bar{\Omega}}\,,\label{eq:Analytical-stepping-T}
\end{equation}
Eq.~\eqref{eq:Analytical-R-polar} is rewritten as
\begin{equation}
R=\rho e^{\mathrm{i}\bar{\omega}\Delta t}=\rho e^{\mathrm{i}2\pi\frac{\Delta t}{\bar{T}}}\,.\label{eq:Analytical-R-dt}
\end{equation}
Using Eq.~\eqref{eq:Analytical-R-dt}, the time-stepping scheme in
Eq.~\eqref{eq:Analytical-stepping-w-eq} is expressed as
\begin{equation}
\tilde{w}_{n}=\tilde{w}_{n-1}\rho e^{\mathrm{i}\bar{\omega}\Delta t}=\tilde{w}_{n-1}\rho e^{\mathrm{i}2\pi\frac{\Delta t}{\bar{T}}}\,.\label{eq:Analytical-stepping-w}
\end{equation}
It becomes apparent that $\bar{\omega}$ and $\bar{T}$ are the angular
frequency and period obtained in the time-stepping formulation. Without
loss of generality, Eq.~\eqref{eq:Analytical-stepping-w} is applied
repeatedly with a uniform size of time step $\Delta t$ from $n=0$
($t=0)$. The discrete-time solution at $t_{n}=n\Delta t$ is obtained
explicitly as 
\begin{equation}
\tilde{w}_{n}=w_{0}\rho^{n}e^{\mathrm{i}2\pi\frac{t_{n}}{\bar{T}}}\,.\label{eq:Analytical-w-sln-discrete}
\end{equation}

Comparing the period in the discrete-time solution $\bar{T}$ in Eq.~\eqref{eq:Analytical-stepping-T}
to the period in the continuous-time solution $T=2\pi/\omega$ (Eq.~\eqref{eq:Analytical-w-sln-general}),
the relative period error is obtained as
\begin{equation}
\dfrac{\bar{T}-T}{T}=\dfrac{\omega\Delta t}{\bar{\Omega}}-1\,.\label{eq:Analytical-stepping-periodError}
\end{equation}

The amount of numerical dissipation can also be measured by the damping
ratio. Rewriting Eq.~\eqref{eq:Analytical-w-sln-discrete} in the
form 
\begin{equation}
\tilde{w}_{n}=w_{0}e^{-\bar{\zeta}\bar{\omega}t_{n}}e^{\mathrm{i}2\pi\frac{t}{\bar{T}}}\,,\label{eq:Analytical-dampingRatio-def}
\end{equation}
the damping ratio is determined as
\begin{equation}
\bar{\zeta}=-\dfrac{\ln\rho}{\bar{\omega}\Delta t}=-\dfrac{\ln\rho}{\bar{\Omega}}\,.\label{eq:Analytical-dampingRatio}
\end{equation}

For the convenience on examining the effects of time step size, the
ratio of amplitude of vibration after one period $T$, i.e., $T/\Delta t$
steps, is expressed using Eq.~\eqref{eq:Analytical-w-sln-discrete}
as
\begin{equation}
\dfrac{|w(T+t)|}{|w(t)|}=\rho^{T/\Delta t}\,.\label{eq:Analytical-amplitudeRatio1}
\end{equation}
After $N_{p}$ period of vibrations, the amplitude ratio becomes
\begin{equation}
\dfrac{|w(N_{p}T+t)|}{|w(t)|}=\rho^{N_{p}T/\Delta t}\,.\label{eq:Analytical-amplitudeRatioNp}
\end{equation}

A \texttt{MATHEMATICA} code of the functions defined above to evaluate
the dissipative and dispersive characteristics of proposed time-stepping
scheme is provided in Fig.~\ref{fig:MathematicaCode} for use in
the subsequent sections. In the function \texttt{\footnotesize{}myArg},
the principal value of the phase angle is shifted to the positive
side. The other functions should be self-evident.
\begin{figure}
\texttt{\small{}f{[}n\_{]}:=Factorial{[}n{]}}{\small\par}

\texttt{\small{}padeP{[}L\_,M\_,x\_{]}:=Sum{[}f{[}(M+L-i){]}/f{[}i{]}/f{[}L-i{]}{*}x\textasciicircum i,\{i,0,L\}{]}}{\small\par}

\texttt{\small{}padeQ{[}L\_,M\_,x\_{]}:=f{[}M{]}/f{[}L{]}{*}Sum{[}f{[}(M+L-i){]}/f{[}i{]}/f{[}M-i{]}{*}(-x)\textasciicircum i,\{i,0,M\}{]}}{\small\par}

\texttt{\small{}P{[}L\_,M\_,a\_,x\_{]}:=a{*}padeP{[}M,M,x{]}+(1-a){*}padeP{[}L,M,x{]}}{\small\par}

\texttt{\small{}Q{[}L\_,M\_,a\_,x\_{]}:=a{*}padeQ{[}M,M,x{]}+(1-a){*}padeQ{[}L,M,x{]}}{\small\par}

\texttt{\small{}R{[}L\_,M\_,a\_,x\_{]}:=P{[}L,M,a,2{*}Pi{*}x{*}I{]}/Q{[}L,M,a,2{*}Pi{*}x{*}I{]}}{\small\par}

\texttt{\small{}rho{[}L\_,M\_,a\_,x\_{]}:=Abs{[}R{[}L,M,a,x{]}{]}}{\small\par}

\texttt{\small{}ph{[}L\_,M\_,a\_,x\_{]}:=R{[}L,M,a,x{]}/Abs{[}R{[}L,M,a,x{]}{]}}{\small\par}

\texttt{\small{}myArg{[}x\_,y\_{]}:=If{[}Im{[}y{]}<0,Arg{[}y{]}+2{*}Pi,If{[}x>1,Arg{[}y{]}+2{*}Pi,Arg{[}y{]}{]}{]}}{\small\par}

\texttt{\small{}pE{[}L\_,M\_,a\_,x\_{]}:=2{*}Pi{*}x/myArg{[}x,R{[}L,M,a,x{]}{]}-1}{\small\par}

\texttt{\small{}dR{[}L\_,M\_,a\_,x\_{]}:=-Log{[}rho{[}L,M,a,x{]}{]}/myArg{[}x,R{[}L,M,a,x{]}{]}}{\small\par}

\texttt{\small{}aR{[}L\_,M\_,a\_,x\_,n\_{]}:=rho{[}L,M,a,x{]}\textasciicircum (n/x)}{\small\par}

\caption{Listing of the \texttt{MATHEMATICA} code for evaluating the dissipative
and dispersive characteristics of present scheme.\label{fig:MathematicaCode}}
\end{figure}

The \texttt{MATHEMATICA} code in Fig.~\ref{fig:MathematicaCode} can be used
to compare the dissipative and dispersive characteristics of the proposed
scheme with those of other time-integration methods reported in the
literature. In this paper, only the comparison with the HHT-$\alpha$
method \citep{Hilber1977} is reported. The HHT-$\alpha$ method is
currently the most widely available scheme for time integration in
commercial finite element software packages. Similar to the proposed
method, the HHT-$\alpha$ method can be used with only two user-specified
parameters, i.e., the time step size and a parameter $\alpha$ controlling
the amount of numerical dissipation. The parameter $\alpha$ is related
to the spectral radius at high-frequency limit $\rho_{\infty}$ by
\citep{Chung1993}
\begin{equation}
\rho_{\infty}=\dfrac{1+\alpha}{1-\alpha}\,.\label{eq:alpha2RhoInfty}
\end{equation}
The parameters $\beta=(1-\alpha)^{2}/4$ and $\gamma=1/2-\alpha$
are chosen to maintain second-order accuracy \citep{Hilber1977}.
The \texttt{MATHEMATICA} codes to evaluate the dissipative and dispersive
characteristics of the HHT-$\alpha$ scheme are listed in Fig.~\ref{fig:MathematicaCode-HHT}.
\begin{figure}
\texttt{\small{}HHT{[}a\_,x\_{]}:=Module{[}\{b,g,o,d,a1,a2,a3,e,r\},}{\small\par}

\texttt{\small{}b = (1-a)\textasciicircum 2/4;}{\small\par}

\texttt{\small{}g = 1/2 - a;}{\small\par}

\texttt{\small{}o= 2{*}Pi{*}x;}{\small\par}

\texttt{\small{}d = 1 + (1+a){*}b{*}o\textasciicircum 2;}{\small\par}

\texttt{\small{}a1 = 1 - o\textasciicircum 2{*}((1+a){*}(g+1/2)-a{*}b)/(2{*}d);}{\small\par}

\texttt{\small{}a2 = 1- o\textasciicircum 2{*}(g-1/2 + 2{*}a{*}(g-b))/d;}{\small\par}

\texttt{\small{}a3 = a{*}o\textasciicircum 2{*}(b-g+0.5)/d;}{\small\par}

\texttt{\small{}e = t\textasciicircum 3-2{*}a1{*}t\textasciicircum 2+a2{*}t-a3;}{\small\par}

\texttt{\small{}r=Sort{[} t/.NSolve{[}e==0,t{]}, Abs{[}\#1{]}>Abs{[}\#2{]}\&{]};}{\small\par}

\texttt{\small{}Complex{[}Re{[} r{[}{[}1{]}{]}{]}, Abs{[}Im{[}r{[}{[}1{]}{]}{]}{]}
{]}}{\small\par}

\texttt{\small{}{]}}{\small\par}

\texttt{\small{}HHTrho{[}a\_,x\_{]}:= Abs{[}HHT{[}a,x{]}{]}}{\small\par}

\texttt{\small{}HHTph{[}a\_,x\_{]}:=HHT{[}a,x{]}/Abs{[}HHT{[}a,x{]}{]}}{\small\par}

\texttt{\small{}HHTmyArg{[}x\_,y\_{]}:=If{[}Im{[}y{]}<0,Arg{[}y{]}+2{*}Pi,If{[}x>1,NaN,Arg{[}y{]}{]}{]}}{\small\par}

\texttt{\small{}HHTpE{[}a\_,x\_{]}:= 2{*}Pi{*}x/HHTmyArg{[}x,HHT{[}a,x{]}{]}-1}{\small\par}

\texttt{\small{}HHTdR{[}a\_,x\_{]}:=-Log{[}HHTrho{[}a,x{]}{]}/HHTmyArg{[}x,HHT{[}a,x{]}{]}}{\small\par}

\texttt{\small{}HHTaR{[}a\_,x\_,n\_{]}:=HHTrho{[}a,x{]}\textasciicircum (n{*}x)}{\small\par}

\caption{Listing of \texttt{MATHEMATICA} code for evaluating the dissipative
and dispersive characteristics of HHT$-\alpha$ scheme. \label{fig:MathematicaCode-HHT}}
\end{figure}

\subsection{Spectral radius\label{subsec:Spectral-radius}}

The dissipative characteristics of the proposed time-stepping scheme
based on Padé expansions are examined by evaluating the spectral radius.
The three cases of (i) diagonal Padé expansions of orders ($M,M$),
(ii) sub-diagonal expansions of orders ($M-1,M$), and (iii) Padé
expansions mixing the orders of ($M,M$) and ($M-1,M$) are considered.

The diagonal Padé expansions are given in Eq.~\eqref{eq:Pade-mixed-PQ}
with $\rho_{\infty}=1$. Using Eq.~\eqref{eq:PQ_L/M}, they are written
as \begin{subequations}\label{eq:Pade-diag-PQ}
\begin{align}
\mathbf{P} & =\mathbf{P}_{M/M}(\mathbf{A})\,,\label{eq:Pade-diag-P}\\
\mathbf{Q} & =\mathbf{Q}_{M/M}(\mathbf{A})=\mathbf{P}_{M/M}(-\mathbf{A})\,.\label{eq:Pade-diag-Q}
\end{align}
\end{subequations}Using Eqs.~\eqref{eq:Pade-diag-PQ}, \eqref{eq:Analytical-amplificationFactor}
and \eqref{eq:Analytical-lambda-T}, the spectral radius given in
Eq.~\eqref{eq:Analytical-rho} is written as
\begin{equation}
\rho=|R(\mathrm{i}2\pi\Delta t/T)|=\left|\dfrac{P_{M/M}(\mathrm{i}2\pi\Delta t/T)}{P_{M/M}(-\mathrm{i}2\pi\Delta t/T)}\right|=1\label{eq:Pade-diag-rho}
\end{equation}
at any size of time step $\Delta t$, since the polynomials in the
numerator and denominator are complex conjugates.

\begin{figure}
\begin{centering}
\includegraphics[width=0.5\textwidth]{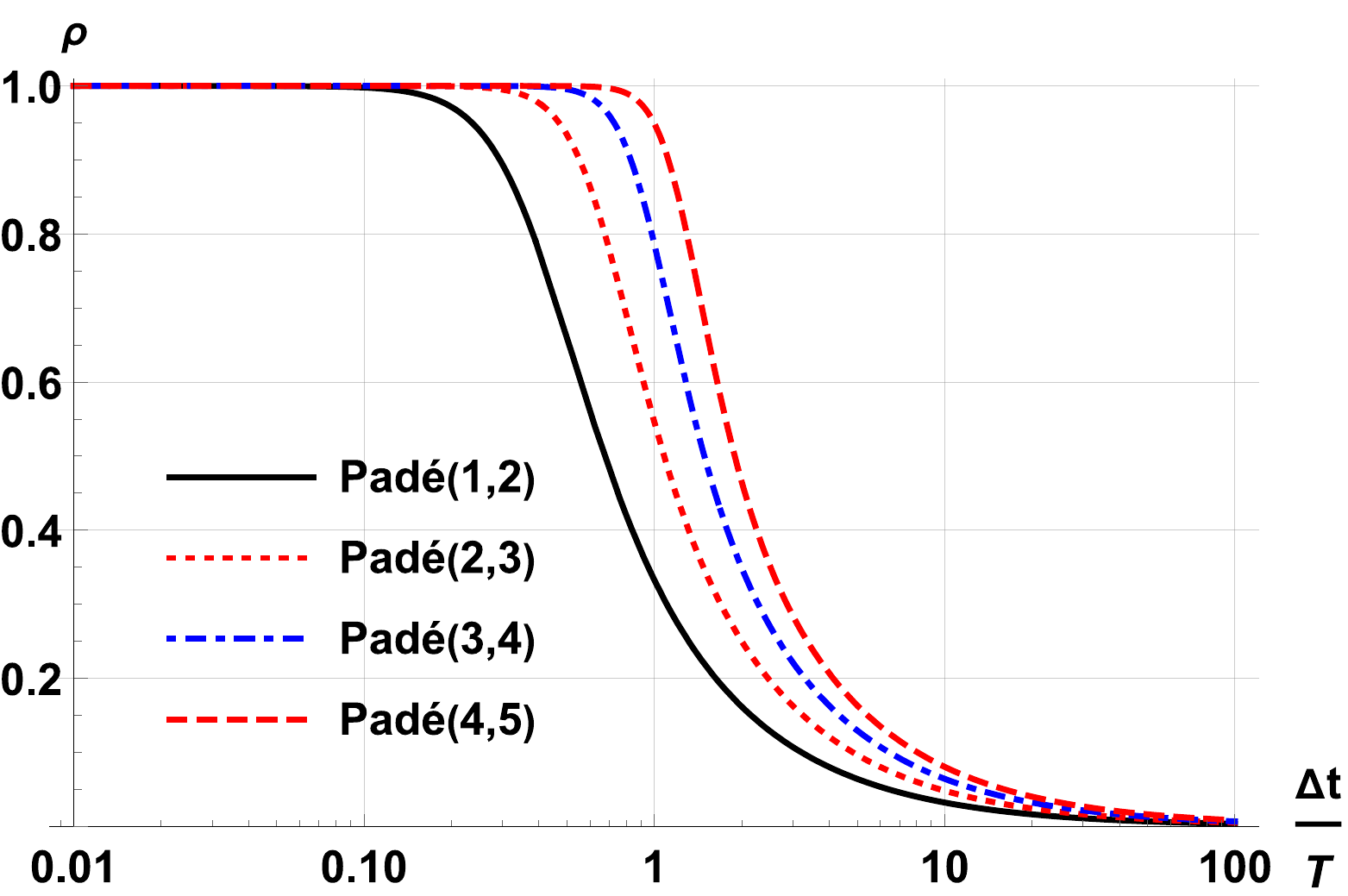}
\par\end{centering}
\medskip{}

\texttt{\footnotesize{}a= 0; }{\footnotesize\par}

\texttt{\footnotesize{}LogLinearPlot{[}\{rho{[}1, 2, a, x{]}, rho{[}2,
3, a, x{]}, rho{[}3, 4, a, x{]}, rho{[}4, 5, a, x{]}\}, \{x, 0.01,
100\}, PlotRange -> \{0, 1.01\}{]}}{\footnotesize\par}

\caption{Spectral radii of sub-diagonal Padé expansions ($M-1,M$).\label{fig:Spectral-radii-M-1}}
\end{figure}

The sub-diagonal Padé expansions of order ($M-1,M$) are given in
Eq.~\eqref{eq:Pade-mixed-PQ} with $\rho_{\infty}=0$. The spectral
radius is expressed as
\begin{equation}
\rho=|R(\mathrm{i}2\pi\Delta t/T)|=\left|\dfrac{P_{(M-1)/M}(\mathrm{i}2\pi\Delta t/T)}{Q_{(M-1)/M}(\mathrm{i}2\pi\Delta t/T)}\right|\,,\label{eq:Pade-sbudiag-rho}
\end{equation}
where the polynomials $P_{(M-1)/M}$ and $Q_{(M-1)/M}$ are given
in Eq.~\eqref{eq:PQ_L/M} with $L=M-1$. The highest orders of $P_{(M-1)/M}$
and $Q_{(M-1)/M}$ are equal to $M-1$ and $M$, respectively. The
limit of the spectral radius as $\Delta t$ increases is obtained
as 
\begin{equation}
\rho_{\infty}=\lim_{\Delta t\rightarrow\infty}\rho=\lim_{\Delta t\rightarrow\infty}\left|\dfrac{(\mathrm{i}2\pi\Delta t/T)^{M-1}}{(\mathrm{i}2\pi\Delta t/T)^{M}}\right|=0\,.\label{eq:Pade-sbudiag-rho-infty}
\end{equation}
The spectral radii of the sub-diagonal Padé expansions ($M-1,M$)
are plotted in Fig.~\ref{fig:Spectral-radii-M-1} as functions of
$\Delta t/T$ for $M=2$, $3$, $4$ and $5$. The \texttt{MATHEMATICA}
command for plotting the curves (without the formatting options) is
provided under the plot. It calls the function \textquotedbl\texttt{rho}\textquotedbl{}
defined in Fig.~\ref{fig:MathematicaCode}. High-order accuracy at
small size of time step $\Delta t/T$ is observed. The curves descend
rapidly to approach zero as the size of time step increases.

\begin{figure}
\begin{centering}
\includegraphics[width=0.5\textwidth]{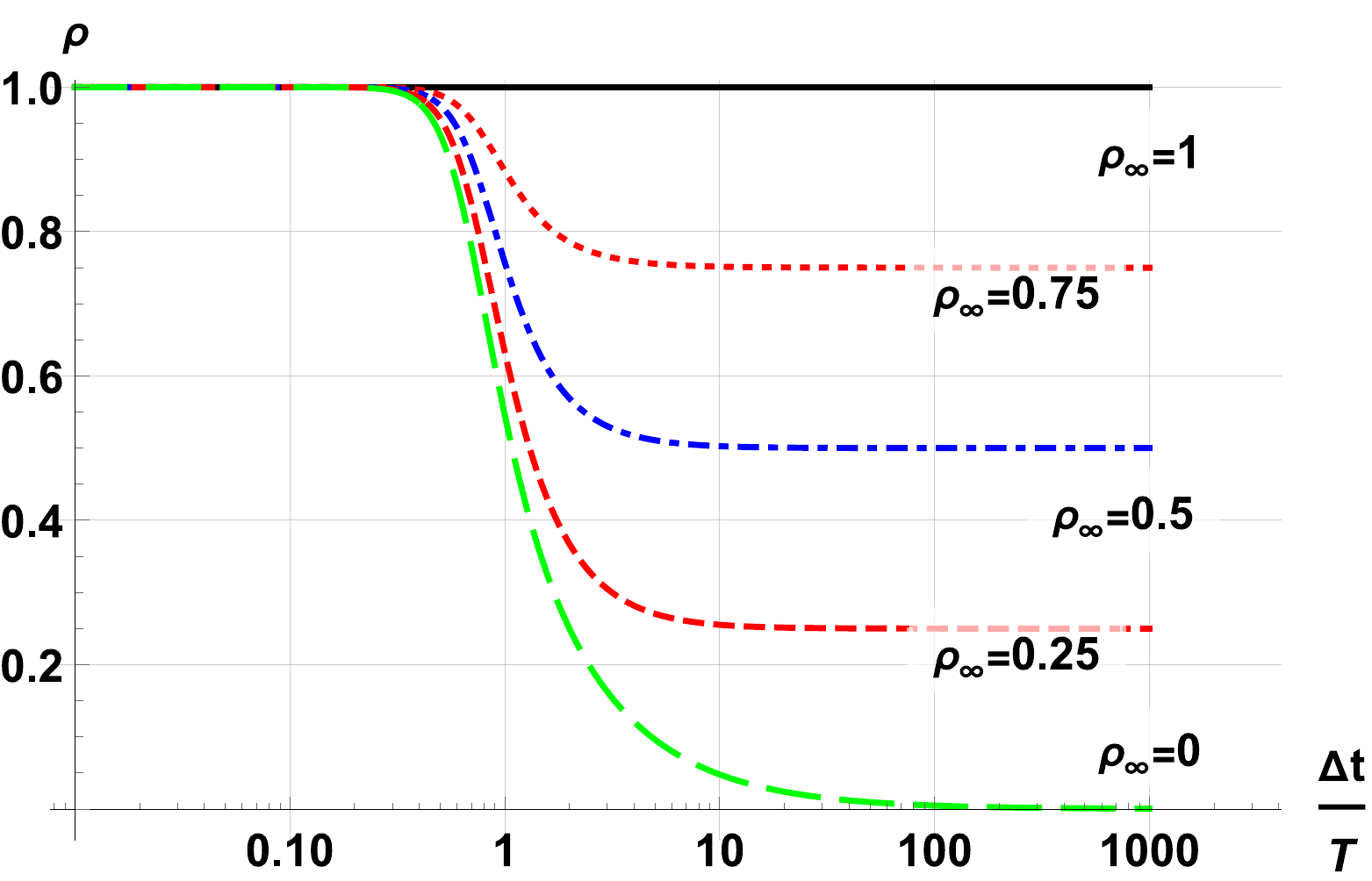}
\par\end{centering}
\medskip{}

\texttt{\small{}M = 3; L = M - 1; }{\small\par}

\texttt{\small{}LogLinearPlot{[}\{rho{[}L, M, 1, x{]}, rho{[}L, M,
3/4, x{]}, rho{[}L, M, 1/2, x{]}, rho{[}L, M, 1/4, x{]}, rho{[}L,
M, 0, x{]}\}, \{x, 0.01, 1000\}, PlotRange -> \{0, 1.01\}{]}}{\small\par}

\caption{Spectral radii of mixed-order Padé expansions with $L=2$ and $M=3$
at $\rho_{\infty}=1$, $0.75$, $0.5$, $0.25,$ and $0$.\label{fig:Spectral-radii-mixed-M=00003D3}}
\end{figure}

The mixed-order Padé expansions are obtained in Eq.~\eqref{eq:Pade-mixed-PQ}
with $0<\rho_{\infty}<1$ and $L<M$. It is identified from Eq.~\eqref{eq:PQ_L/M}
that the highest orders of $P_{L/M}$ and $P_{M/M}$ are equal to
$L$ and $M$, respectively, and the highest orders of $Q_{L/M}$
and $Q_{M/M}$ are both equal to $M$. The limit of the spectral radius
as $\Delta t$ increases is obtained as 
\begin{equation}
\lim_{\Delta t\rightarrow\infty}\left|\frac{\rho_{\infty}(\mathrm{i}2\pi\Delta t/T)^{M}+(1-\rho_{\infty})(\mathrm{i}2\pi\Delta t/T)^{L}}{\rho_{\infty}(\mathrm{i}2\pi\Delta t/T)^{M}+(1-\rho_{\infty})(\mathrm{i}2\pi\Delta t/T)^{M}}\right|=\rho_{\infty}\,.\label{eq:Pade-sbudiag-rho-infty-1}
\end{equation}
It indeed tends to the specified $\rho_{\infty}$. As an example,
the spectral radii of Padé expansions mixing the orders $(2,3)$ and
$(3,3)$ are plotted in Fig.~\ref{fig:Spectral-radii-mixed-M=00003D3}
at $\rho_{\infty}=1$, $0.75$, $0.5$, $0.25,$ and $0$. It is observed
that the numerical dissipation is controlled by the specified value
of $\rho_{\infty}$ over the range between $0$ and $1$. The spectral
radius is strictly less than or equal to $1$. 

\begin{figure}
\begin{centering}
a) \includegraphics[width=0.45\textwidth]{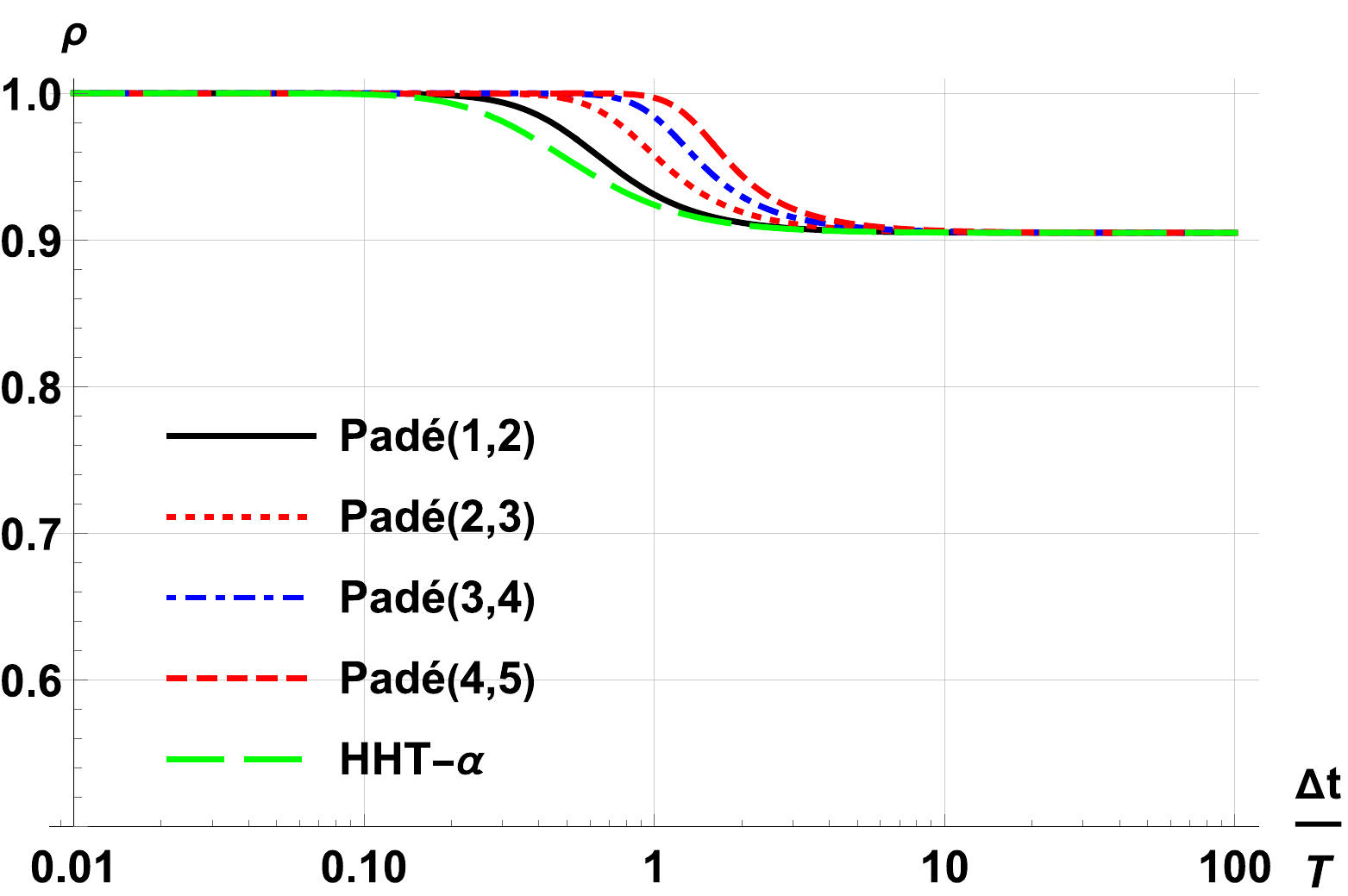}\hfill{}b)
\includegraphics[width=0.45\textwidth]{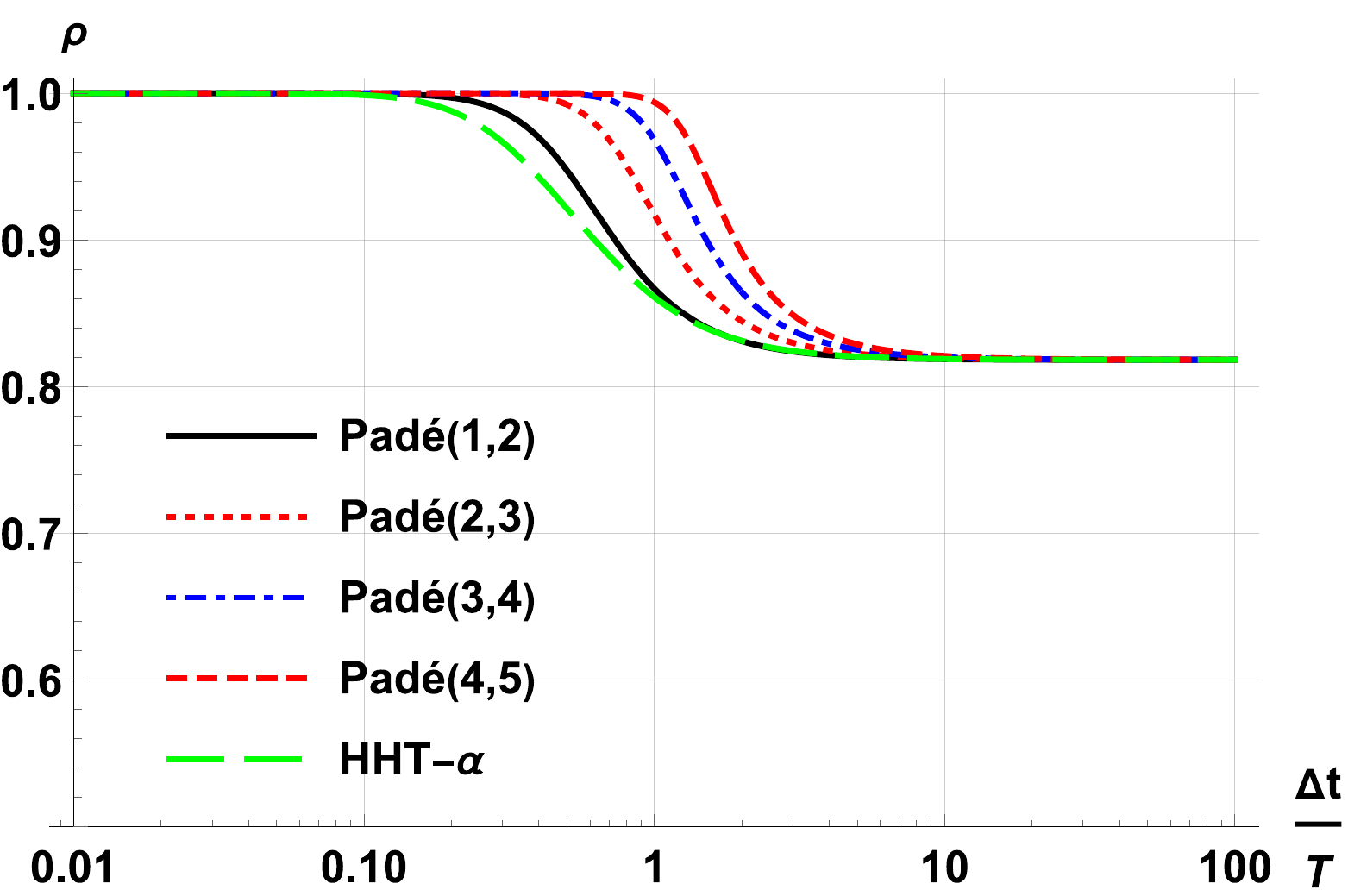}
\par\end{centering}
\begin{centering}
\hfill{}c) \includegraphics[width=0.45\textwidth]{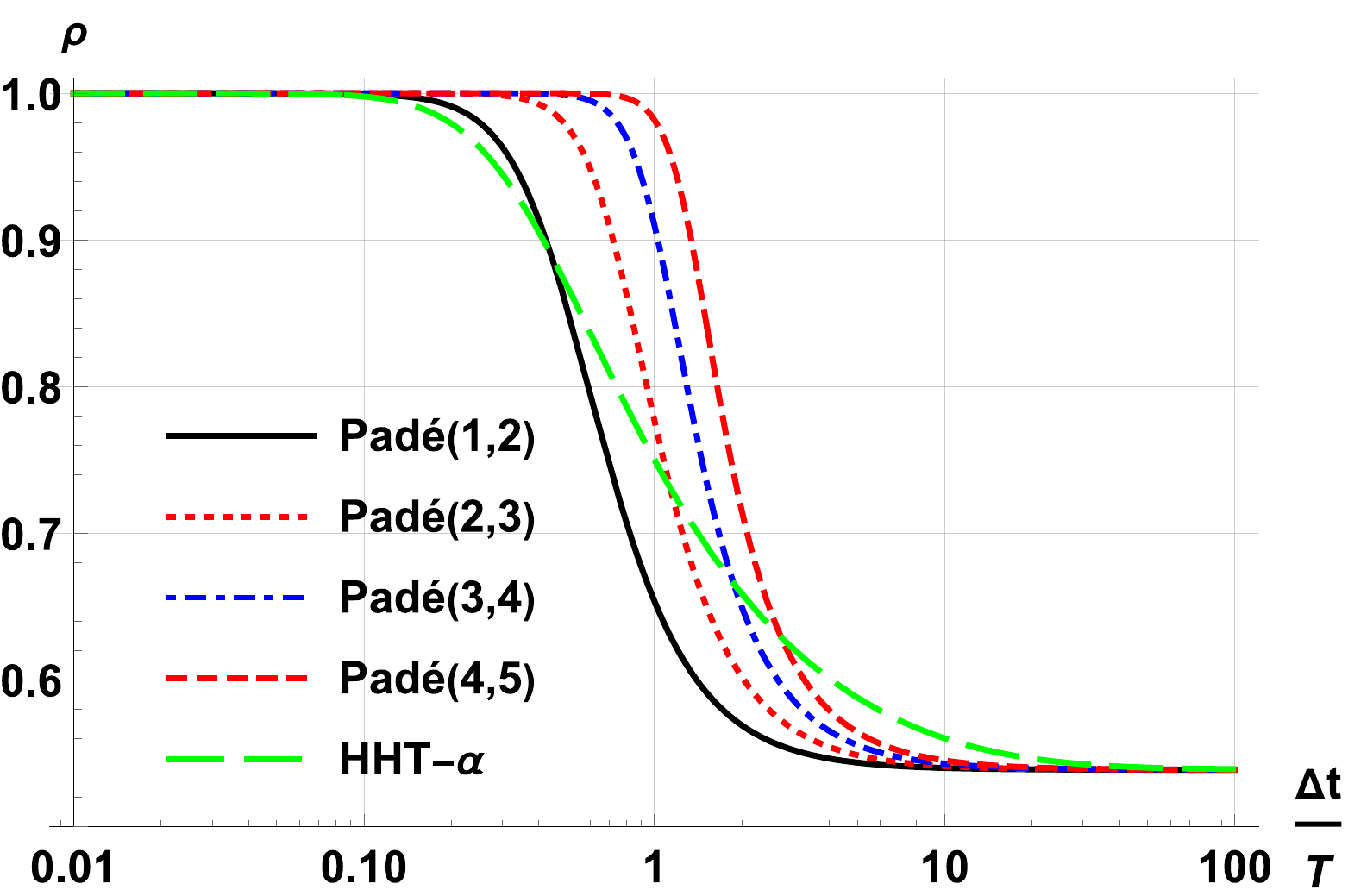}\hfill{}
\par\end{centering}
\medskip{}

\texttt{\small{}alpha = -0.05; a = (1 + alpha)/(1 - alpha);}{\small\par}

\texttt{\small{}LogLinearPlot{[}\{rho{[}1, 2, a, x{]}, rho{[}2, 3,
a, x{]}, rho{[}3, 4, a, x{]}, rho{[}4, 5, a, x{]}, HHTrho{[}alpha,
x{]}\}, \{x, 0.01, 100\}, PlotRange -> \{0.5, 1.01\}{]}}{\small\par}

\caption{Spectral radii of mixed-order Padé expansions and HHT-$\alpha$ scheme
for a) $\rho_{\infty}=0.90476$ at $\alpha=-0.05$, b) $\rho_{\infty}=0.81818$
at $\alpha=-0.1$, and c) $\rho_{\infty}=0.53846$ at $\alpha=-0.3$.\label{fig:Spectral-radii-HHT}}
\end{figure}

The numerical dissipative characteristics of the proposed time-stepping
schemes are compared with those of the HHT-$\alpha$ scheme for the
cases of $\alpha=-0.05$, $\alpha=-0.1$, and $\alpha=-0.3$. In the
high-frequency limit $\Delta t/T\rightarrow\infty$, the spectral
radius tends to $\rho_{\infty}=0.90476$ at $\alpha=-0.05$, $\rho_{\infty}=0.81818$
at $\alpha=-0.1,$ and $\rho_{\infty}=0.53846$ at $\alpha=-0.3$.
For each given value of $\alpha$, the corresponding value of $\rho_{\infty}$
is taken as an input in the proposed scheme. The spectral radii are
compared in Fig.~\ref{fig:Spectral-radii-HHT}a for $\alpha=-0.05$,
Fig.~\ref{fig:Spectral-radii-HHT}b for $\alpha=-0.1$, and Fig.~\ref{fig:Spectral-radii-HHT}c
for $\alpha=-0.3$, respectively. The \texttt{MATHEMATICA} code is
given with $\alpha=-0.05$, but it can be used for other values of
$\rho_{\infty}$ (represented by the variable \texttt{a}) and $\alpha$.
It is observed that the spectral radii of the proposed high-order
scheme are closer to $1$ than the HHT-$\alpha$ scheme when the value
of $\Delta t/T$ is small and decrease faster when the value of $\Delta t/T$
becomes large. This behaviour indicates that the proposed scheme is
more accurate for lower modes and exhibits stronger numerical dissipation
for higher modes.

\subsection{Period error\label{subsec:Period-error}}

\begin{figure}
\begin{centering}
\includegraphics[width=0.45\textwidth]{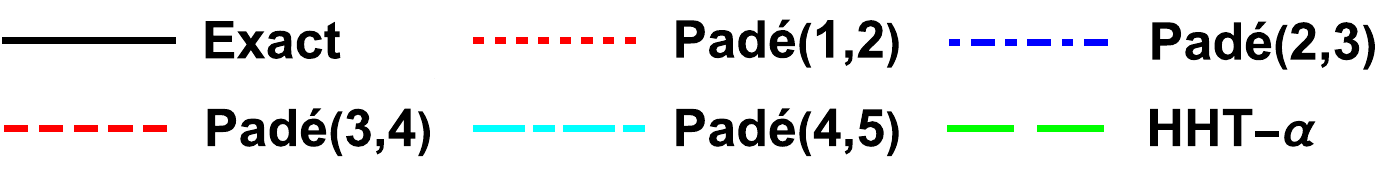}
\par\end{centering}
\medskip{}

\begin{centering}
a) \includegraphics[width=0.45\textwidth]{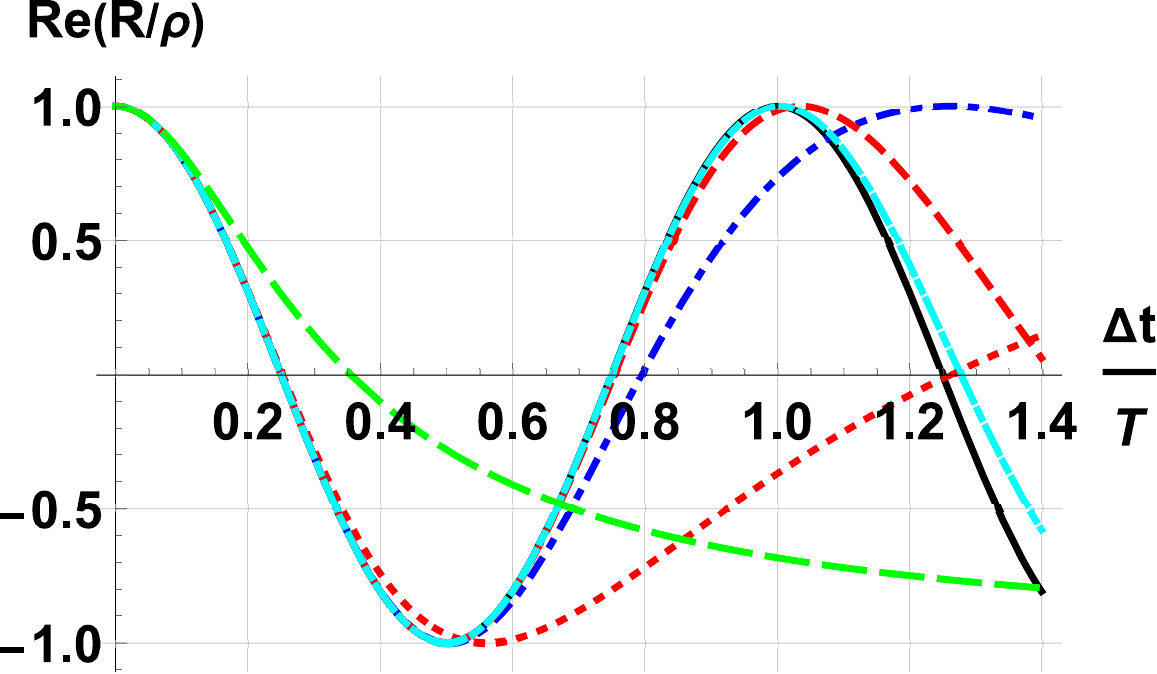}\hfill{}b)
\includegraphics[width=0.45\textwidth]{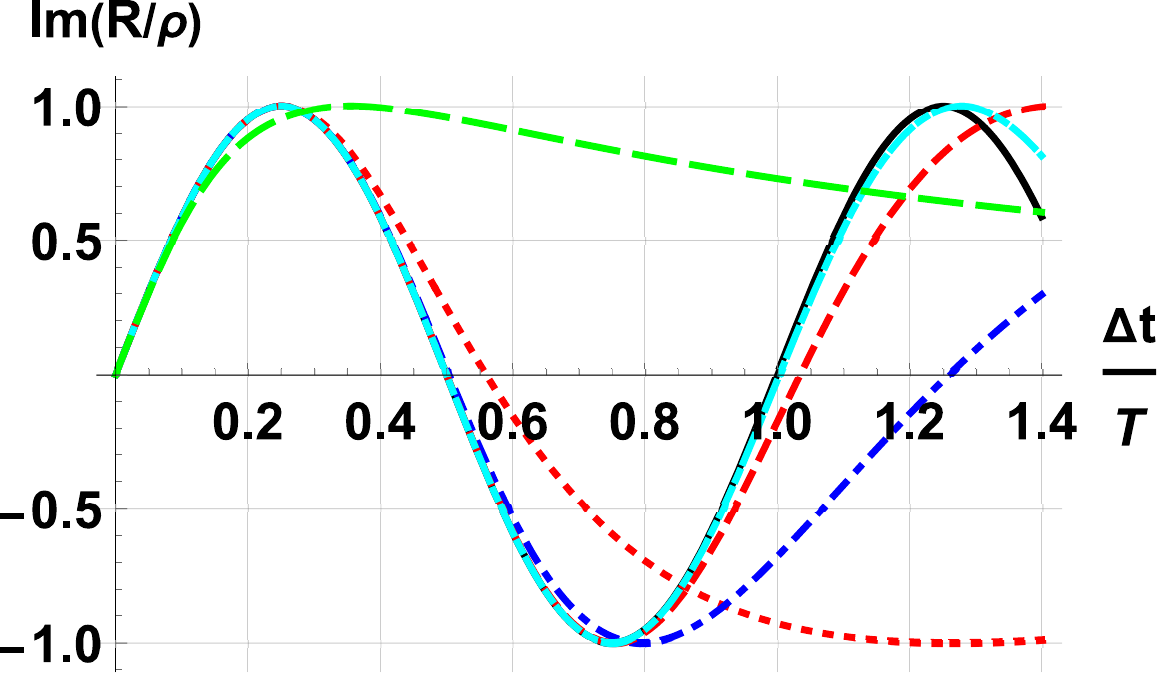}
\par\end{centering}
\medskip{}

\texttt{\small{}alpha = -0.3; a = (1 + alpha)/(1 - alpha); }{\small\par}

\texttt{\small{}Plot{[}\{Re{[}Exp{[}2{*}Pi{*}x{*}I{]}{]}, Re{[}ph{[}1,
2, a, x{]}{]}, Re{[}ph{[}2, 3, a, x{]}{]}, Re{[}ph{[}3, 4, a, x{]}{]},
Re{[}ph{[}4, 5, a, x{]}{]}, Re{[}HHTph{[}alpha, x{]}{]}\}, \{x, 0.0,
1.4\}{]}}{\small\par}

\caption{Phase of vibration of mixed-order Padé expansions at $\rho_{\infty}=0.53846$
and HHT-$\alpha$ with $\alpha=-0.3$: a) Real part, b) Imaginary
part.\label{fig:phase}}
\end{figure}

As it is observed by comparing the discrete-time solution in Eq.~\eqref{eq:Analytical-stepping-w-eq}
with Eq.~\eqref{eq:Analytical-R-polar} to the continuous-time solution
in Eq.~\eqref{eq:Analytical-stepping-exact}, numerical errors arise
not only from the amplitude but also from the phase of vibration.
The exact solution of the phase $e^{\mathrm{i}2\pi\Delta t/T}$ is
approximated by $e^{\mathrm{i}\bar{\Omega}}=R/\rho$ in the time-stepping
scheme. The approximation~$R/\rho$ of the proposed scheme with $\rho_{\infty}=0.53846$
and the HHT-$\alpha$ scheme with $\alpha=-0.3$ are compared to the
exact solution in Fig.~\ref{fig:phase}. It shows that the proposed
high-order scheme is significantly more accurate than the HHT-$\alpha$
scheme. As the order increases, the accuracy of approximation improves
rapidly.

\begin{figure}
\begin{centering}
\includegraphics[width=0.45\textwidth]{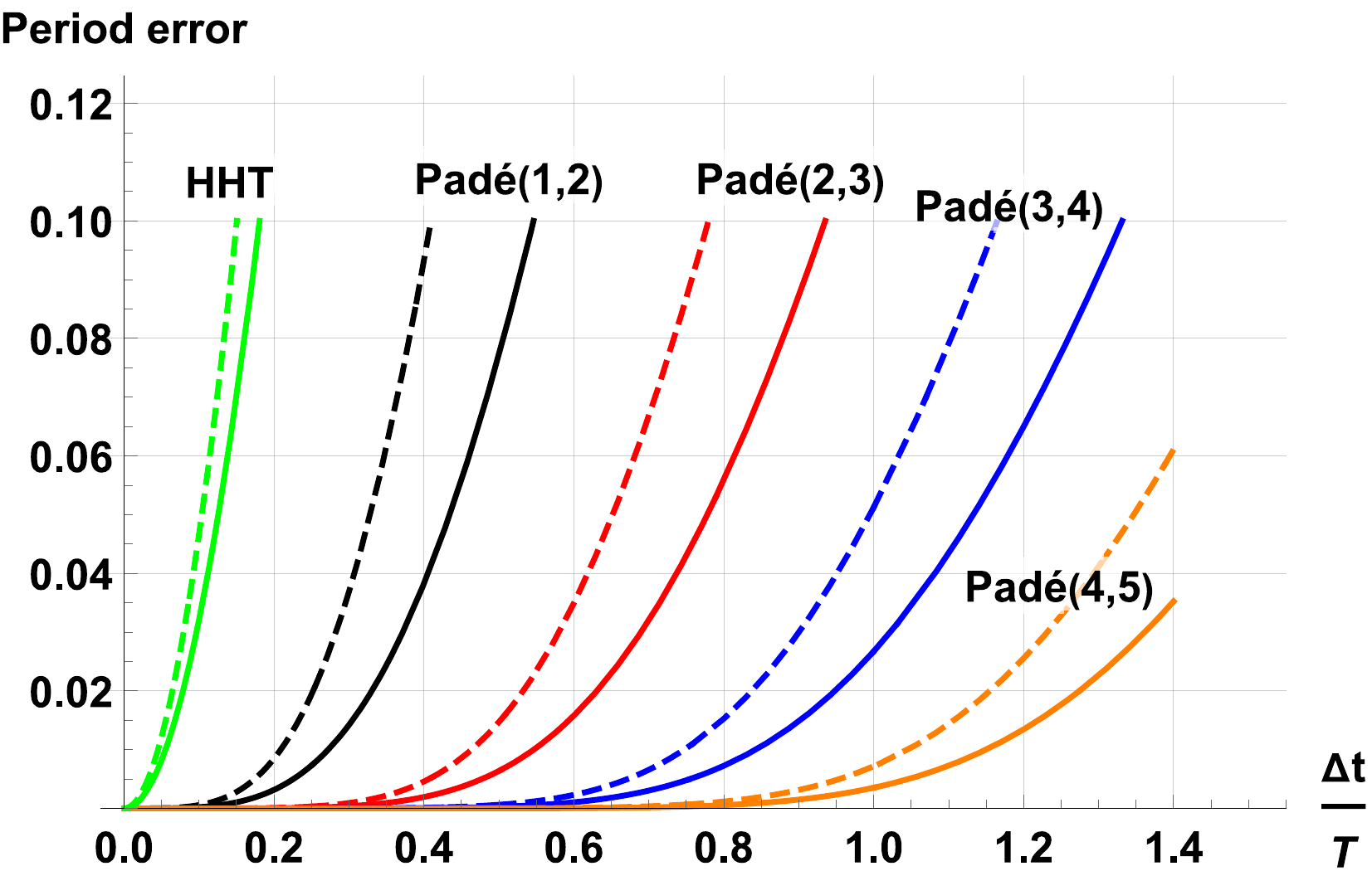}
\par\end{centering}
\medskip{}
\texttt{\small{}a1 = 1; a2 = 0; Plot{[}\{pE{[}1, 2, a1, x{]}, pE{[}1,
2, a2, x{]}, pE{[}2, 3, a1, x{]}, pE{[}2, 3, a2, x{]}, pE{[}3, 4,
a1, x{]}, pE{[}3, 4, a2, x{]}, pE{[}4, 5, a1, x{]}, pE{[}4, 5, a2,
x{]}, HHTpE{[}0, x{]}, HHTpE{[}-0.3, x{]}\}, \{x, 0.001, 1.4\}, PlotRange
-> \{0, 0.1\}{]}}{\small\par}

\caption{Relative period errors. The present high-order scheme with $\rho_{\infty}=1$
and the HHT-$\alpha$ scheme with $\alpha=0$ are indicated by the
solid lines and $\rho_{\infty}=0.53846$ and $\alpha=-0.3$ by the
dashed lines.\label{fig:Relative-period-errors}}
\end{figure}

The relative period errors given by Eq.~\eqref{eq:Analytical-stepping-periodError}
for orders $(1,2)$, $(2,3)$, $(3,4)$, and $(4,5)$ of the mixed
Padé scheme are plotted in Fig.~\ref{fig:Relative-period-errors}
at two levels of numerical dissipation, together with those of the
HHT-$\alpha$ scheme. The two levels of numerical dissipation are
specified by $\rho_{\infty}=1$ and $\rho_{\infty}=0.53846$ corresponding
to $\alpha=0$ and $\alpha=-0.3$, respectively. The relative period
errors at $\rho_{\infty}=1$ are depicted by solid lines and those
at $\rho_{\infty}=0.53846$ by dashed lines. It is observed that the
relative period error decreases rapidly with the increase of the order
of the scheme. As expected, introducing numerical dissipation leads
to larger relative period errors since the order of accuracy of the
Padé expansion is reduced by one. 

\subsection{Damping ratio}

\begin{figure}
\begin{centering}
a) \includegraphics[width=0.45\textwidth]{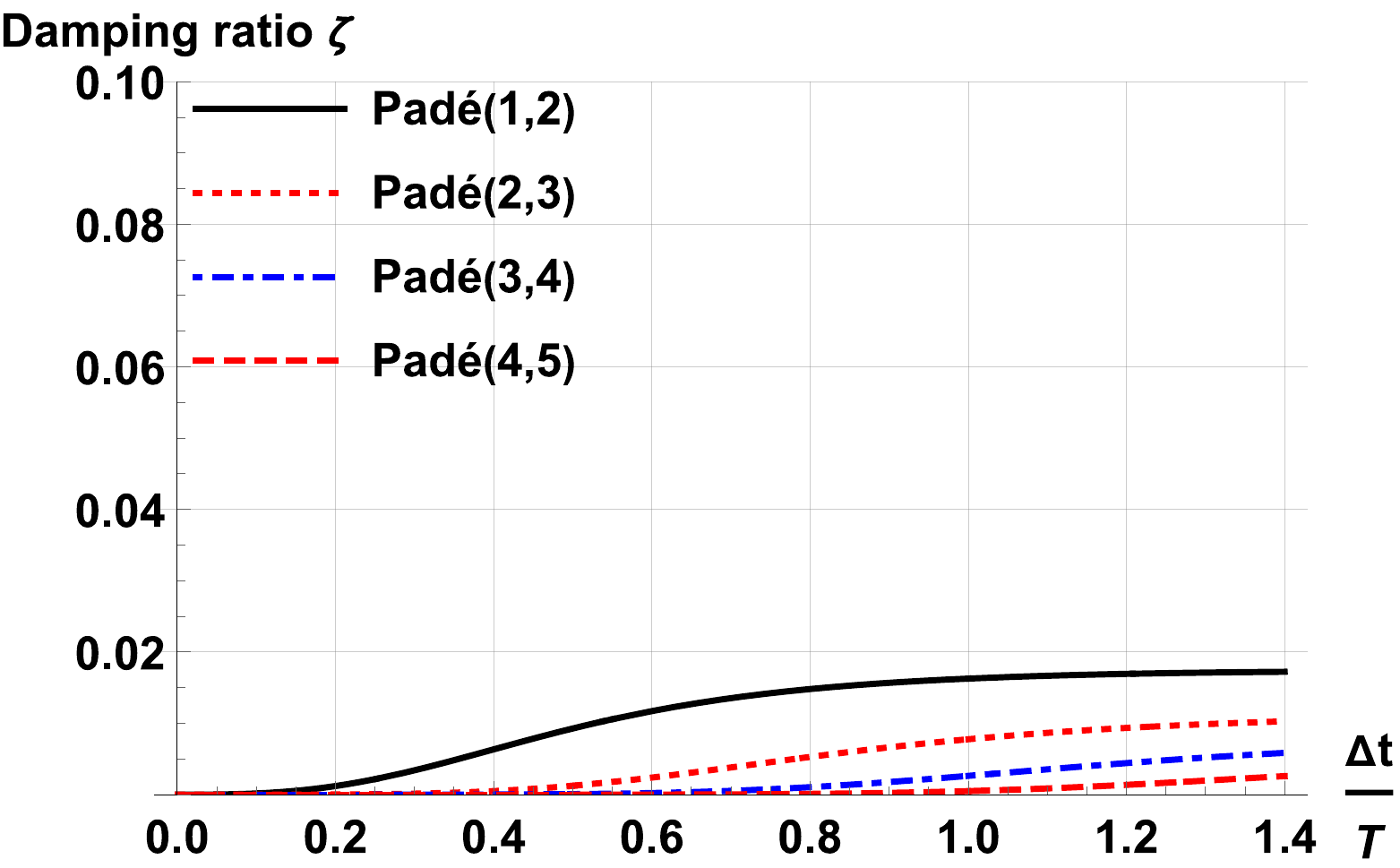}\hfill{}b)
\includegraphics[width=0.45\textwidth]{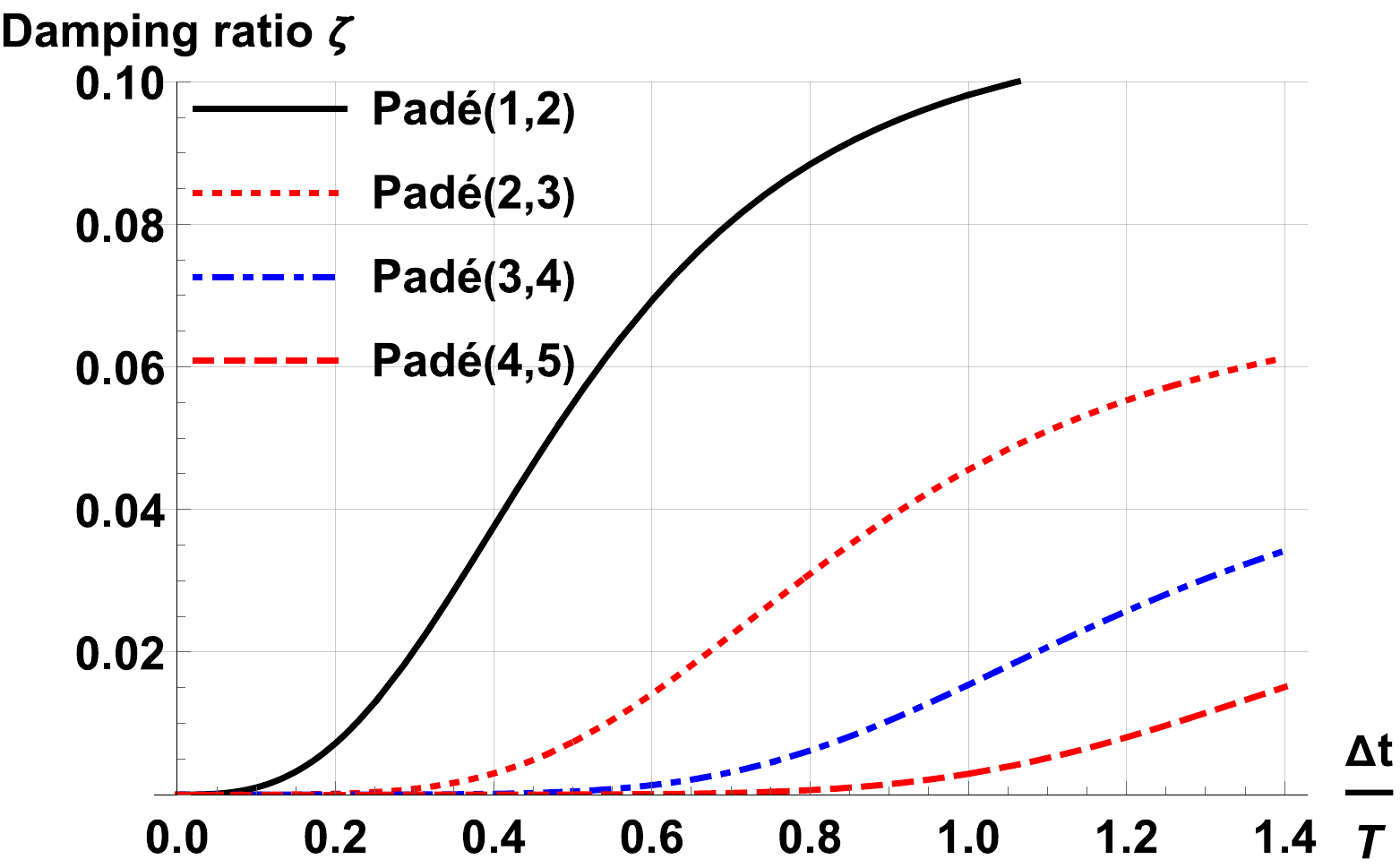}
\par\end{centering}
\medskip{}

\texttt{\small{}a = 0.90476; }{\small\par}

\texttt{\small{}Plot{[}\{dR{[}1, 2, a, x{]}, dR{[}2, 3, a, x{]}, dR{[}3,
4, a, x{]}, dR{[}4, 5, a, x{]}\}, \{x, 0.001, 1.4\}, PlotRange ->
\{0, 0.1\} {]}}{\small\par}

\caption{Damping ratios of mixed-order Padé expansions at: a) $\rho_{\infty}=0.90476$,
b) $\rho_{\infty}=0.53846$. \label{fig:Damping-ratios-of-mixed}}
\end{figure}

The damping ratios using the mixed-order Padé expansions are plotted
in Fig.~\ref{fig:Damping-ratios-of-mixed} as functions of $\Delta t/T$.
The numerical damping is specified with $\rho_{\infty}=0.90476$ in
Fig.~\ref{fig:Damping-ratios-of-mixed}a and $\rho_{\infty}=0.53846$
in Fig.~\ref{fig:Damping-ratios-of-mixed}b. It is found that the
damping ratio tends to zero when the size of time step $\Delta t/T$
is small. This ensures that the effect of numerical damping on low-frequency
modes is small. 

\subsection{Effect of the time step size on low-frequency modes \label{subsec:Discussion-on-selection}}

The choice of the time step size is critical in the effective use
of numerical dissipation. To reduce the temporal discretization error,
a small time step size is desirable. On the other hand, when the time
step size is too small, vibrations of spurious high-frequency modes
may not be sufficiently damped and pollute the solution. Therefore,
it is necessary to choose a suitable time step size that will lead
to the desired accuracy for the lower modes (below the maximum frequency
of interest), and suppress the responses of the higher modes (above
the maximum frequency of interest) at the same time. The effect of
the time step size on the vibrations of lower modes is examined in
this section in order to provide a guideline on the selection of the
time step size.

The HHT-$\alpha$ scheme is one of the most widely used time integration
methods in commercial software and in practice. It is addressed in
Section~\ref{subsec:HHT-scheme-timeStep} to provide a reference
case for the discussion on the proposed high-order scheme in Section~\ref{subsec:High-order-timeStep}.

\subsubsection{HHT-$\alpha$ scheme\label{subsec:HHT-scheme-timeStep}}

\begin{figure}
\begin{centering}
\includegraphics[width=0.45\textwidth]{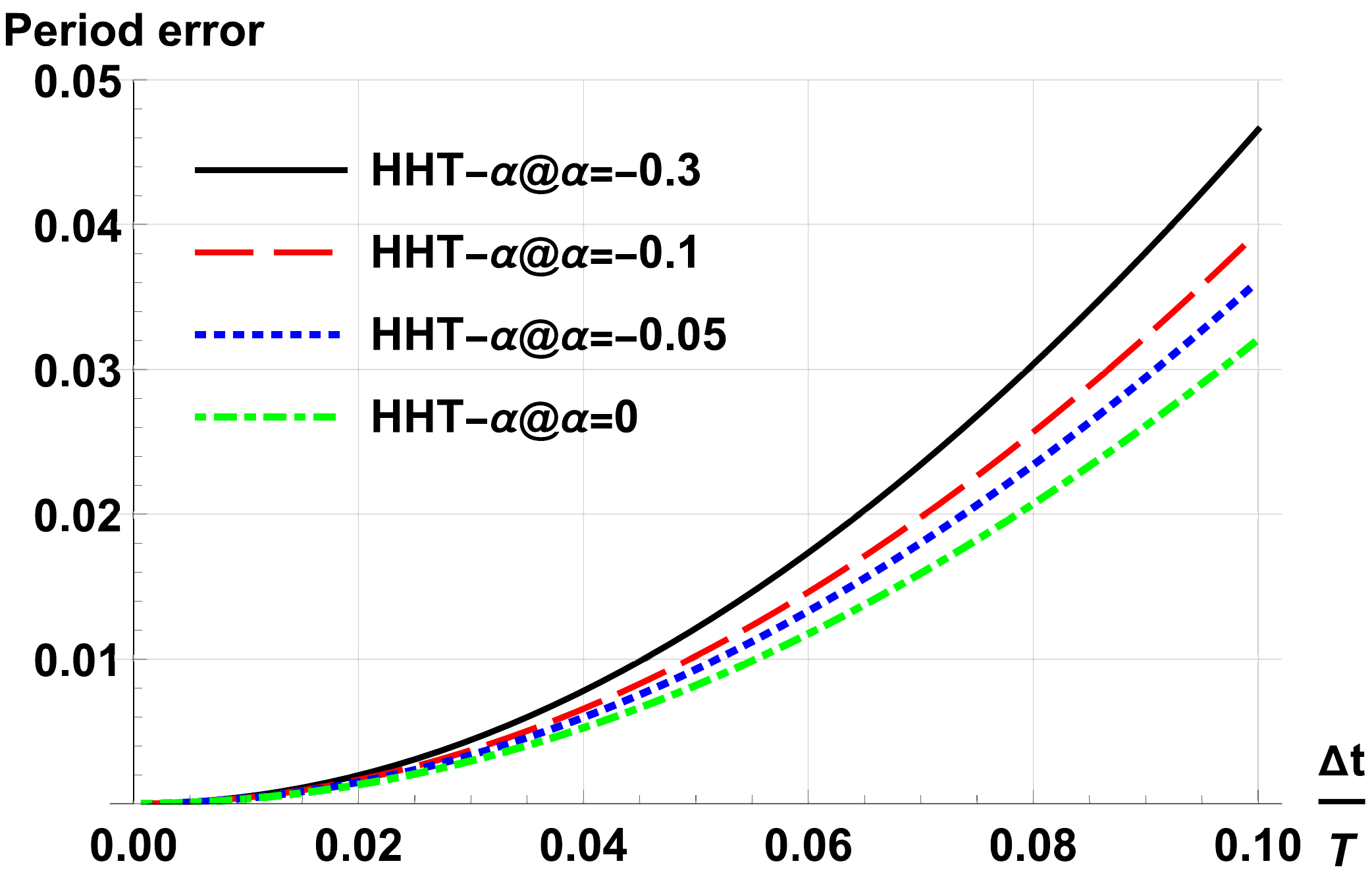}
\par\end{centering}
\medskip{}

\texttt{\small{}Plot{[}\{HHTpE{[}0, x{]}, HHTpE{[}-0.05, x{]}, HHTpE{[}-0.1,
x{]}, HHTpE{[}-0.3, x{]}\}, \{x, 0.001, 0.1\}, PlotRange -> \{0, 0.05\}{]}}{\small\par}

\caption{Relative period error of HHT-$\alpha$ method with the parameter $\alpha=0$,
$-0.05$, $-0.1$ and $-0.3$. \label{fig:HHT-alpha-method}}
\end{figure}

The HHT-$\alpha$ scheme is of second-order accuracy. The amount of
numerical dissipation is controlled by selecting the parameter $\alpha$
in the range of $0\geq\alpha\geq-1/3$. The period errors of $\alpha=0$,
$-0.05$, $-0.1$, and $-0.3$ are shown in Fig.~\ref{fig:HHT-alpha-method}.
 It is observed that numerical dissipation increases the relative
period error. In the following, the case of $\alpha=-0.3$ is considered.
In order to limit the error to $3\%$ at a given period \textbf{$T$},
the step size should not be larger than $\Delta t=0.08T$, i.e., about
$12$ steps per period, which is a common choice in practice. To reduce
the relative period error to $1\%$, the time increment must be below
$\Delta t=0.04T$, which corresponds to 25 steps per period.

It is worthwhile to note that the error in the overall time integration
also depends on the frequency contents of the excitation and will
be smaller than the error at the maximum frequency of interest.

\subsubsection{High-order scheme\label{subsec:High-order-timeStep}}

It is shown in Section~\ref{subsec:Spectral-radius} that the proposed
high-order scheme exhibits better properties of numerical dissipation
than the HHT-$\alpha$ method in the range of $-0.05\geq\alpha\geq-0.3$
(i.e., $0.90476\geq\rho_{\infty}\geq0.53846$). Parametric studies
on the numerical examples in Section~\ref{sec:Numerical-examples}
demonstrate that spurious high-frequency oscillations can be effectively
dissipated when the user-specified parameter is chosen between $0.90476\geq\rho_{\infty}\geq0$.
Only the case of $\rho_{\infty}=0.53846$, corresponding to $\alpha=-0.3$
in the HHT-$\alpha$ scheme, is considered in this section. When the
value of $\rho_{\infty}$ is smaller, the effect of numerical dissipation
on lower modes will be smaller at the same size of time step. Other
cases such as the $L$-stable scheme ($\rho_{\infty}=0$) can be examined
using the provided \texttt{MATHEMATICA} code.

\begin{figure}
\begin{centering}
\includegraphics[width=0.45\textwidth]{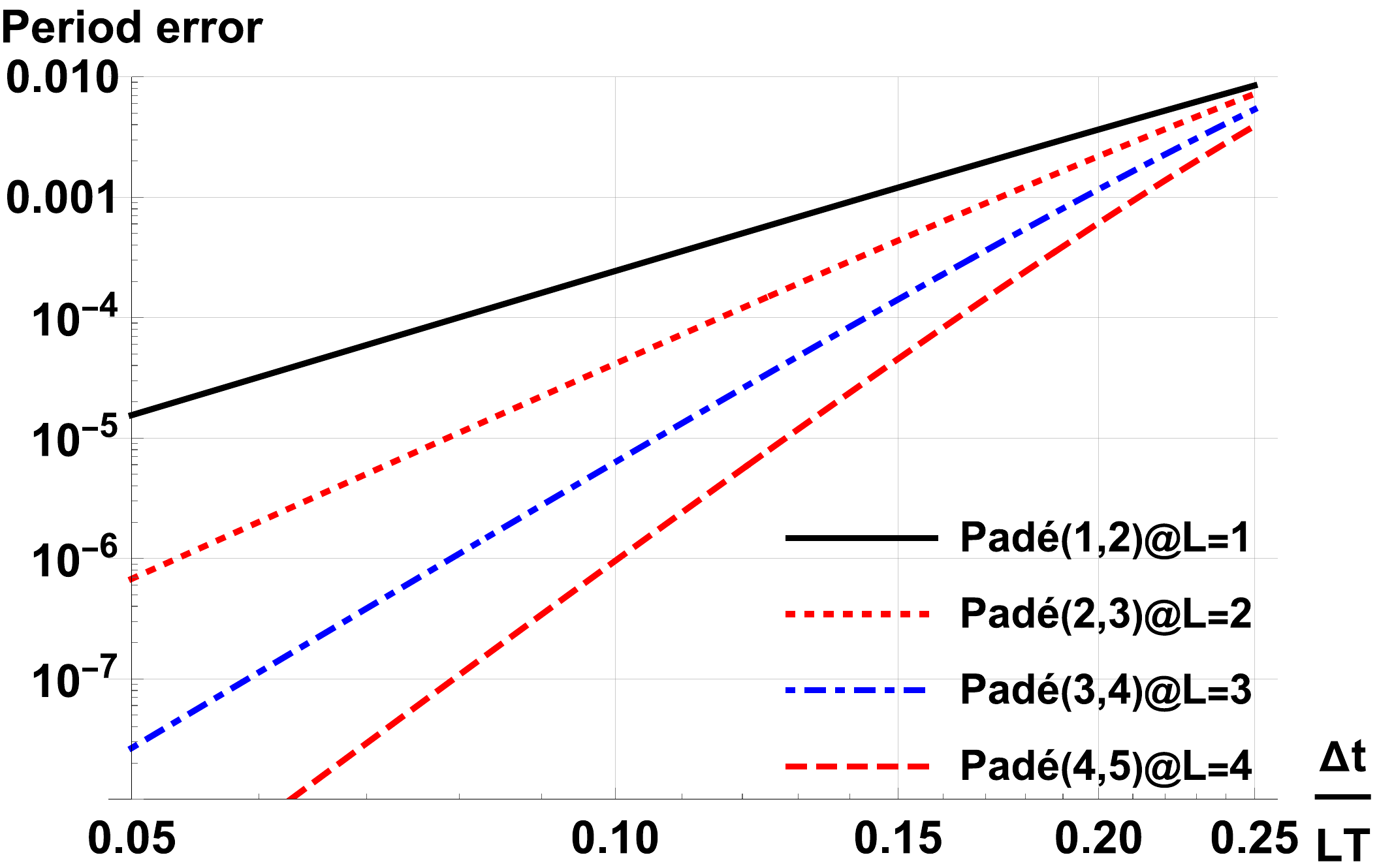}
\par\end{centering}
\medskip{}

\texttt{\small{}a = 0.53846; }{\small\par}

\texttt{\small{}LogLogPlot{[}\{pE{[}1, 2, a, x{]}, pE{[}2, 3, a, 2{*}x{]},
pE{[}3, 4, a, 3{*}x{]}, pE{[}4, 5, a, 4{*}x{]}\}, \{x, 0.05, 0.2\},
PlotRange -> \{10\textasciicircum (-8), 0.01\}{]}}{\small\par}

\caption{Relative period error of high-order scheme with $\rho_{\infty}=0.53846$
(corresponding to $\alpha=-0.3$ for HHT-$\alpha$ scheme). The \texttt{MATHEMATICA}
commands are listed below the plots.\label{fig:TimeStepA053}}
\end{figure}

Similar to the discussion on the HHT-$\alpha$ method, the relative
period errors are shown in Fig\@.~\ref{fig:TimeStepA053} for orders
$(1,2)$, $(2,3)$, $(3,4)$, and $(4,5)$. 
To summarize the observations and conclusions concisely, the order
of the high-order scheme is denoted as $(L=M-1,M)$. Note that the
time step size $\Delta t/T$ is scaled according to the order $L$
by a factor $1/L$ in the horizontal axis of the plot. For example,
the period error of the order $(4,5)$ scheme with $L=4$ is found
from Fig\@.~\ref{fig:TimeStepA053} to be about $10^{-6}$ at $\Delta t/T=4\times0.1=0.4$.
It is noted that the relative period error of the high-order scheme
is much smaller than that of the HHT-$\alpha$ scheme and should result
in significantly higher accuracy. 

\section{Numerical examples\label{sec:Numerical-examples}}

The dissipative and dispersive characteristics of the proposed scheme
with controllable numerical damping have been analyzed using the free
vibration of a single-degree-of-freedom problem in Section~\ref{subsec:Discussion-on-selection},
where the selection of time step size is also discussed. In this section,
the performance of solving a model problem with a high stiffness ratio
is examined first. Wave propagation problems discretized with finite
elements are then addressed. Finally, a simple guideline for the selection
of $\rho_{\infty}$ and time step size $\Delta t$ is proposed for
wave propagation problems.

We focus on evaluating the performance of the proposed scheme on the
numerical dissipation of spurious high-frequency oscillations without
affecting the low-frequency responses. Linear (1st-order) elements
are utilized for the spatial discretization, which are known to lead
to significant errors of spatially unresolved high-frequency modes
in the semi-discretized equation of motion \citep{Cottrell2006}.
The use of high-order formulations to reduce the spatial discretization
error is out of the scope of the present work and will be addressed
in forthcoming publications.

The first five examples in this section are selected from the literature.
The results obtained with the Bathe and other high-order methods are
available in the cited references and can be directly used for comparison.

In our previous work \citep{Song2022}, the convergence of the present
method with respect to the time step size (and no numerical damping)
has been studied. The computational times and accuracy are compared
with the Newmark method, which has a similar computational cost to
the HHT-$\alpha$ method. Since the difference in computer times taken
by the present method with or without numerical dissipation is minor,
the conclusions related to the computer times reached in \citep{Song2022}
are still valid. As this paper focuses on the effectiveness of numerical
dissipation, the evaluation of effectiveness and computer time in
comparison with the HHT-$\alpha$ method is performed for given finite
element models (semi-discretized systems). The evaluation is discussed
in Section~\ref{subsec:Three-dimensional-wave-propagati} on the
large-scale simulation of a 3D sandwich panel.

Source codes of the proposed scheme written in \texttt{MATLAB} and
\texttt{FORTRAN} are available for download at \href{https://github.com/ChongminSong/HighOrderTimeIntegration}{https://github.com/ChongminSong/HighOrderTimeIntegration}.
Interested readers may use the source codes to compare the computational
cost and accuracy with other time-integration methods on their computer
systems. 

\subsection{A three-degree-of-freedom model problem}

\begin{figure}
\begin{centering}
\includegraphics[width=0.5\textwidth]{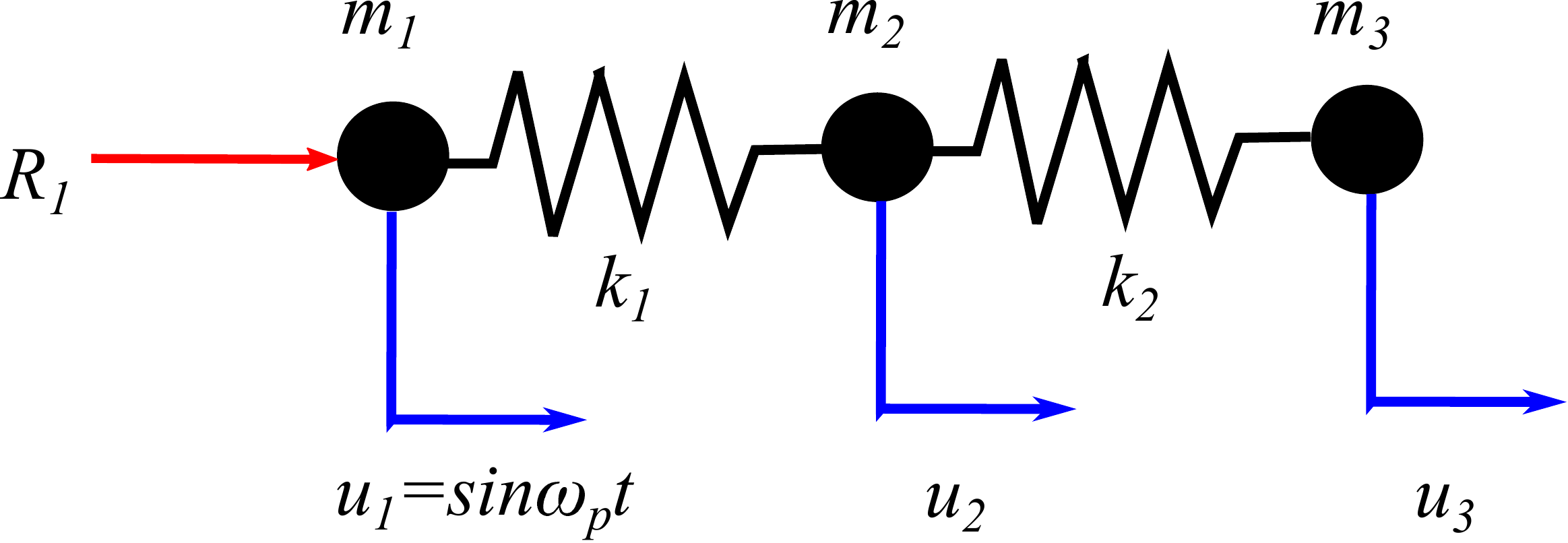}
\par\end{centering}
\caption{A three-degree-of-freedom model problem\label{fig:A-linear-three-degree-of-freedom}}
\end{figure}

The model problem studied in \citep{Noh2018} and \citep{Choi2022}
is shown in Fig.~\ref{fig:A-linear-three-degree-of-freedom}. It
consists of three masses connected by two springs with a high stiffness
ratio. The spring and mass coefficients are given as $k_{1}=10^{7}$,
$k_{1}=1$, $m_{1}=0$, $m_{2}=1$ and $m_{3}=1$. The displacement
of $m_{1}$ is prescribed as 
\begin{equation}
u_{1}=\sin\omega_{p}t\label{eq:3dof0}
\end{equation}
with $\omega_{p}=1.2$ corresponding to the period of vibration $T_{p}=5.236$.
The equation of motion of the system is expressed 
\begin{equation}
\left[\begin{array}{cc}
m_{2} & 0\\
0 & m_{3}
\end{array}\right]\left\{ \begin{array}{c}
\ddot{u}_{2}\\
\ddot{u}_{3}
\end{array}\right\} +\left[\begin{array}{cc}
k_{1}+k_{2} & -k_{2}\\
-k_{2} & k_{2}
\end{array}\right]\left\{ \begin{array}{c}
u_{2}\\
u_{3}
\end{array}\right\} =\left\{ \begin{array}{c}
k_{1}u_{1}\\
0
\end{array}\right\} \label{eq:3dof1}
\end{equation}
with the displacements $u_{2}$ of $m_{2}$ and $u_{3}$ of $m_{3}$.
The reaction force at $m_{1}$ is equal to
\begin{equation}
R_{1}=m_{1}\ddot{u}_{1}+k_{1}(u_{1}-u_{2})\,.\label{eq:3dof2}
\end{equation}
The system is initially at rest, $u_{2}(0)=\dot{u}_{2}(0)=u_{3}(0)=\dot{u}_{3}(0)=0$.

Equation~\eqref{eq:3dof1} can be solved by mode superposition. The
two natural frequencies are approximately equal to $\omega_{1}=1$
and $\omega_{2}=3162$, which correspond to the periods of vibration
$T_{1}=6.283$ and $T_{2}=0.002$. The reference solution is obtained
by excluding the participation of the high-frequency mode $\omega_{2}$
from the mode superposition solution.

The time integration of the present high-order scheme is performed
with the parameter $\rho_{\infty}=0$ and the time step size of $\Delta t=0.14$
as being used in \citep{Choi2022}. This choice results in $\Delta t/T_{p}=0.0267$,
$\Delta t/T_{1}=0.0223$ and $\Delta t/T_{2}=70.5$. Each period of
excitation is divided into about 37 time steps. It is expected from
Fig.~\ref{fig:Spectral-radii-M-1} that the high-frequency oscillations
with the period $T_{2}$ are rapidly damped. The analysis is performed
for a long duration of $t=5000\approx967T_{p}$.

\begin{figure}
\includegraphics[width=0.32\textwidth]{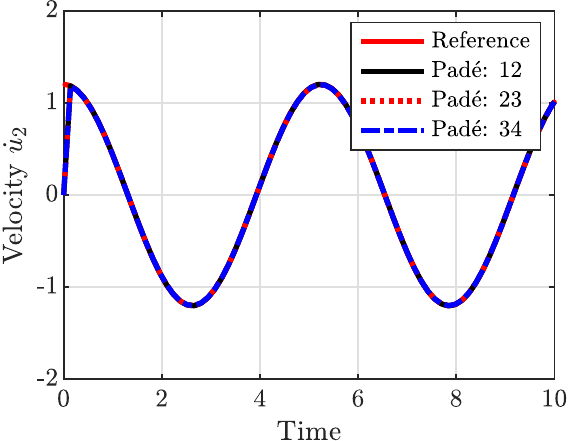}\hfill{}\includegraphics[width=0.32\textwidth]{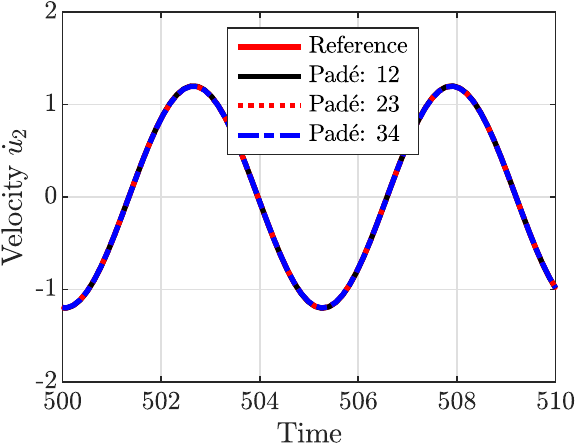}\hfill{}\includegraphics[width=0.32\textwidth]{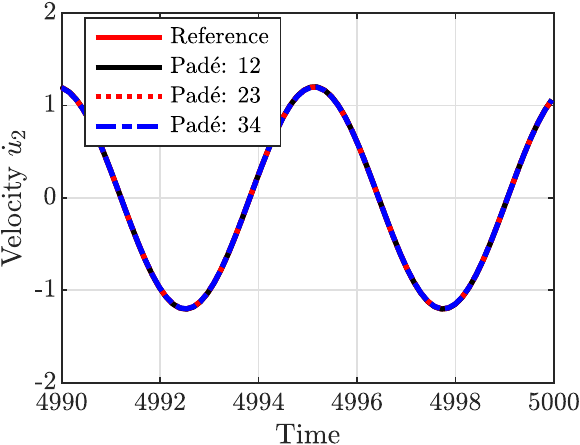}

\includegraphics[width=0.32\textwidth]{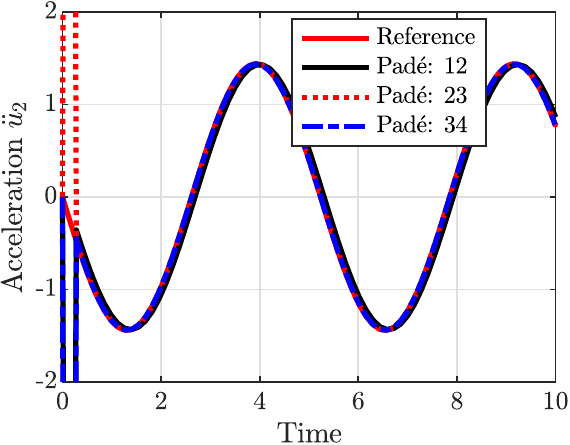}\hfill{}\includegraphics[width=0.32\textwidth]{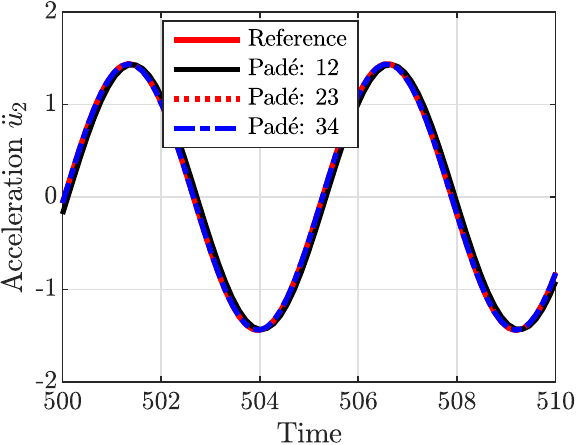}\hfill{}\includegraphics[width=0.32\textwidth]{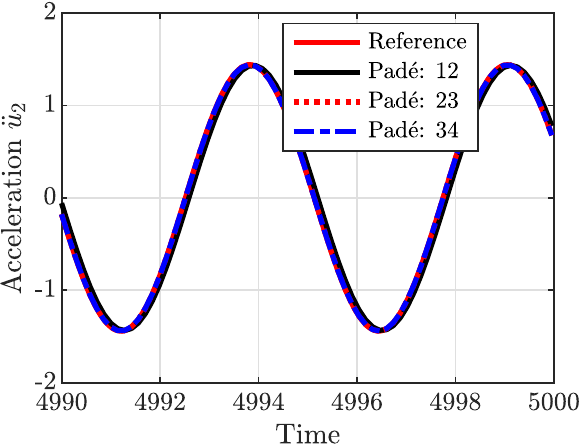}

\caption{Velocity (top row) and acceleration (bottom row) responses of $m_{2}$
during $0\protect\leq t\protect\leq10$ (left column), $500\protect\leq t\protect\leq510$
(middle column) and $4900\protect\leq t\protect\leq5000$ (right column).
\label{fig:3dofs_m2}}
\end{figure}
\begin{figure}
\includegraphics[width=0.32\textwidth]{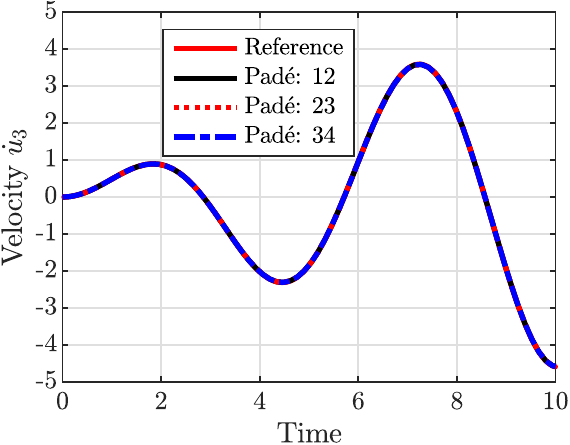}\hfill{}\includegraphics[width=0.32\textwidth]{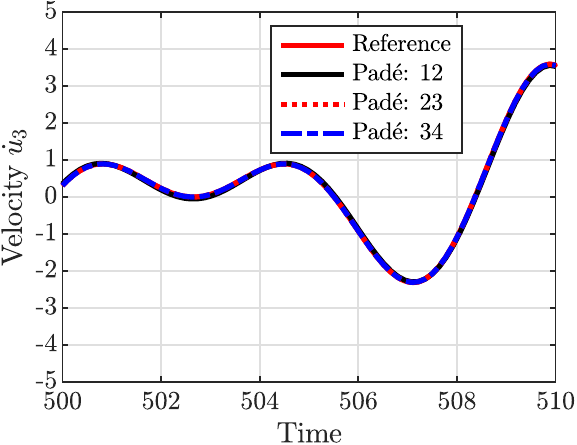}\hfill{}\includegraphics[width=0.32\textwidth]{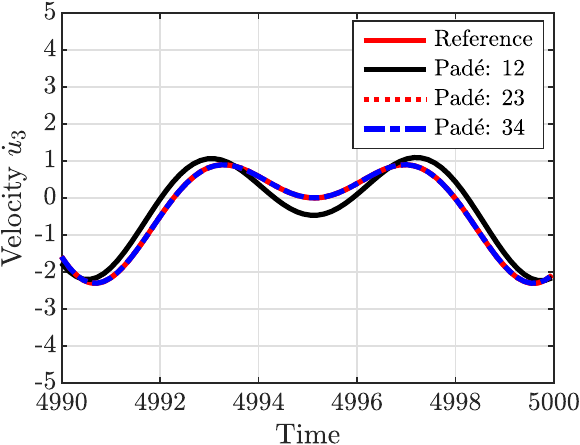}

\includegraphics[width=0.32\textwidth]{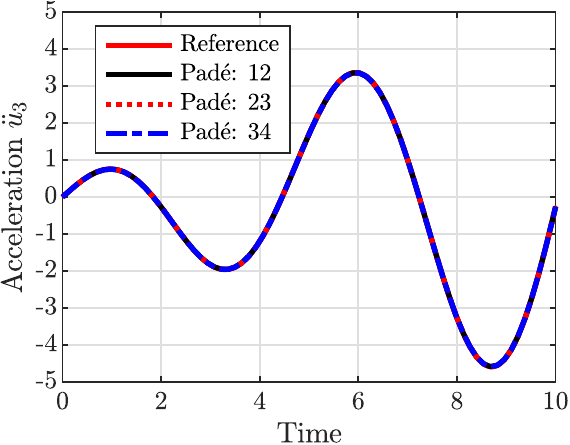}\hfill{}\includegraphics[width=0.32\textwidth]{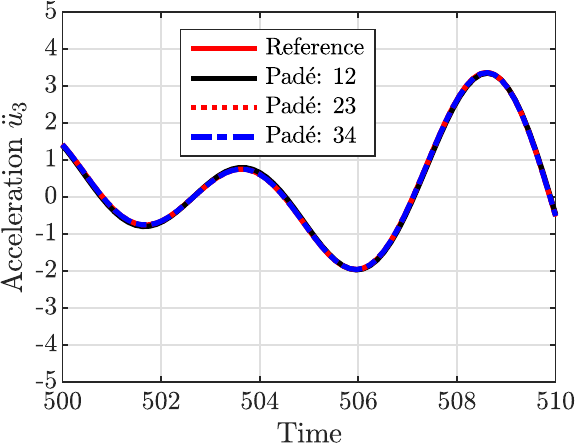}\hfill{}\includegraphics[width=0.32\textwidth]{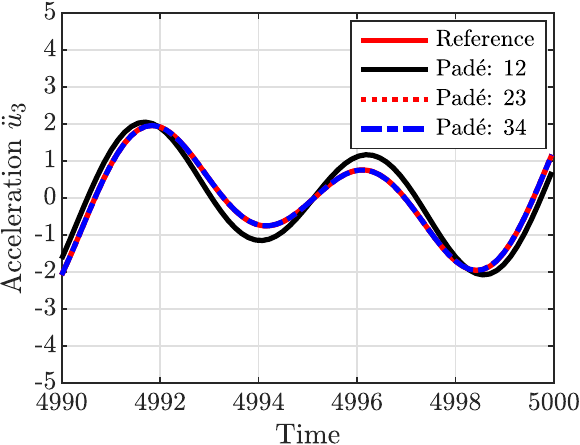}

\hfill{}

\caption{Velocity (top row) and acceleration (bottom row) responses of $m_{3}$
during $0\protect\leq t\protect\leq10$ (left column), $500\protect\leq t\protect\leq510$
(middle column) and $4900\protect\leq t\protect\leq5000$ (right column).
\label{fig:3dofs_m3}}
\end{figure}

\begin{figure}
\hfill{}

\hfill{}\includegraphics[width=0.32\textwidth]{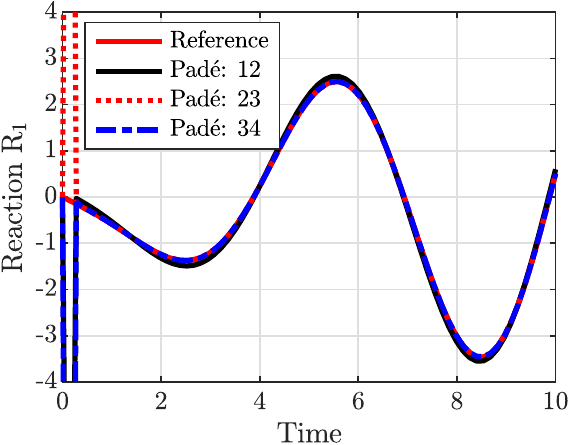}\hfill{}\includegraphics[width=0.32\textwidth]{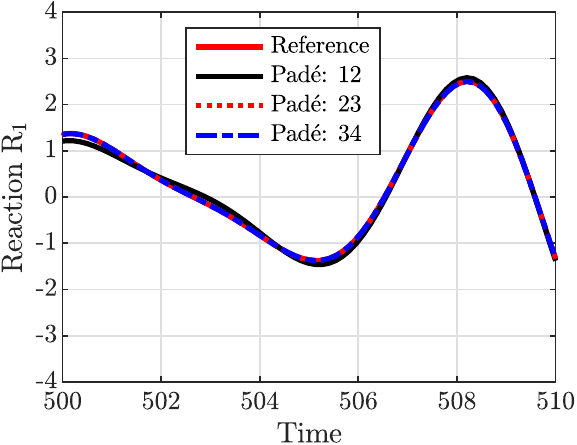}\hfill{}\includegraphics[width=0.32\textwidth]{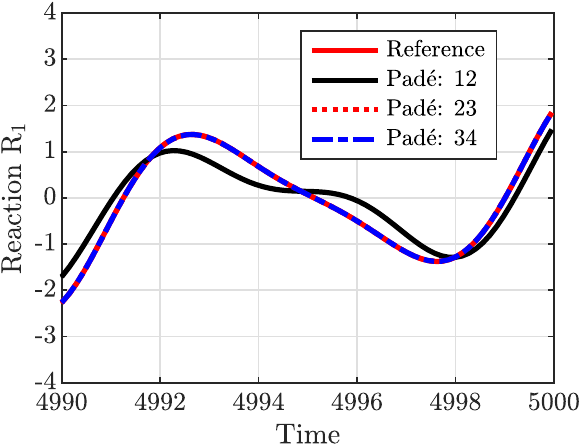}

\caption{Reaction force response during $0\protect\leq t\protect\leq10$ (left
column), $500\protect\leq t\protect\leq510$ (middle column) and $4900\protect\leq t\protect\leq5000$
(right column). \label{fig:3dofs_R1}}
\end{figure}

Figures~\ref{fig:3dofs_m2} and \ref{fig:3dofs_m3} show the velocity
(in the top row) and acceleration (in the bottom row) responses of
$m_{2}$ and $m_{3}$, respectively. The responses of the reaction
force $R_{1}(t)$ are plotted in Fig.~\ref{fig:3dofs_R1}. The three
columns of each figure show the responses at three different time
intervals: $0\leq t\leq10$ (left column), $500\leq t\leq510$ (middle
column) and $4900\leq t\leq5000$ (right column). It is observed from
Fig.~\ref{fig:3dofs_m2} that the initial velocity $\dot{u}_{2}(0)=1$
is inconsistent with the initial condition ($\dot{u}_{2}(0)=0$) as
the result of excluding the high-frequency vibrations. This inconsistency
leads to a spike in the first time step of the acceleration $\ddot{u}_{2}$
in Fig.~\ref{fig:3dofs_m2} and the reaction force $R_{1}$ in Fig.~\ref{fig:3dofs_R1}.
After the first step, the result obtained at order $(1,2)$ differs
slightly from the reference resolution. The increase of the difference
with time (from the left column with $0\leq t\leq10$ to the right
column with $4900\leq t\leq5000$) is appreciable in the responses
of $m_{3}$ (Fig.~\ref{fig:3dofs_m3}) and reaction force $R_{1}$
(Fig.~\ref{fig:3dofs_R1}). The results obtained at orders $(2,3)$
and $(3,4)$ are indistinguishable from the reference solution throughout
the whole duration, showing negligible numerical dissipation and phase
error of the low-frequency mode. This example illustrates that high-order
schemes are advantageous for analyses of long duration. The results
reported in \citep{Noh2018} using the Newmark method and \citep{Noh2018,Choi2022}
using the Bathe method are available for comparison with the present
results.

\subsection{One-dimensional wave propagation in a homogeneous rod \label{subsec:1D rod}}

\begin{figure}
\centering{}\includegraphics{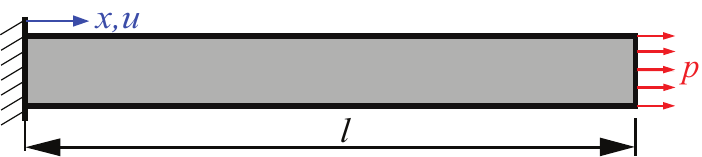}\caption{A homogeneous elastic rod subjected to a step loading. \label{fig:geometry beam}}
\end{figure}

The problem of elastic wave propagation in a one-dimensional prismatic
rod, as sketched in Fig.~\ref{fig:geometry beam}, is frequently
used in the literature when studying the numerical dissipation properties
of time integration methods. The material and geometrical parameters
are adopted from Ref.~\citep{Malakiyeh2019} with a consistent set
of units: length of the rod $l=200,$ Young's modulus $E=3\times10^{7}$,
Poisson's ratio $\nu=0.0$, and mass density $\rho=0.00073$. The
longitudinal wave speed is $c=\sqrt{E/\rho}=\sqrt{3\times10^{7}/0.00073}=2.0272\times10^{5}$.
We assume that no physical damping is present. The left end of the
rod is fixed and the right end is subjected to a step loading $pH(t)$
($H(t)$ denotes the Heaviside function) with the amplitude of the
pressure $p=10^{4}$. Since the step loading includes high-frequency
components, it will excite some spurious high-frequency modes of the
finite element model and therefore, numerical dissipation is desired.

The effect of numerical dissipation on the response at a frequency
(or a period $T$) depends on the time step size $\Delta t$. For
a given mesh, the time step size is often expressed as the Courant-Friedrichs-Lewy
(CFL) number defined in 1D as
\begin{equation}
\mathrm{CFL}=c\frac{\Delta t}{\Delta x}\,,\label{eq:CFL}
\end{equation}
where $\Delta x$ denotes the element size. The CFL number measures
the number of elements that the wave can travel in each time step.
In a numerical analysis of wave propagation by direct time integration,
the error mainly comes from two sources, the time discretization error
controlled by the time step size $\Delta t$ and the spatial discretization
error controlled by the element size $\Delta x$. For a wave with
a period $T$, Eq.~\eqref{eq:CFL} can be rewritten as
\begin{equation}
\dfrac{cT}{\Delta x}=\mathrm{CFL}\dfrac{T}{\Delta t}\,.
\end{equation}
It represents the ratio of the number of elements in one wavelength
($cT$) to the number of time steps in one period. Generally speaking,
the CFL number reflects the relative amount of errors in time and
spatial discretizations. In the remainder of this section, the selection
of the weight factor $\rho_{\infty}$ and the CFL number for linear
finite elements will be discussed.

\subsubsection{Control of high-frequency numerical dissipation by $\rho_{\infty}$\label{subsec:Control-of-high-frequency}}

The effect of the user-specified control parameter $\rho_{\infty}$
on the numerical dissipation is investigated in the following. The
analysis is performed for a time duration of $0.02$. A uniform spatial
discretization with 1,000 elements along the length (element size
of $\Delta x=0.2$) is considered. The selection of the time step
size (the CFL number for the given element size) will be discussed
later in Section~\ref{subsec:1D-Size-of-time}. In this section,
the time step size for the order $(1,2)$ scheme is chosen to have
a CFL number of $10$ leading to $\Delta t=9.8658\times10^{-6}$,
which means in each time step the wave travels through $10$ linear
elements. Following the discussions in Section~\ref{subsec:High-order-timeStep},
the time step sizes of orders $(2,3)$, $(3,4)$ and $(4,5)$ are
chosen by multiplying that of order $(1,2)$ with a factor of 2, 3
and 4, respectively, to introduce a similar amount of numerical dissipation.
These time step sizes correspond to CFL numbers of 20 for order $(2,3)$,
30 for order $(3,4)$, and 40 for order $(4,5)$. In the subsequent
analyses, the velocity at the middle point of the rod will be examined.
For all the cases considered below, accurate results for displacement
responses are obtained and will not be reported explicitly.

\begin{figure}
\centering{}\subfloat[Time histories. \label{fig:Velocity-pade-rho1}]{\includegraphics[width=0.45\textwidth]{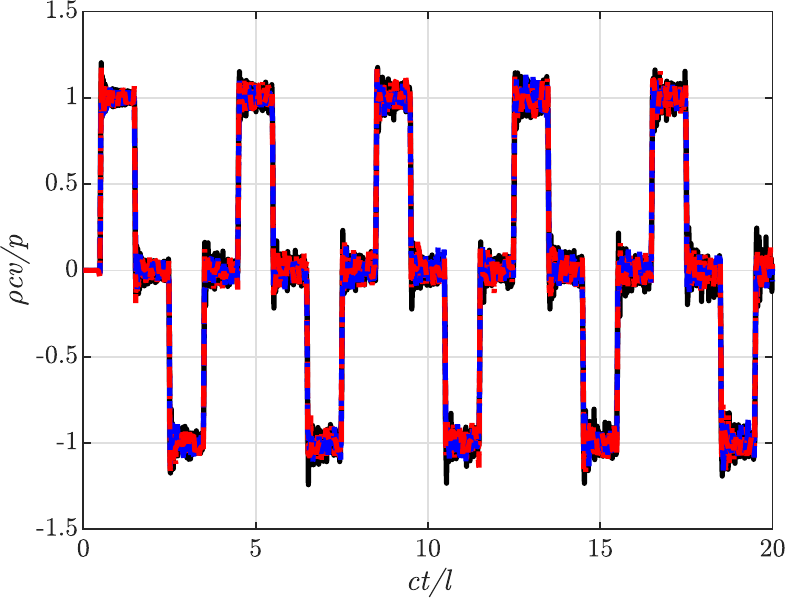}}\hfill{}\subfloat[Close-up view from $ct/l=16.3$ to $ct/l=17.7$. \label{fig:Section-A-A-pade-rho1}]{\includegraphics[width=0.45\textwidth]{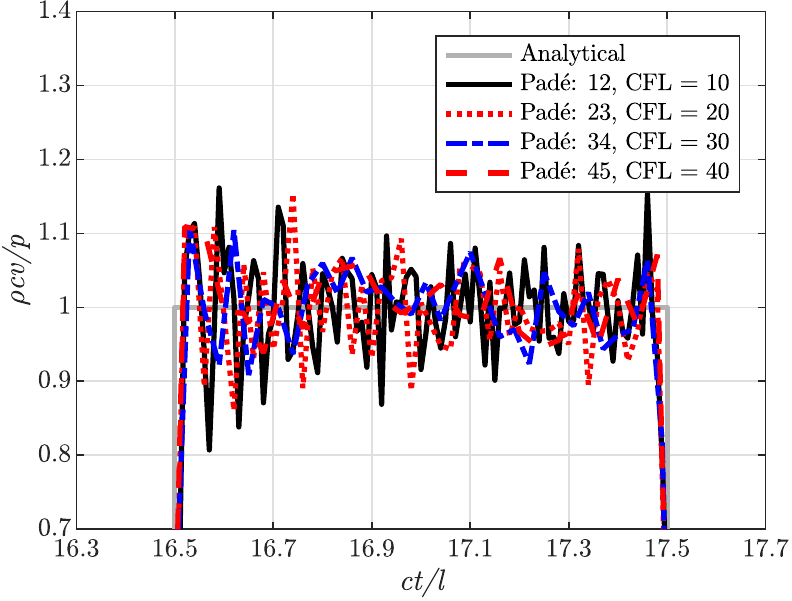}}\caption{Velocity responses at the middle point of 1D rod obtained by the high-order
scheme with $\rho_{\infty}=1$. \label{fig:Velocity-rho1}}
\end{figure}

\begin{figure}
\centering{}\subfloat[Time histories. \label{fig:Velocity-pade-rho0}]{\includegraphics[width=0.45\textwidth]{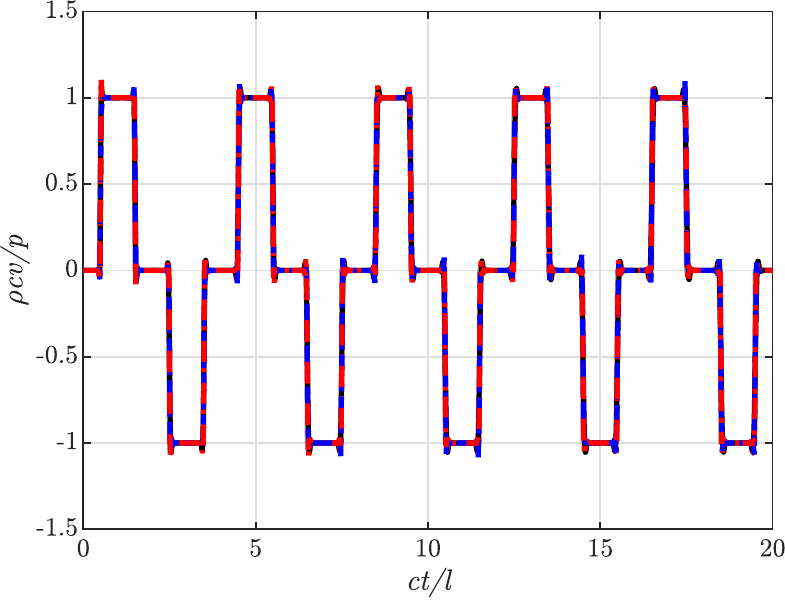}}\hfill{}\subfloat[Close-up view from $ct/l=16.3$ to $ct/l=17.7$. \label{fig:Section-A-A-pade-rho0}]{\includegraphics[width=0.45\textwidth]{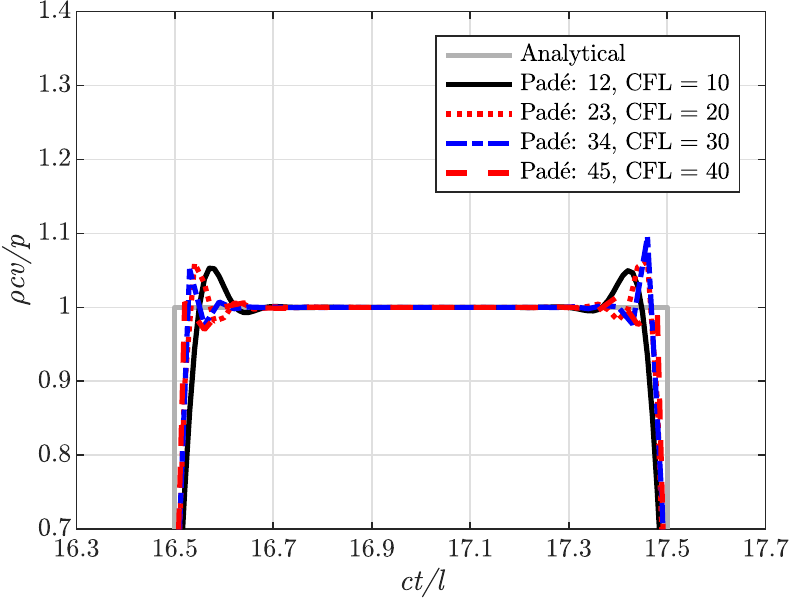}}\caption{Velocity responses at the middle point of 1D rod obtained by the high-order
scheme with $\rho_{\infty}=0$. \label{fig:Velocity-rho0}}
\end{figure}

When the parameter $\rho_{\infty}=1$ (diagonal Padé expansions) is
chosen, no numerical dissipation is introduced. The time histories
of the dimensionless velocity $\rho cv/p$ at the middle of the rod
are plotted in Fig.~\ref{fig:Velocity-pade-rho1} versus the dimensionless
time $ct/l$. In addition, the analytical solution is shown by the
solid gray line. The peak value of the response is equal to $\rho cv/p=1$
for the analytical solution. A close-up view of the results from the
dimensionless time $ct/l=16.3$ to $ct/l=17.7$ (corresponding to
$t=0.016$ to $t=0.0175$) is depicted in Fig.~\ref{fig:Section-A-A-pade-rho1}.
Strong spurious high-frequency oscillations can be observed at all
orders of the scheme. When $\rho_{\infty}=0$ is chosen, the scheme
is $L-$stable and the maximum amount of numerical dissipation is
introduced. The velocity responses are plotted in Fig.~\ref{fig:Velocity-rho0}.
It is observed that the high-frequency oscillations are largely suppressed
by the addition of numerical dissipation.

\begin{figure}
\centering{}\subfloat[Time histories\label{fig:Velocity-pade-rho09}]{\includegraphics[width=0.45\textwidth]{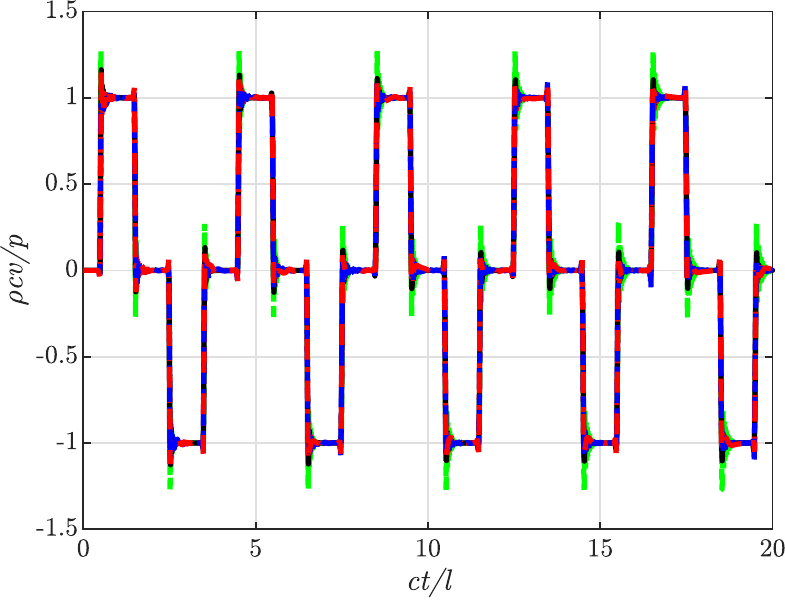}}\hfill{}\subfloat[Close-up view from $ct/l=16.3$ to $ct/l=17.7$.\label{fig:Section-A-A-pade-rho09}]{\includegraphics[width=0.45\textwidth]{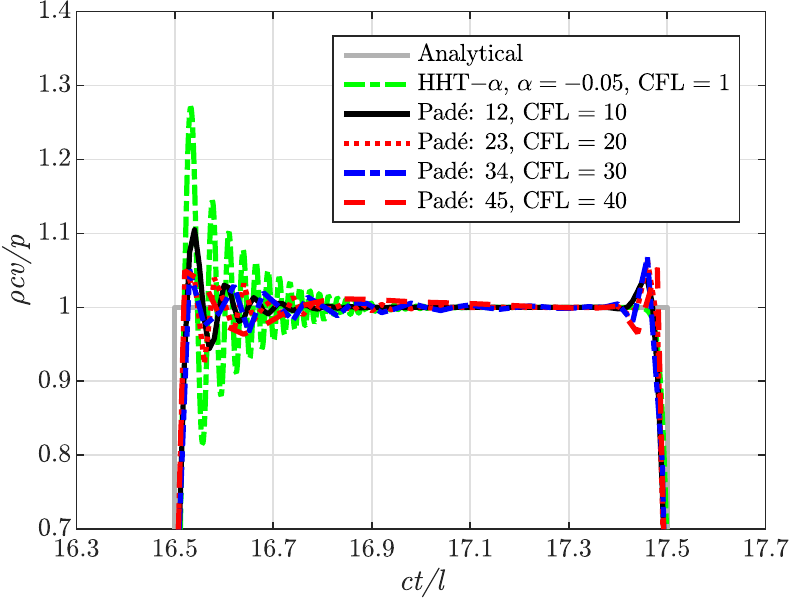}

}\caption{Velocity responses at the middle point of 1D rod obtained by the proposed
high-order scheme with $\rho_{\infty}=0.9$ and HHT-$\alpha$ method
with $\alpha=-0.05$.\label{fig:Velocity-rho09}}
\end{figure}

\begin{figure}
\centering{}\subfloat[Time histories. \label{fig:Velocity-pade-rho08}]{\includegraphics[width=0.45\textwidth]{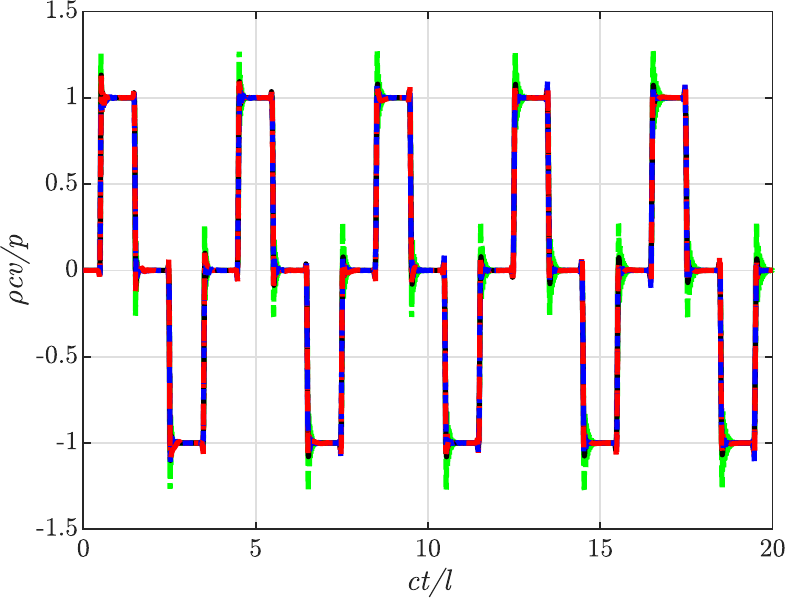}

}\hfill{}\subfloat[Close-up view from $ct/l=16.3$ to $ct/l=17.7$. \label{fig:Section-A-A-pade-rho08}]{\includegraphics[width=0.45\textwidth]{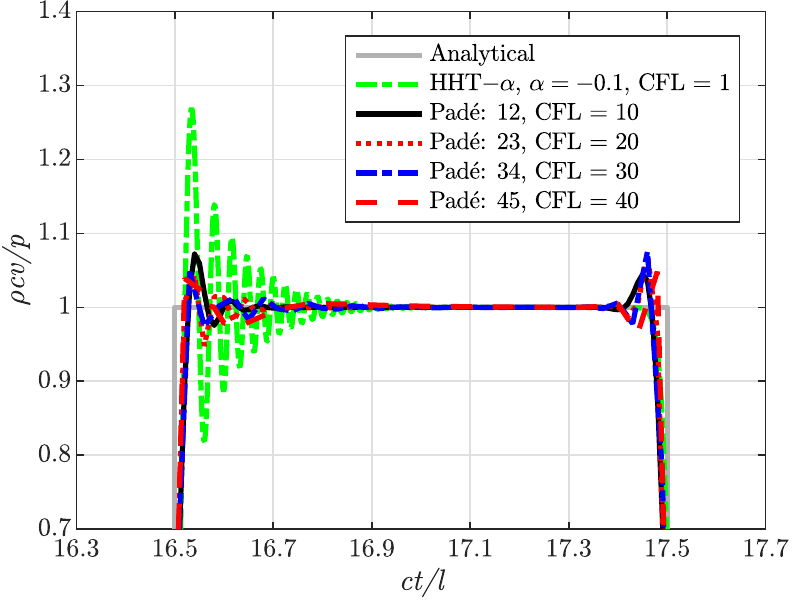}

}\caption{Velocity responses at the middle point of 1D rod obtained by the proposed
high-order scheme with $\rho_{\infty}=0.8$ and HHT-$\alpha$ method
with $\alpha=-0.1$. \label{fig:Velocity-rho08}}
\end{figure}

\begin{figure}
\centering{}\subfloat[Time histories. \label{fig:Velocity-pade-rho05}]{\includegraphics[width=0.45\textwidth]{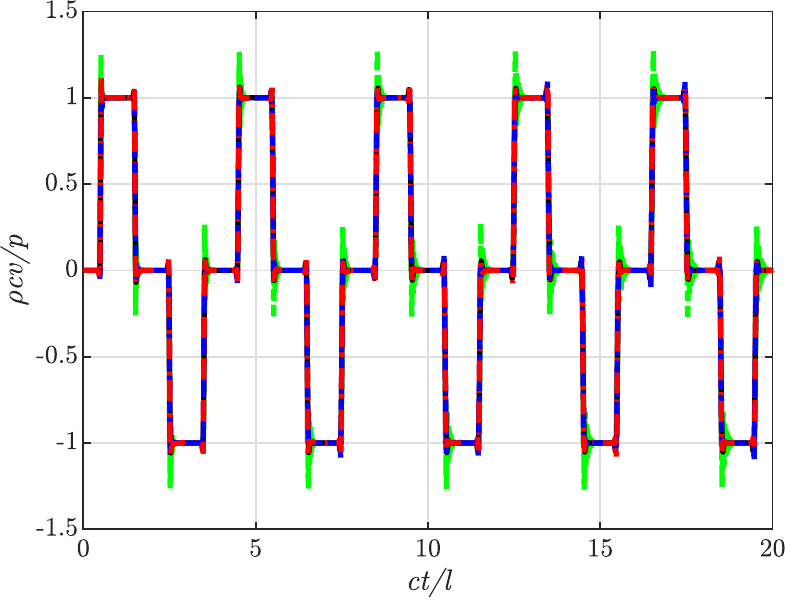}}\hfill{}\subfloat[Close-up view from $ct/l=16.3$ to $ct/l=17.7$. \label{fig:Section-A-A-pade-rho05-2}]{\includegraphics[width=0.45\textwidth]{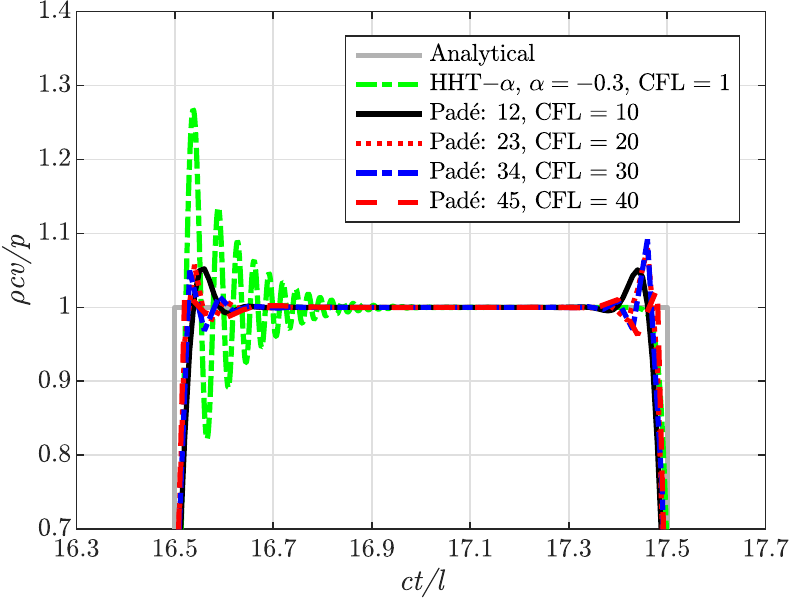}}\caption{Velocity responses at the middle point of 1D rod obtained by the proposed
high-order scheme with $\rho_{\infty}=0.5$ and HHT-$\alpha$ methods
with $\alpha=-0.3$. \label{fig:Velocity-rho05}}
\end{figure}

The choice of the parameter $\rho_{\infty}$ is investigated by performing
the analyses using three additional values of $\rho_{\infty}$ that
correspond to typical values of the parameter $\alpha$ in the HHT-$\alpha$
method (For simplicity, we round the value $\rho_{\infty}$ to a single
significant digit):
\begin{enumerate}
\item $\rho_{\infty}=0.9$, approximately corresponds to $\alpha=-0.05$
in the HHT-$\alpha$ method. This case is commonly regarded as lightly
dissipative. The response histories of velocity at the middle of the
rod are shown in Fig.~\ref{fig:Velocity-rho09}. 
\item $\rho_{\infty}=0.8$ corresponding to $\alpha=-0.1$. The response
histories are shown in Fig.~\ref{fig:Velocity-rho08}. 
\item $\rho_{\infty}=0.5$ corresponding to $\alpha=-0.3$. This case is
close to the maximum amount of numerical dissipation that can be introduced
with the HHT-$\alpha$ method and is regarded as heavily dissipative.
The velocity histories are shown in Fig.~\ref{fig:Velocity-rho05}.
\end{enumerate}

The results obtained with the HHT-$\alpha$ method at $\textrm{CFL}=1$
are also shown for comparison. It is observed that the proposed method
at any order is able to suppress the spurious high-frequency oscillations.
The duration and peak of spurious oscillations are smaller than those
in the results obtained with HHT-$\alpha$ method. The effect of increasing
the numerical dissipation by varying the parameter from $\rho_{\infty}=0.9$
to $\rho_{\infty}=0$ is minor. Therefore, $\rho_{\infty}$ can be
selected from a rather wide range, say between $0.9$ and $0$, to
effectively suppress spurious high-frequency oscillations. In the
remainder of this section, a value of $\rho_{\infty}=0.8$ will be
used as the controlling parameter for the high-order time-stepping
scheme.

\subsubsection{Effect of the time step size\label{subsec:1D-Size-of-time} }

When using a finite element model, the spatial discretization error,
more specifically the highest frequency that can be accurately resolved
by a given mesh, needs to be considered when choosing a suitable time
step size. If $\Delta t$ is too small, spurious high-frequency oscillations
will not be sufficiently damped. If $\Delta t$ is too large, low-frequency
modes that can be accurately resolved by the mesh will be unnecessarily
damped. As discussed previously, the time step size $\Delta t$ is
represented by the CFL number. In the following, both a uniform mesh
and non-uniform mesh are considered.

\begin{figure}
\centering{}\subfloat[Padé order $(1,2)$\label{fig:Velocity-pade12-CFL}]{\includegraphics[width=0.45\textwidth]{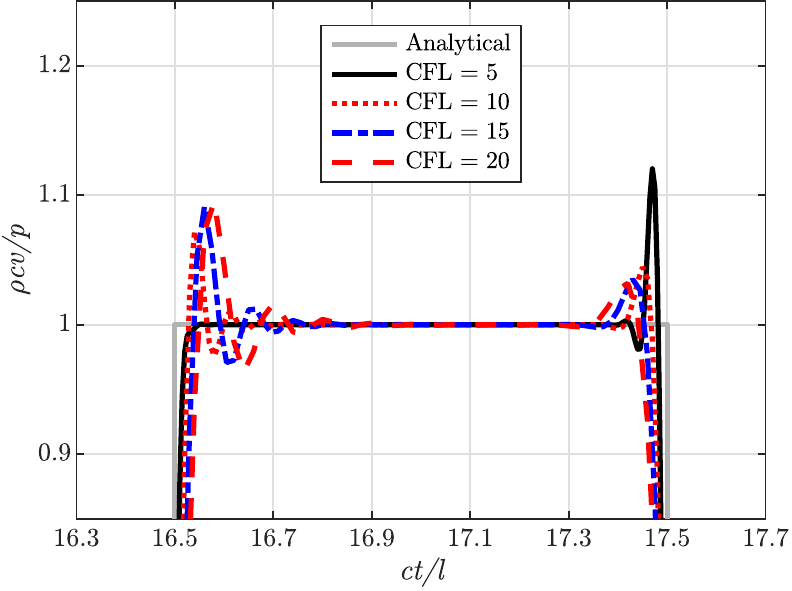}}\hfill{}\subfloat[Padé order $(2,3)$ \label{fig:Velocity-pade23-CFL}]{\includegraphics[width=0.45\textwidth]{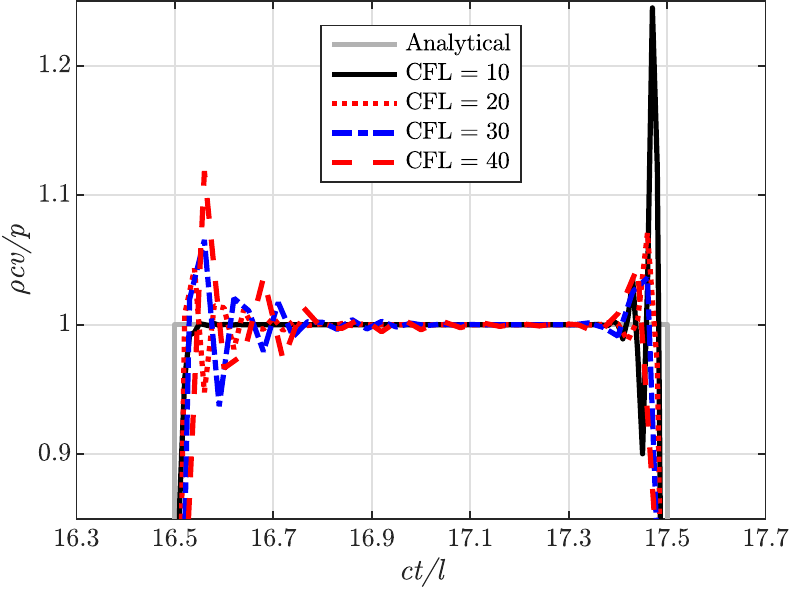}}\medskip{}
\subfloat[Padé order $(3,4)$\label{fig:Velocity-pade34-CFL}]{\includegraphics[width=0.45\textwidth]{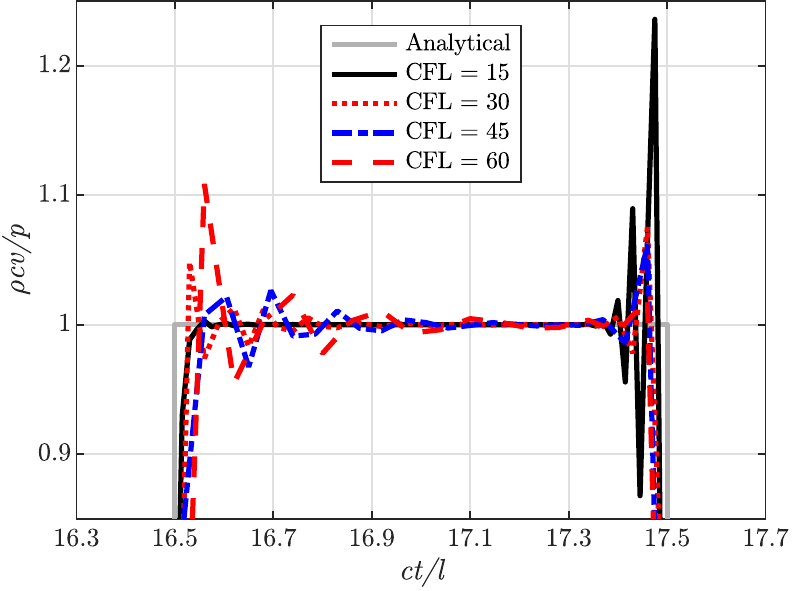}}\hfill{}\subfloat[Padé order $(4,5)$ \label{fig:Velocity-pade45-CFL}]{\includegraphics[width=0.45\textwidth]{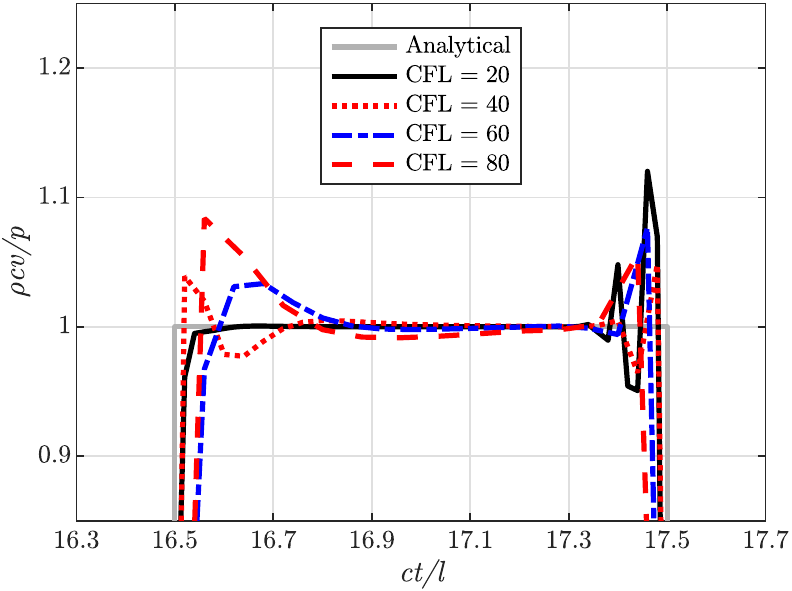}}\caption{Velocity responses at middle point of 1D rod obtained by the proposed
high-order scheme at various CFL numbers with $\rho_{\infty}=0.8$.\label{fig:Velocity-pade-CFL}}
\end{figure}

The analyses are performed using a uniform mesh consisting of $1{,}000$
elements for the orders $(1,2)$, $(2,3)$, $(3,4)$ and $(4,5)$.
The parameter $\rho_{\infty}=0.8$ is chosen. The velocity response
histories at the middle point of the rod are similar to those in Fig.~\ref{fig:Velocity-rho08}.
Therefore, only the close-up views of the response histories from
$ct/l=16.3$ to $ct/l=17.7$ (i.e., $t=0.016$ to $t=0.0175$) are
provided in Fig.~\ref{fig:Velocity-pade-CFL}. Following the observations
made in Section~\ref{subsec:High-order-timeStep}, the CFL numbers
are chosen according to the order of the scheme. Denoting the order
as $(L=M-1,M)$, four different values of CFL number $5L$, $10L$,
$15L$, and $20L$ are considered for each order $L=2$, $3$, $4,$
and $5$. It is observed by comparing Fig.~\ref{fig:Velocity-rho08}b
to Fig.~\ref{fig:Velocity-pade-CFL} that all the CFL numbers lead
to better results than that of the HHT-$\alpha$ method. The results
corresponding to the lowest CFL number $5L$ show less than desirable
numerical dissipation, while all the other CFL numbers lead to similar
results. Therefore, the CFL numbers can be chosen from a rather wide
range, and there are no obvious benefits in identifying optimal CFL
numbers. In the remainder of this section, the CFL number is chosen
as $10L$ for the order $(L=M-1,M)$ scheme. It is confirmed by a
parametric study that this choice of the CFL number is also suitable
for any value of $\rho_{\infty}$ between 0 and 1.

\begin{figure}
\centering{}\subfloat[Time histories. \label{fig:Velocity-1D-nonuniform-rho05-1}]{\includegraphics[width=0.45\textwidth]{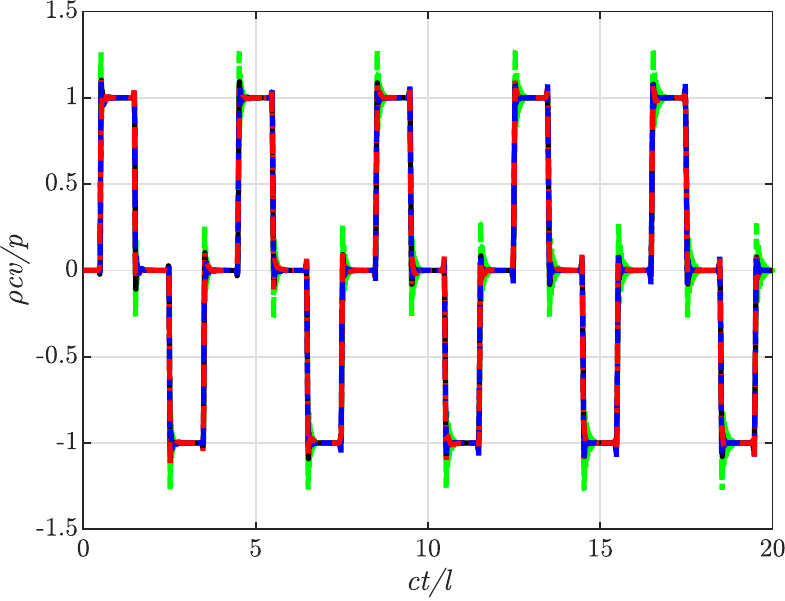}}\hfill{}\subfloat[Close-up view from from $ct/l=16.3$ to $ct/l=17.7$. \label{fig:Velocity-1D-nonuniform-rho08-2}]{\includegraphics[width=0.45\textwidth]{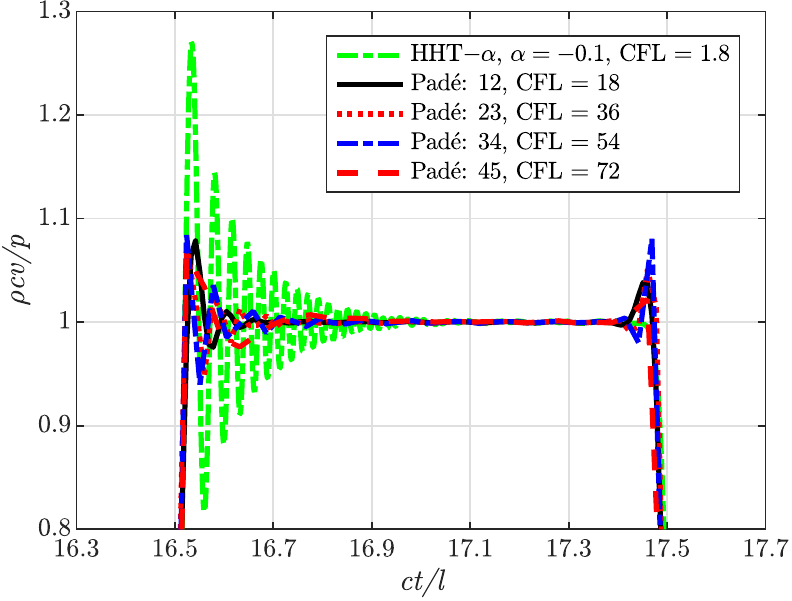}}\caption{Velocity responses at the middle point of 1D rod obtained with a non-uniform
mesh by the proposed high-order scheme with $\rho_{\infty}=0.8$ and
HHT-$\alpha$ method with $\alpha=-0.1$. The CFL numbers shown are
calculated using the average length of all elements.\label{fig:Velocity-1D-nonuniform-rho08}}
\end{figure}

In the next step, a non-uniform mesh is considered as often encountered
in practical 2D or 3D finite element analyses. The coordinate $x_{i}$
of the node $i$ in a mesh of $n_{e}$ elements is given by 
\[
x_{i}=l\times\left(\frac{i-1}{n_{e}}+\dfrac{9}{11\times20\pi}\sin^{2}\left(20\pi\frac{i-1}{n_{e}}\right)\right);\qquad(i=1,2,\ldots,n_{e}+1)\,.
\]
The length of the elements varies as a sinusoidal function with a
mean value of $l/n_{e}$ and a period of $l/20$. The length ratio
of the largest element to the smallest element is about 10. The rod
is divided into $n_{e}=2{,}000$ elements. The size of the largest
and smallest element is equal to $9.088\times10^{-4}l$ and $9.12\times10^{-5}l$,
respectively. Since the highest frequency that the mesh can accurately
resolve is controlled by the largest elements, the time step size
is chosen for the order $(L=M-1,M)$ scheme in such a way that the
CFL number is about $10L$ according to the largest elements. This
choice corresponds to a CFL number of $18L$ in terms of the average
length of elements. The velocity response histories at the middle
point of the rod are shown in Fig.~\ref{fig:Velocity-1D-nonuniform-rho08}.
The spurious oscillations are largely suppressed by the proposed scheme.
The result of the HHT-$\alpha$ method with $\alpha=-0.1$ and CFL$=1.8$
is also shown in the figures, in which strong oscillations are present.

\subsection{One-dimensional wave propagation in a bi-material rod}

\begin{figure}
\centering{}\includegraphics{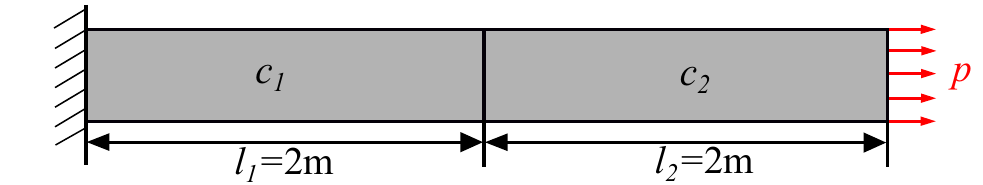}\caption{A bi-material rod subjected to a step loading. \label{fig:A-bi-material-rod}}
\end{figure}

The selection of $\rho_{\infty}$ and time step size $\Delta t$ is
further evaluated using the bi-material rod (taken from Ref.~\citep{Kim2021})
shown in Fig.~\ref{fig:A-bi-material-rod}. The rod consists of two
segments with a length of $2\,\mathrm{m}$ each. The wave speeds of
the the left and right segments are equal to $c_{1}=\unit[40\sqrt{5}]{\,m/s}$
and $c_{2}=\unit[20\sqrt{2}]{\,m/s}$, respectively. The left end
of the rod is fixed. A traction of $p=\unit[1]{Pa}$ is applied as
a step function to the right end.

\begin{figure}
\centering{}\subfloat[Proposed high-order scheme with $\rho_{\infty}=0.8$. \label{fig:Velocity-1D-bimaterial-vel-pade}]{\includegraphics[width=0.45\textwidth]{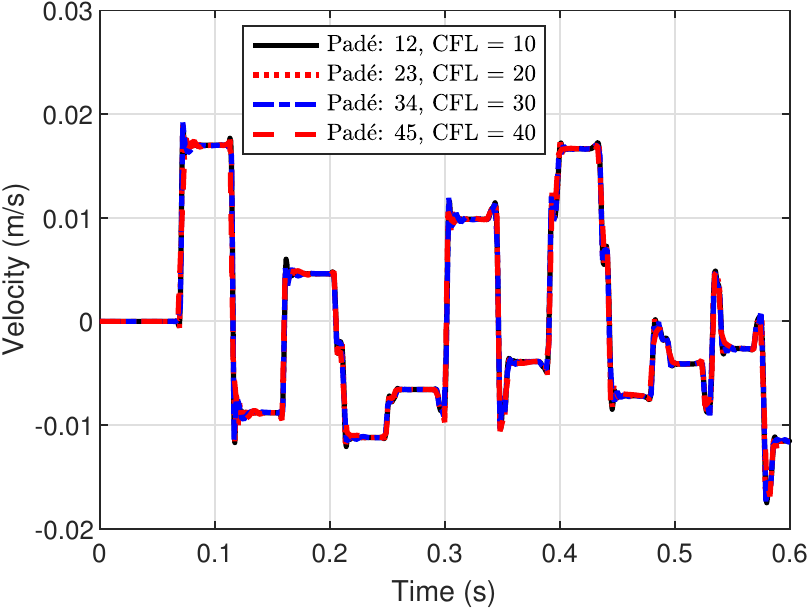}}\hfill{}\subfloat[HHT-$\alpha$ method with $\alpha=-0.1$. \label{fig:Velocity-1D-bimaterial-vel-hht}]{\includegraphics[width=0.45\textwidth]{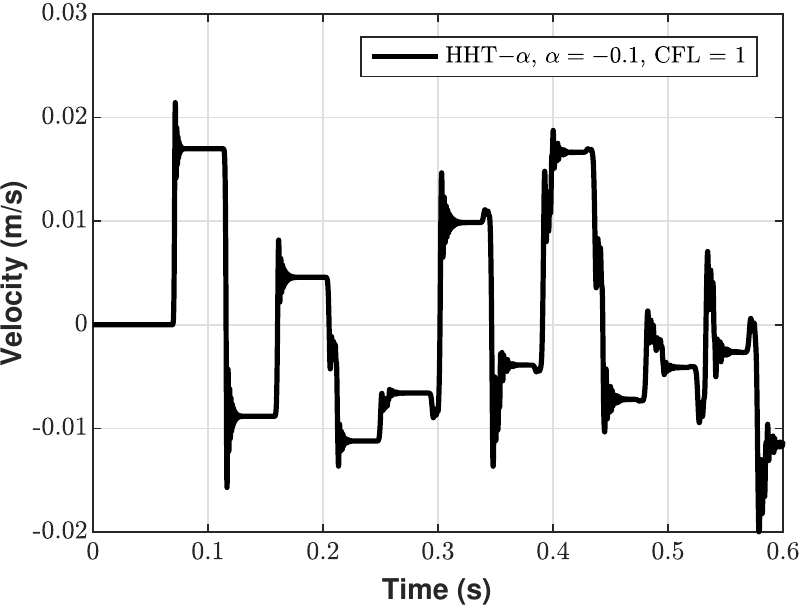}}\caption{Velocity responses at middle point of 1D bi-material rod. The CFL
numbers shown are calculated using the lower wave speed.\label{fig:Velocity-1D-bimaterial-vel}}
\end{figure}

\begin{figure}
\centering{}\subfloat[Proposed high-order scheme with $\rho_{\infty}=0.8$. \label{fig:Velocity-1D-bimaterial-stress-pade}]{\includegraphics[width=0.45\textwidth]{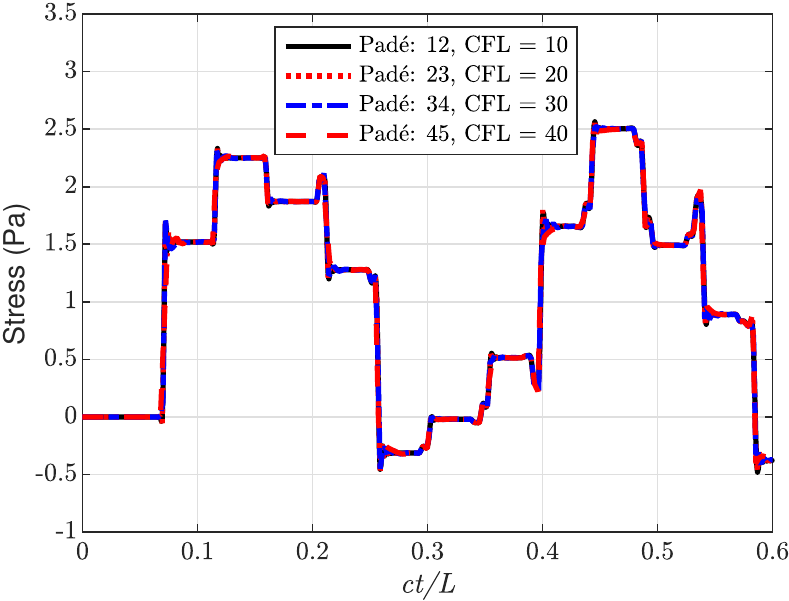}}\hfill{}\subfloat[HHT-$\alpha$ method with $\alpha=-0.1$. \label{fig:Velocity-1D-bimaterial-stress-hht}]{\includegraphics[width=0.45\textwidth]{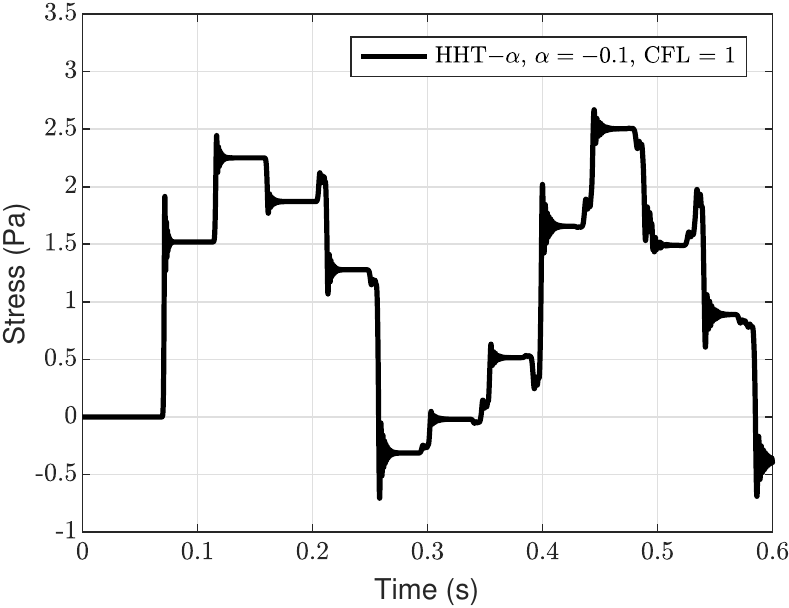}}\caption{Axial stress responses at middle point of 1D bi-material rod. The
CFL numbers shown are calculated using the lower wave speed.\label{fig:Velocity-1D-bimaterial-stress}}
\end{figure}

For the subsequent analyses, a uniform mesh is used, where each segment
is divided into $1{,}000$ linear finite elements. Since the wave
travels within a time step through fewer elements in the right segment
with a lower speed than in the left segment, the spatial discretization
error is relatively higher. Thus, the time step size is determined
from the CFL numbers of the right segment chosen as $\mathrm{CFL}=10L$
for the $(L=M-1,M)$ order scheme. The parameter $\rho_{\infty}=0.8$
is adopted. The velocity and axial-stress response histories at the
interface of the materials (i.e., the middle point of the rod) are
plotted in Figs.~\ref{fig:Velocity-1D-bimaterial-vel-pade} and \ref{fig:Velocity-1D-bimaterial-stress-pade},
respectively. The results obtained using the HHT-$\alpha$ method
with corresponding $\alpha=-0.1$ and $\mathrm{CFL}=1$ are shown
in Figs.~\ref{fig:Velocity-1D-bimaterial-vel-hht} and \ref{fig:Velocity-1D-bimaterial-stress-hht}
for comparison. It is observed that the proposed scheme is more effective
in suppressing spurious high-frequency oscillations. The results obtained
with the Bathe method and an overlapping finite element scheme are
reported in Ref.~\citep{Kim2021}.

\subsection{Scalar wave propagation in a square domain\label{subsec:Scalar-wave-propagation}}

\begin{figure}
\centering{}\includegraphics[width=0.3\textwidth]{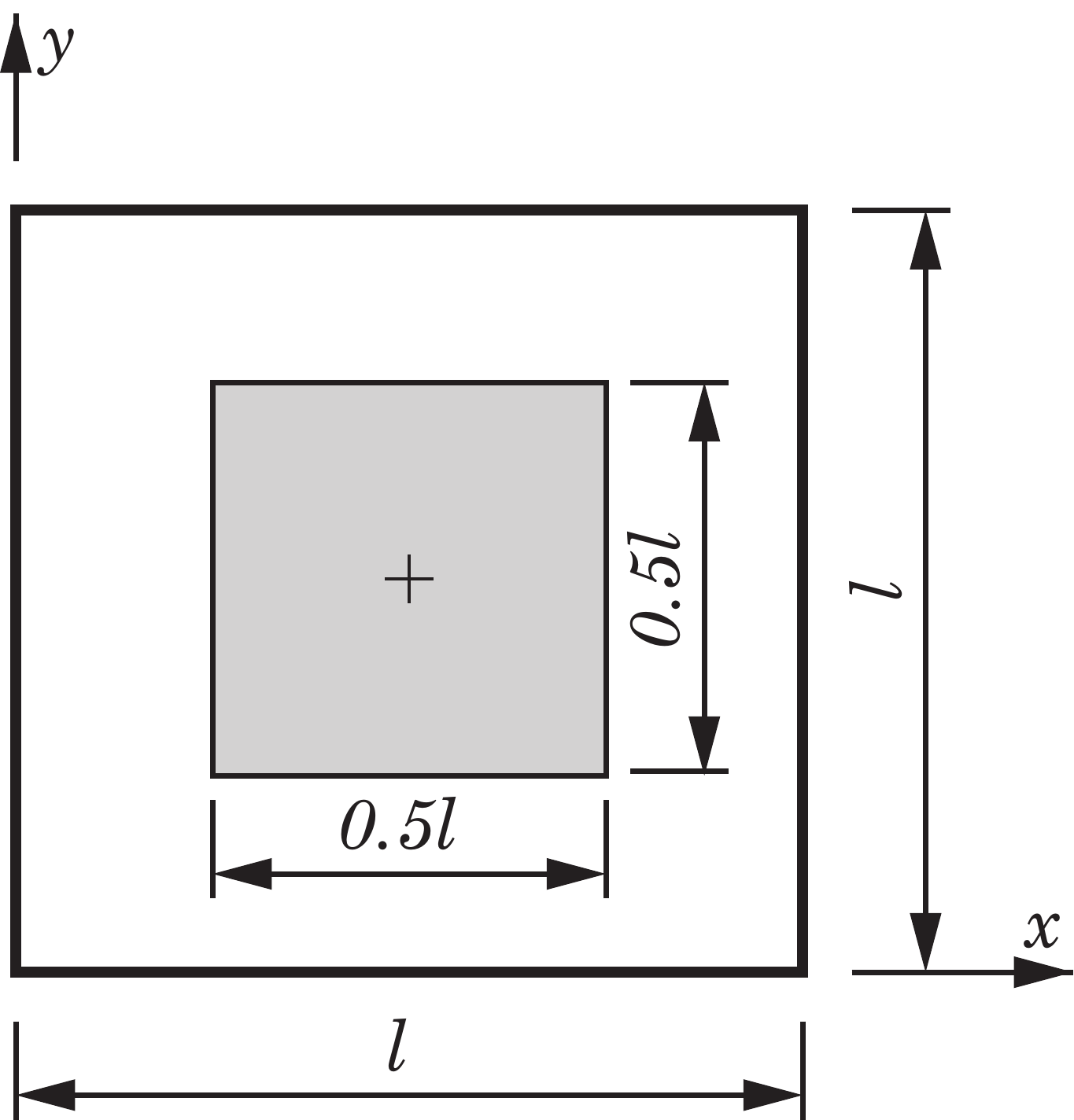}\caption{Scalar wave propagation in a square domain.\label{fig:Scalar-wave-propagation}}
\end{figure}

The problem of scalar wave propagation in a square domain \citep{Soares2020}
is shown in Fig.~\ref{fig:Scalar-wave-propagation}. The edges of
the square domain of dimension $l\times l$ are fixed. The wave speed
is denoted as $c$. An initial velocity $A$ is prescribed over an
area of $0.5l\times0.5l$ (shaded square in Fig.~\ref{fig:Scalar-wave-propagation})
at the middle of the domain. 
The analytical solutions for the displacement and velocity responses
are obtained by the method of separation of variables as\begin{subequations}
\label{eq:scalarWave-sln}
\begin{align}
u(x,y,t) & =\dfrac{16A}{\pi^{2}}\sum_{m=1}^{\infty}\sum_{n=1}^{\infty}\dfrac{1}{mn}\dfrac{1}{\mu_{mn}}\sin(\mu_{mn}t)\sin\frac{m\pi}{2}\sin\frac{m\pi}{4}\sin\frac{n\pi}{2}\sin\frac{n\pi}{4}\sin\frac{m\pi x}{l}\sin\frac{n\pi y}{l}\,,\label{eq:scalarWave-eign-sln-disp}\\
\dot{u}(x,y,t) & =\dfrac{16A}{\pi^{2}}\sum_{m=1}^{\infty}\sum_{n=1}^{\infty}\dfrac{1}{mn}\cos(\mu_{mn}t)\sin\frac{m\pi}{2}\sin\frac{m\pi}{4}\sin\frac{n\pi}{2}\sin\frac{n\pi}{4}\sin\frac{m\pi x}{l}\sin\frac{n\pi y}{l}\,,\label{eq:scalarWave-eign-sln-vel}
\end{align}
\end{subequations}with
\begin{equation}
\mu_{mn}=\frac{c\pi}{l}\sqrt{m^{2}+n^{2}}\,.\label{eq:scalarWave-eign}
\end{equation}

\begin{figure}
\begin{centering}
\subfloat[Proposed high-order scheme with $\rho_{\infty}=0.8$. \label{fig:Velocity-square-1}]{\includegraphics[width=0.45\textwidth]{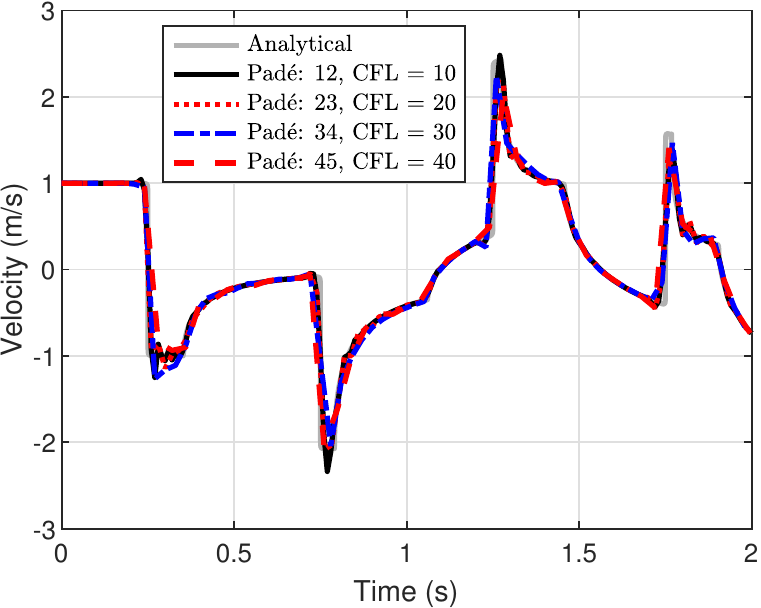}}\hfill{}\includegraphics[width=0.45\textwidth]{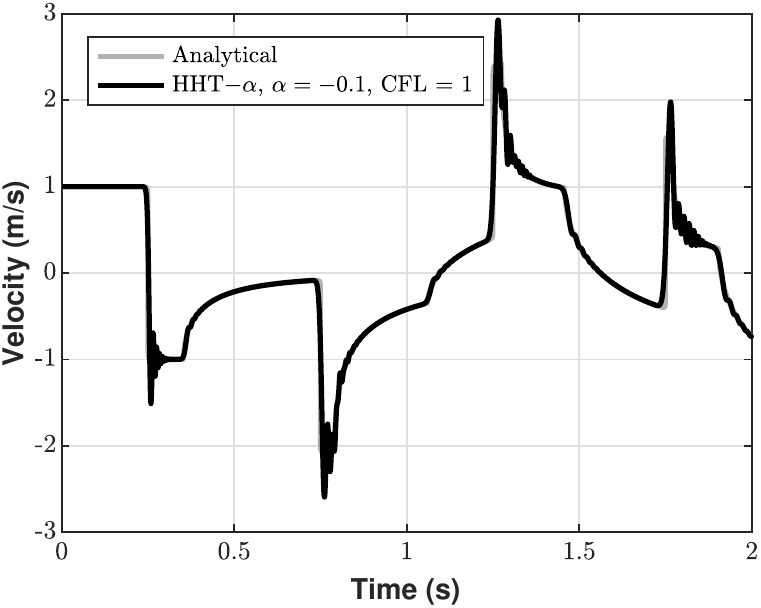}\caption{HHT-$\alpha$ method with $\alpha=-0.1$. \label{fig:Velocity-square-2}}
\par\end{centering}
\centering{}\caption{Velocity responses at center of square domain.\label{fig:Velocity-square}}
\end{figure}

In the following analyses, $l=\unit[1\,]{m}$, $c=\unit[1]{\,m/s}$,
and $A=\unit[1]{\,m/s}$ are chosen. Considering symmetry, only a
quarter of the domain ($0.5l\leq x\leq l$, $0.5l\leq y\leq l$) is
modeled. A mesh of $1000\times1000$ linear finite elements is generated.
The parameter $\rho_{\infty}=0.8$ is used to introduce numerical
dissipation. The CFL number is calculated using Eq.~\eqref{eq:CFL},
where $\Delta x$ is chosen as the length of a finite element equal
to $\Delta x=5\times10^{-4}l.$ The time step size is determined to
attain $\mathrm{CFL}=10L$ for the $(L=M-1,M)$ order scheme. The
velocity response histories at the center of the square domain, indicated
by a cross in Fig.~\ref{fig:Scalar-wave-propagation}, are plotted
in Fig.~\ref{fig:Velocity-square-1} together with the analytical
solution. The results of HHT-$\alpha$ method with $\alpha=-0.1$
and $\mathrm{CFL}=1$ are shown Fig.~\ref{fig:Velocity-square-2}.
It is again observed that the proposed high-order scheme is more effective
in dissipating spurious high-frequency oscillations. Results for this
example obtained with a high-order method are reported in Ref.~\citep{Soares2020}.

\subsection{Two-dimensional wave propagation in a semi-infinite elastic domain
- Lamb problem\label{subsec:Two-dimensional-wave-propagation}}

\begin{figure}
\begin{centering}
\includegraphics[width=0.4\textwidth]{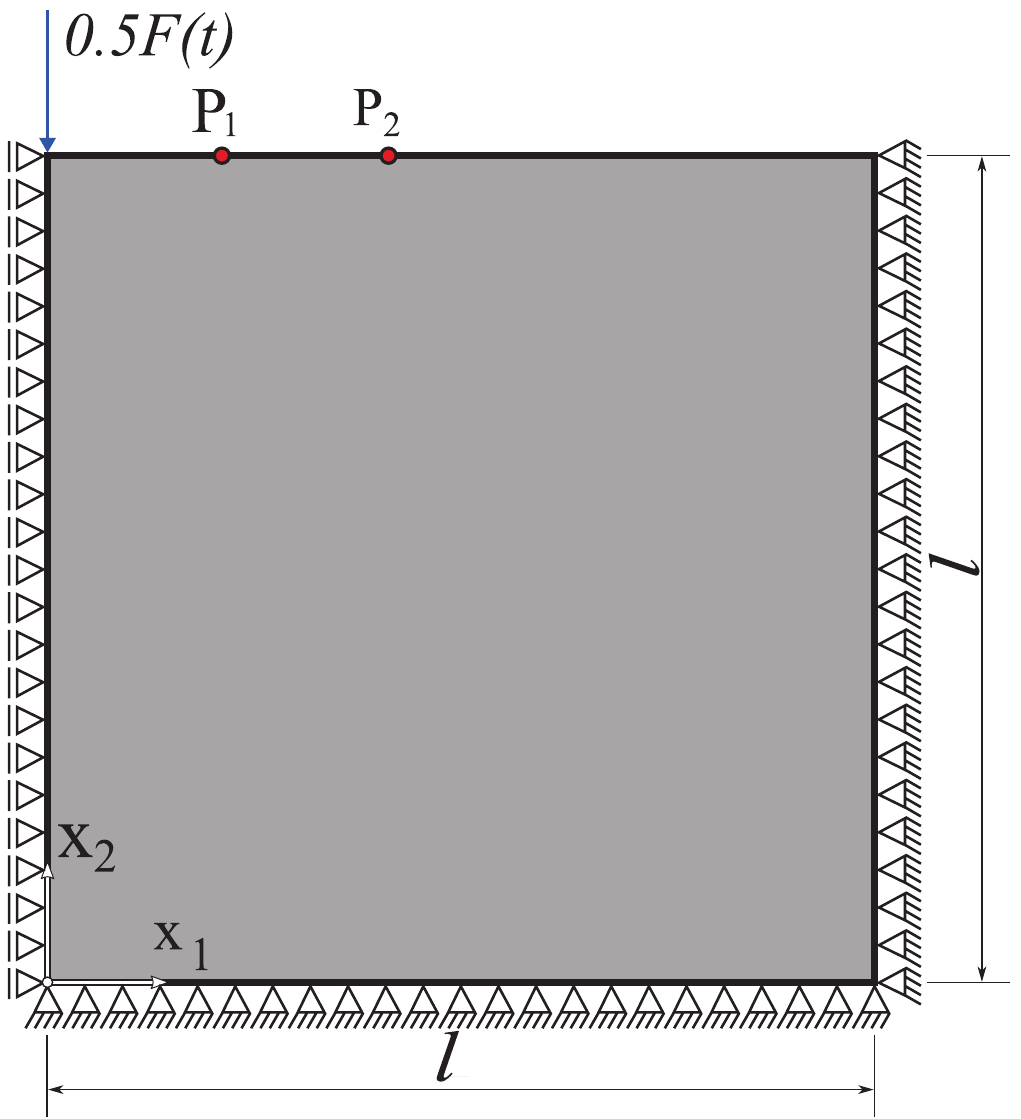}
\par\end{centering}
\centering{}\caption{A semi-infinite elastic domain in plane strain conditions \label{fig:geometry Lamb}}
\end{figure}

Lamb's problem with a vertical point load is analyzed in this section.
The square domain of dimension $l\times l$ representing the semi-infinite
elastic plane and boundary conditions are depicted in Fig.~\ref{fig:geometry Lamb}.
Due to symmetry, only the right side to the point load $F(t)$ is
considered. The geometry and materials parameters presented in Refs.~\citep{Kim2021}
and \citep{Kwon2020} are adopted: Length of the square domain $l=\unit[3{,}200]{\,m},$
Young's modulus $E=\unit[18.77\times10^{9}]{\,Pa},$ Poisson's ratio
$\nu=0.25$, and mass density $\rho=\unit[2{,}200]{\,kg/m^{3}}$.
Plane strain conditions are assumed. The P-wave, S-wave and Rayleigh-wave
velocities are equal to $c_{p}=\unit[3{,}200\,]{m/s}$, $c_{s}=\unit[1{,}847.5]{\,m/s}$,
and $c_{R}=\unit[1{,}698.6\,]{m/s}$, respectively. The time history
of the point load consists of three step functions: $F(t)=2\times10^{6}(H(0.15-t)-3H(0.1-t)+3H(0.05-t))\,\unit{N}$.
The analytical solution of the horizontal and vertical displacements
at the two points, $P_{1}\,(640,3{,}200)$ and $P_{2}\,(1{,}280,3{,}200)$,
indicated in Fig.~\ref{fig:geometry Lamb} are given in~\citep{Kim2021}
until $t=\unit[1\,]{s},$ when the P-waves reach the right side of
the domain.

The square domain is divided into a uniform mesh of $2{,}000\times2{,}000$
linear finite elements. The side length $\Delta x$ of each square
element is $\unit[1.6]{m}$. In the calculation of the CFL number,
the P-wave velocity $c_{p}$ and side length $\Delta x$ are used.

\begin{figure}
\subfloat[Horizontal displacement at $P_{1}$]{\includegraphics[width=0.45\textwidth]{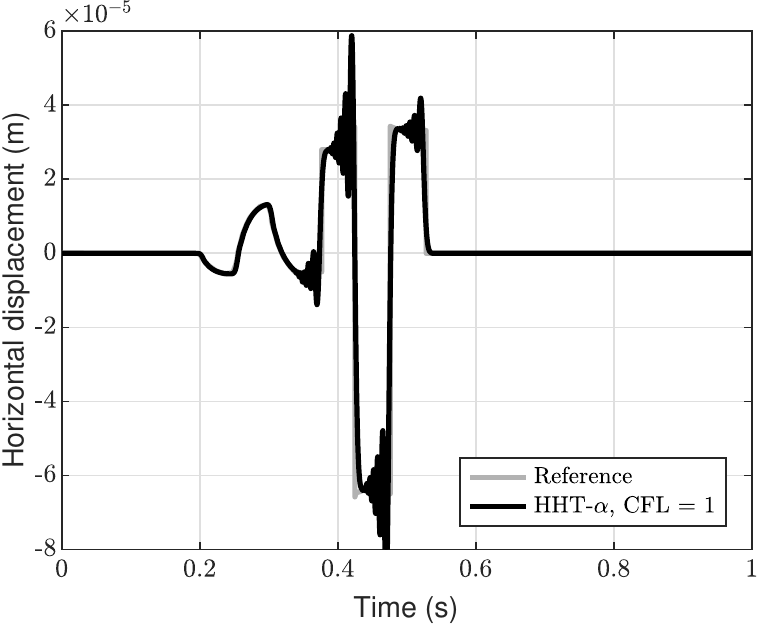}}\hfill{}\subfloat[Horizontal displacement at $P_{2}$]{\includegraphics[width=0.45\textwidth]{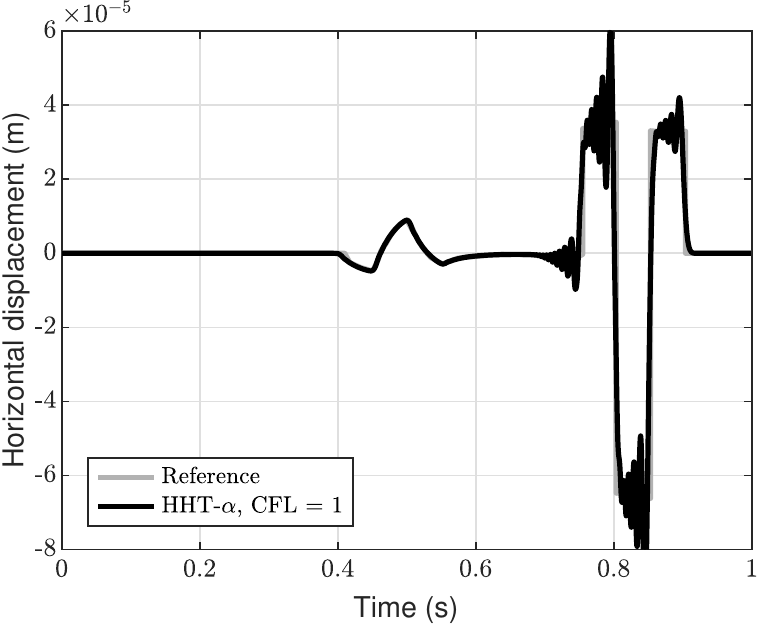}}

\vspace*{\medskipamount}
\subfloat[Vertical displacement at $P_{1}$]{\includegraphics[width=0.45\textwidth]{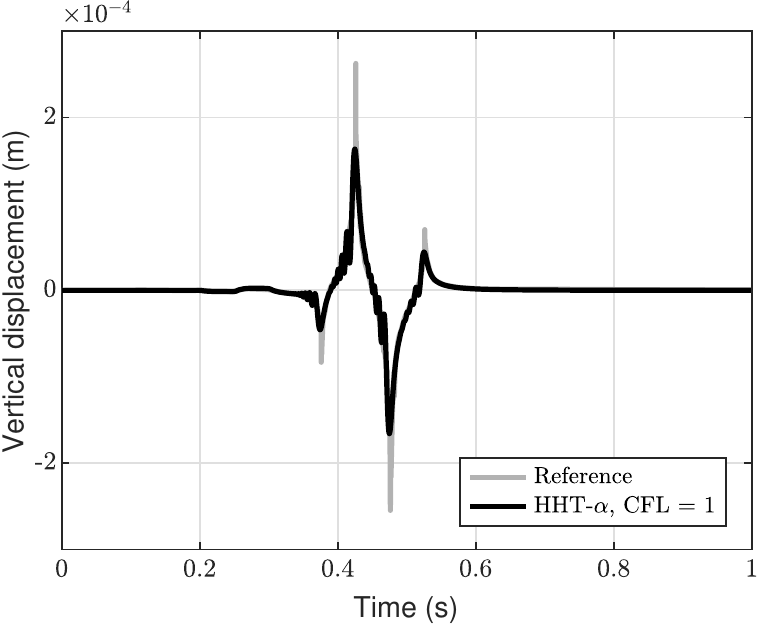}}\hfill{}\subfloat[Vertical displacement at $P_{2}$]{\includegraphics[width=0.45\textwidth]{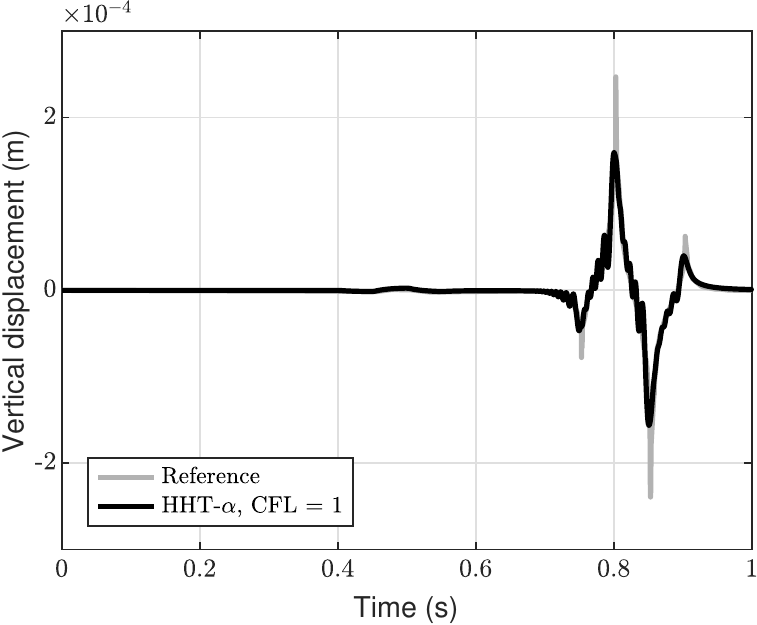}}

\caption{Horizontal and vertical displacements of an semi-infinite elastic
domain obtained by HHT-$\alpha$ method. \label{fig:LambdisplacementHHT}}
\end{figure}

\begin{figure}
\subfloat[Horizontal velocity at $P_{1}$]{\includegraphics[width=0.45\textwidth]{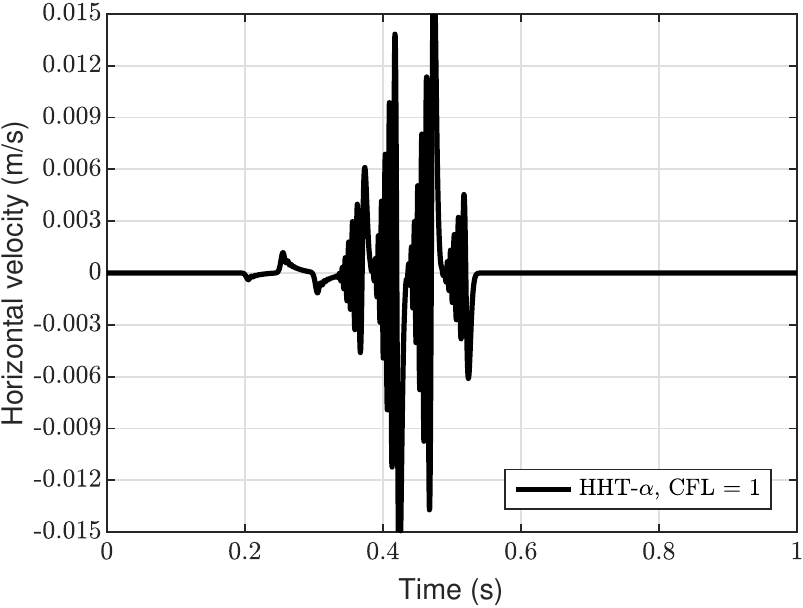}}\hfill{}\subfloat[Horizontal velocity at $P_{2}$]{\includegraphics[width=0.45\textwidth]{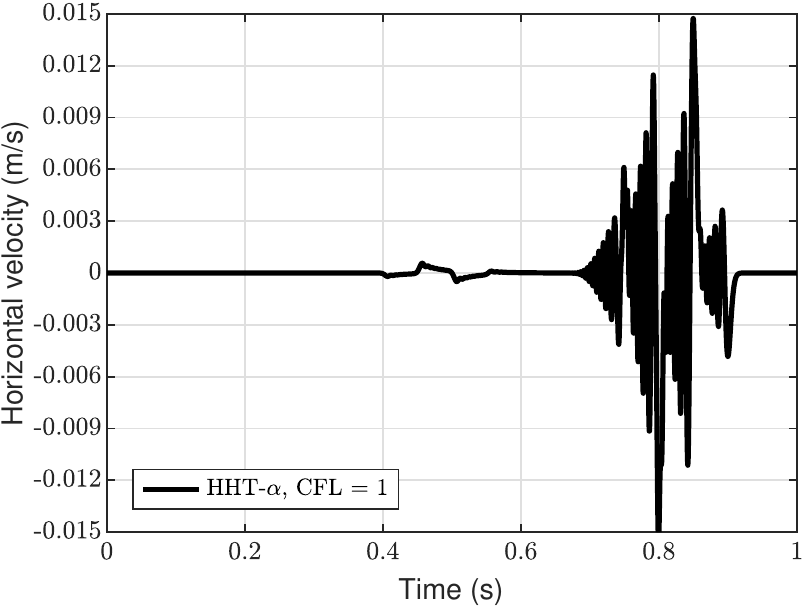}}

\vspace*{\medskipamount}
\subfloat[Vertical velocity at $P_{1}$]{\includegraphics[width=0.45\textwidth]{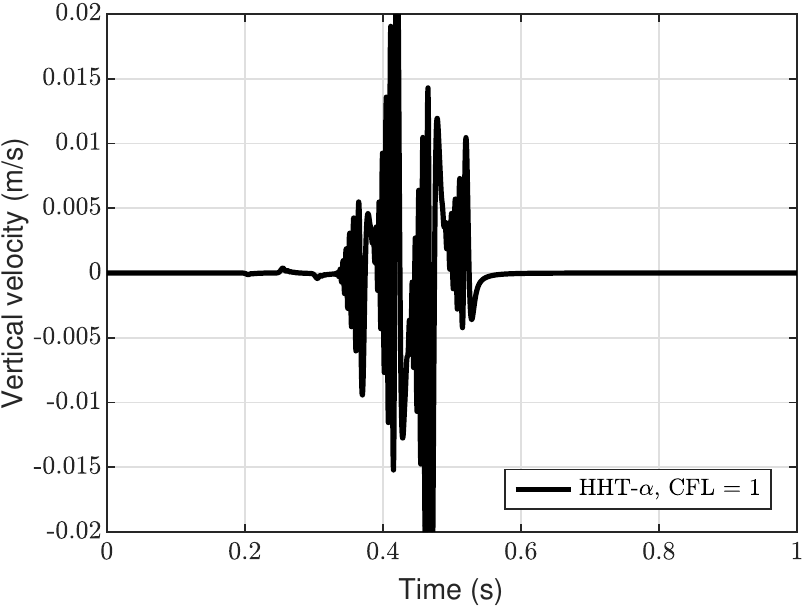}}\hfill{}\subfloat[Vertical velocity at $P_{2}$]{\includegraphics[width=0.45\textwidth]{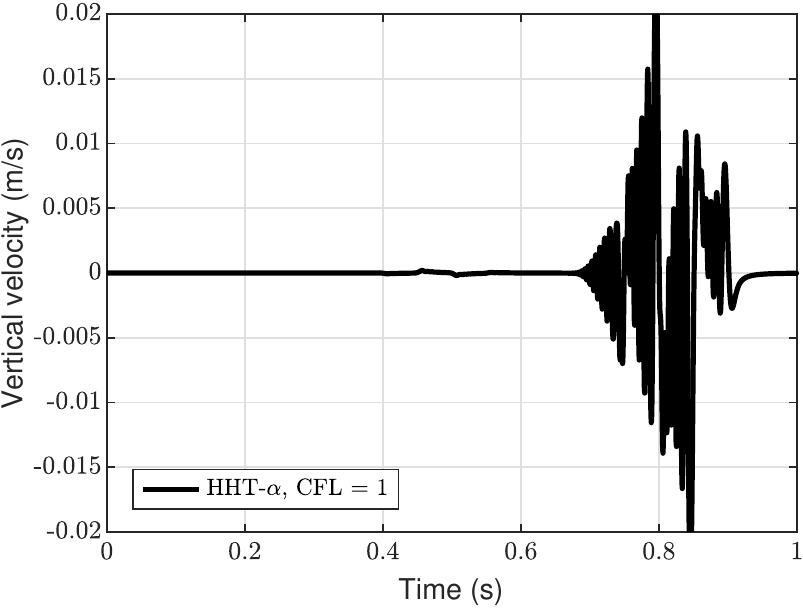}}

\caption{Horizontal and vertical velocities of a semi-infinite elastic domain
obtained by HHT-$\alpha$ method. \label{fig:LambvelocityHHT}}
\end{figure}

An analysis using the HHT-$\alpha$ method with $\alpha=-0.1$ and
$\mathrm{CFL=1}$ is also performed. The displacements and velocity
responses are plotted in Figs.~\ref{fig:LambdisplacementHHT} and
\ref{fig:LambvelocityHHT}, respectively. The peak responses of vertical
displacements are about 5 times larger than those of horizontal displacements.
The displacements show good agreement with the analytical solution.
Some spurious oscillations of the smaller horizontal displacements
are observed. However, the velocity responses in Fig.~\ref{fig:LambvelocityHHT}
exhibit very strong high-frequency oscillations. The presence of waves
traveling at different speeds renders the numerical dissipation of
the HHT-$\alpha$ method much less effective than in the example of
scalar waves (Section~\ref{subsec:Scalar-wave-propagation}). The
displacement responses obtained with the Bathe method are available
in \citep{Kim2021} and \citep{Kwon2020} for comparison. 

\begin{figure}
\subfloat[Horizontal displacement at $P_{1}$]{\includegraphics[width=0.45\textwidth]{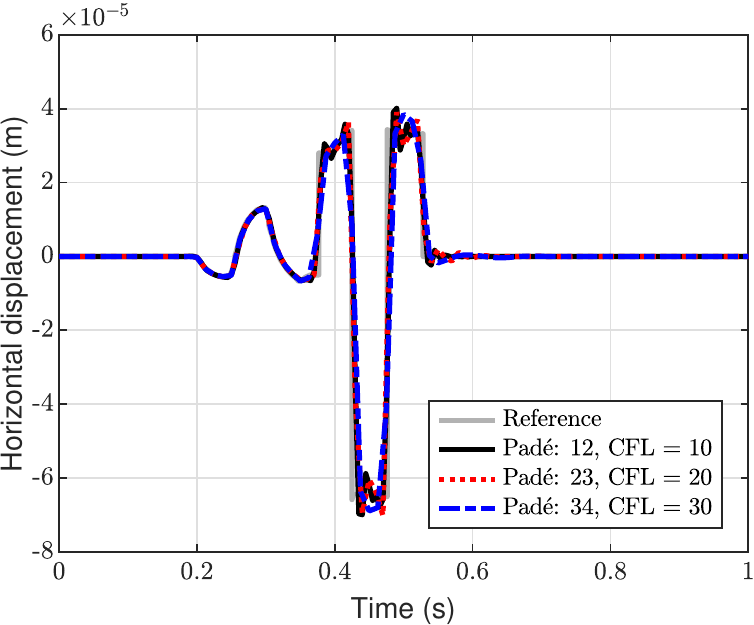}}\hfill{}\subfloat[Horizontal displacement at $P_{2}$]{\includegraphics[width=0.45\textwidth]{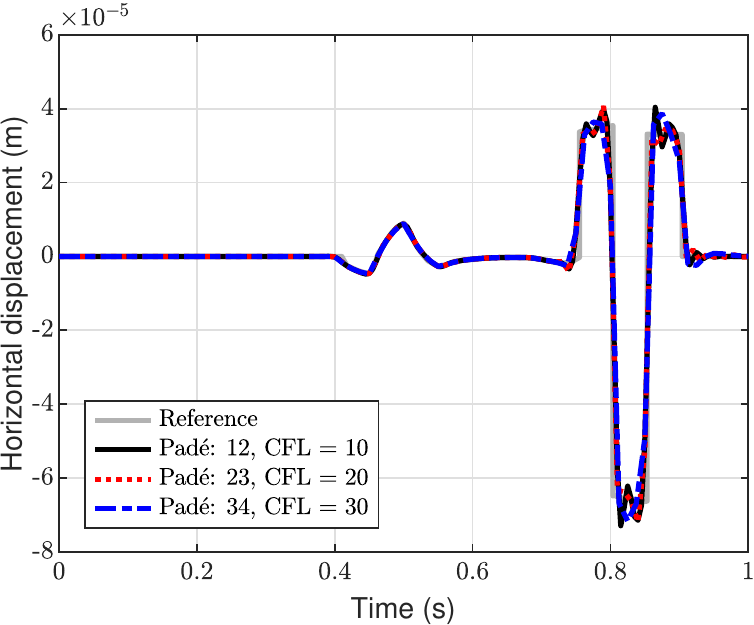}}

\vspace*{\medskipamount}
\subfloat[Vertical displacement at $P_{1}$]{\includegraphics[width=0.45\textwidth]{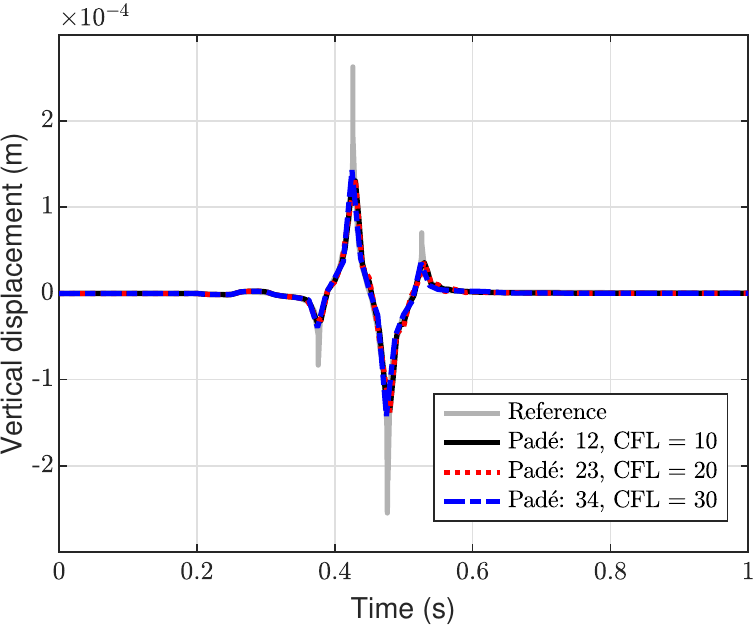}}\hfill{}\subfloat[Vertical displacement at $P_{2}$]{\includegraphics[width=0.45\textwidth]{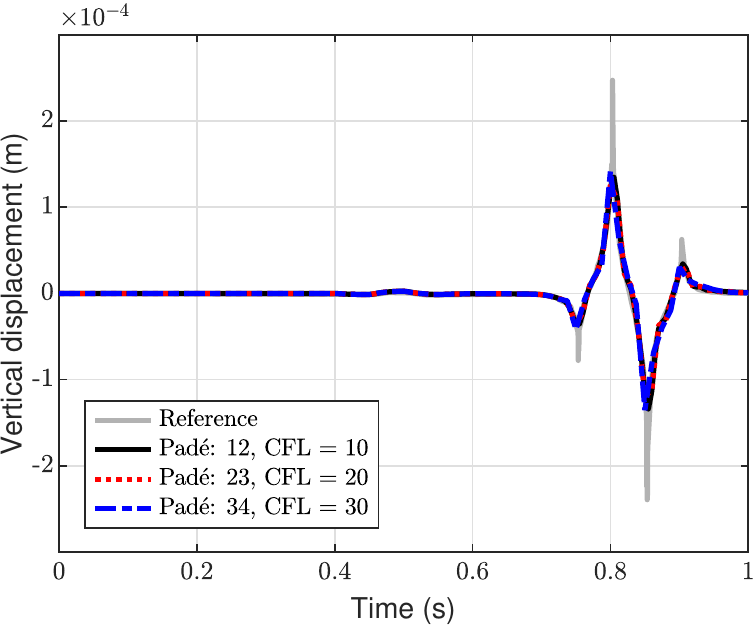}}

\caption{Horizontal and vertical displacements of a semi-infinite elastic domain
obtained by the proposed high-order scheme. \label{fig:Lambdisplacement}}
\end{figure}

\begin{figure}
\subfloat[Horizontal velocity at $P_{1}$]{\includegraphics[width=0.45\textwidth]{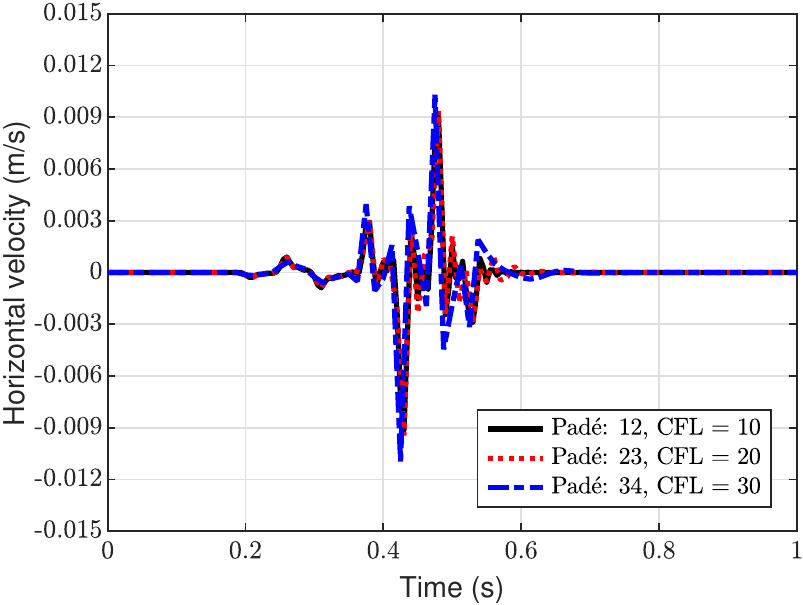}}\hfill{}\subfloat[Horizontal velocity at $P_{2}$]{\includegraphics[width=0.45\textwidth]{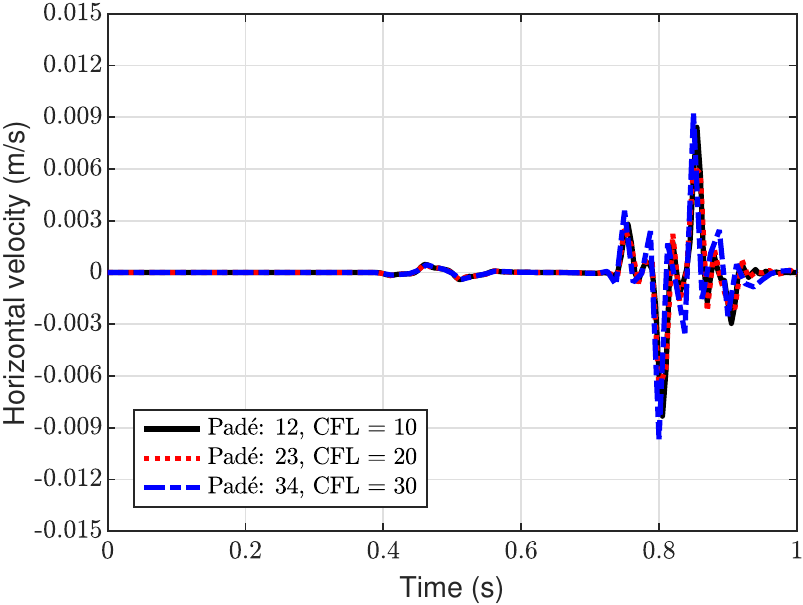}}

\vspace*{\medskipamount}
\subfloat[Vertical velocity at $P_{1}$]{\includegraphics[width=0.45\textwidth]{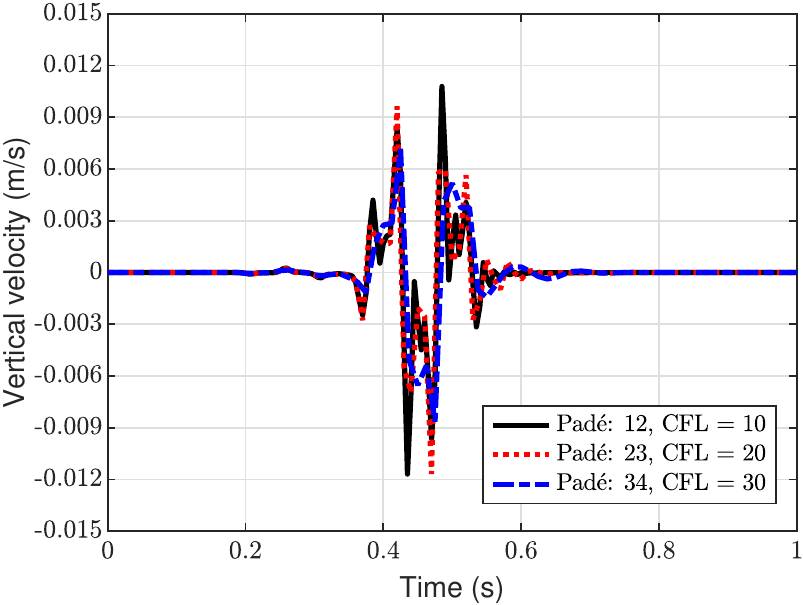}}\hfill{}\subfloat[Vertical velocity at $P_{2}$]{\includegraphics[width=0.45\textwidth]{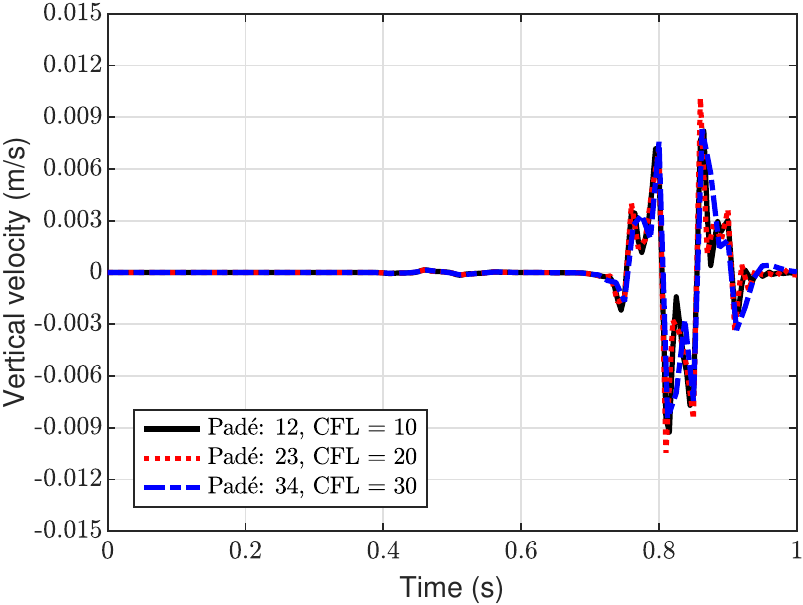}}

\caption{Horizontal and vertical velocities of the Lamb problem obtained by
the proposed high-order scheme. \label{fig:Lambvelocity}}
\end{figure}

In the analysis using the high-order scheme, the parameter $\rho_{\infty}=0.8$
is adopted to introduce numerical dissipation. The CFL numbers for
the time-integration scheme at orders $(1,2)$, $(2,3)$ and $(3,4)$
are selected as $10$, $20$ and $30$, respectively. Correspondingly,
the time step sizes are given by $\Delta t=\Delta x\times\mathrm{CFL}/c_{p}$
as $\unit[0.005]{s}$, $\unit[0.01]{s}$ and $\unit[0.015]{s}$. The
displacement and velocity responses are plotted in Fig.~\ref{fig:Lambdisplacement}
and Fig.~\ref{fig:Lambvelocity}, respectively. Very good agreement
of the displacement response with the analytical solution is observed.
Spurious oscillations in the velocity response are largely suppressed.


\FloatBarrier

\subsection{Three-dimensional wave propagation in 3D sandwich panel\label{subsec:Three-dimensional-wave-propagati}}

\begin{figure}
\centering{}\subfloat[Geometry\label{fig:geometry_3d}]{\includegraphics[width=0.45\textwidth]{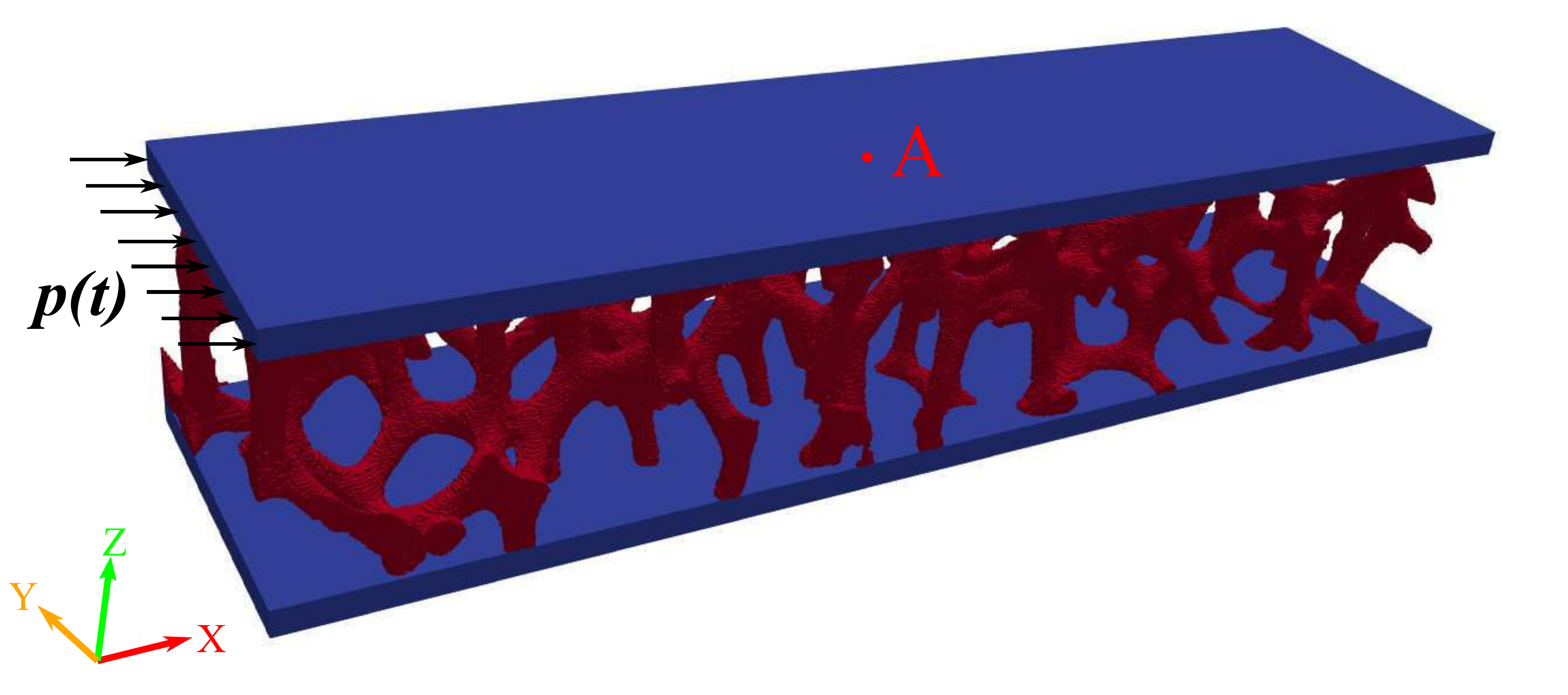}}\hfill{}\subfloat[Octree mesh\label{fig:Octree-mesh_3d}]{\includegraphics[width=0.45\textwidth]{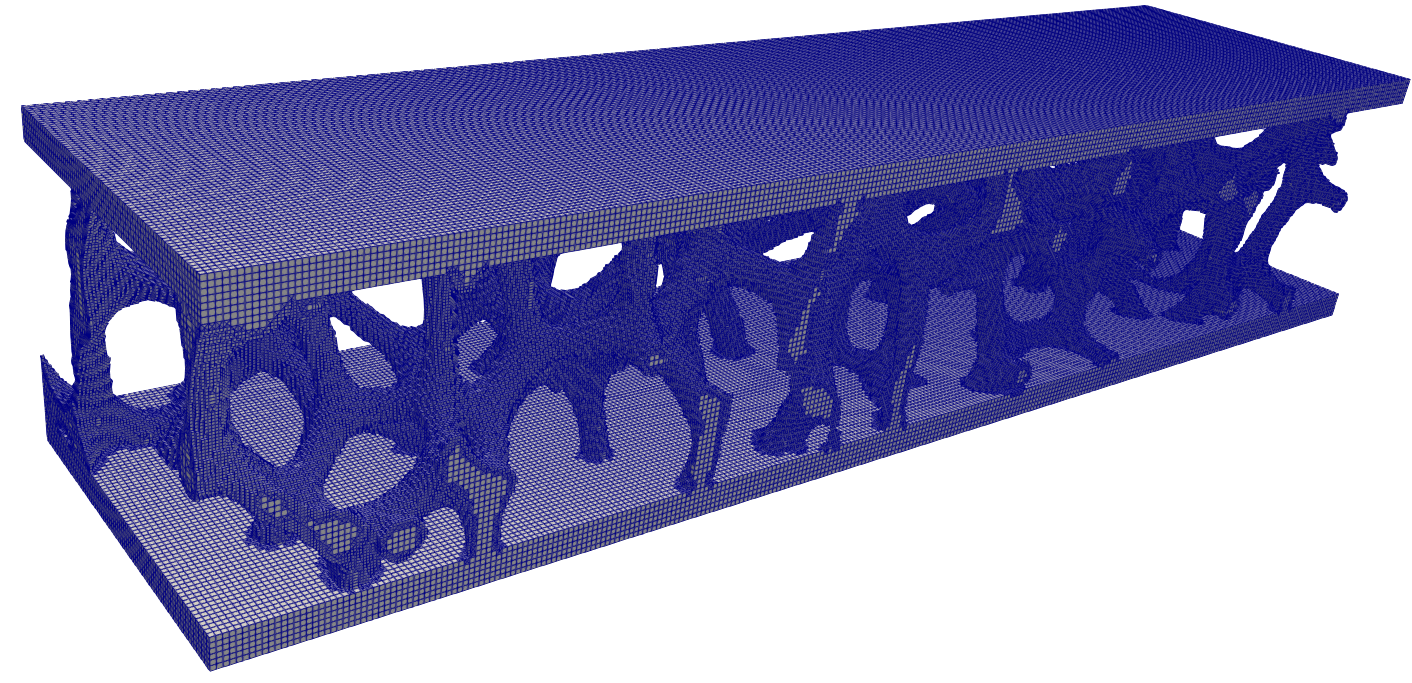}}\caption{\foreignlanguage{australian}{Sandwich panel.\label{fig:sandwich panel}}}
\end{figure}

A sandwich panel with two cover-sheets and a foam-core is shown in
Fig.~\ref{fig:geometry_3d}. The outer dimension of the panel is
$288\,\unit{mm}\times72\,\unit{mm}\times59.5\,\unit{mm}$ and the
thickness of the cover sheets is $5.76\,\unit{mm}$. The cover-sheets
are made of steel with the following properties: Young's modulus $E_{\mathrm{s}}=210\,\unit{GPa}$,
Poisson's ratio $\nu_{\mathrm{s}}=0.3$, and the mass density $\rho_{\mathrm{s}}=8050\,\unitfrac{kg}{m^{3}}$.
Therefore, the P- and S-wave speeds are $c_{p}=5926\,\mathrm{m/s}$
and $c_{s}=3168\,\mathrm{m/s}$, respectively. The foam-core is given
as a digital image obtained by X-ray CT scans. The material of the
foam is aluminium with the property $E_{\mathrm{a}}=70\,\unit{GPa}$,
$\nu_{\mathrm{a}}=0.3$ and $\rho_{\mathrm{a}}=\unit[2700]{kg/m^{3}}$.
The corresponding wave speeds are $c_{p}=6198\,\mathrm{m/s}$ and
$c_{s}=3122\,\mathrm{m/s}$, respectively. The right end of the panel
is fixed in the normal direction and the left end of the top cover-sheet
is subjected to a uniformly distributed pressure $p(t)$. The time
history of the pressure consists of two step functions: $p(t)=1\times(2H(1.5\times10^{-5}-t)-H(3\times10^{-5}-t))\,\unit{kPa}$.

The sandwich panel is discretized by an octree mesh as shown in Fig.~\ref{fig:Octree-mesh_3d}
and modelled by the scaled boundary finite element method \citep{Song2018,Zhang2021}.
The smallest and largest element sizes are $0.48\,\unit{mm}$ and
$1.92\,\unit{mm}$, respectively. Overall, the mesh consists of 597,325
elements, 1,099,242 nodes, and consequently, 3,297,726 degrees of
freedom.

The time integration using the proposed high-order schemes is performed
with the parameter $\rho_{\infty}=0.8$ to suppress the oscillations
in the response due to the high-frequency components in the excitation.
The time step size is selected based on the CFL number in the steel
cover-plates, where the element size is $0.96\,\unit{mm}$. The simulation
is performed until $t_{\mathrm{sim}}=6\times10^{-4}\,{\textstyle \mathrm{s}}$,
which allows the waves to be reflected at the two ends for several
times. As for the examples in the previous sections, the CFL number
is chosen as 10 for order (1,2) and 20 for order (2,3). The corresponding
time step size is equal to $\Delta t=1.5\times10^{-6}\,{\textstyle \mathrm{s}}$
for order (1,2), resulting in 400 time steps, and $\Delta t=3\times10^{-6}\,{\textstyle \mathrm{s}}$
for order (2,3) with 200 time steps. 

\begin{figure}
\subfloat[Displacement in $x-$direction\label{fig:Displacement-responses}]{\begin{centering}
\includegraphics[width=0.95\textwidth]{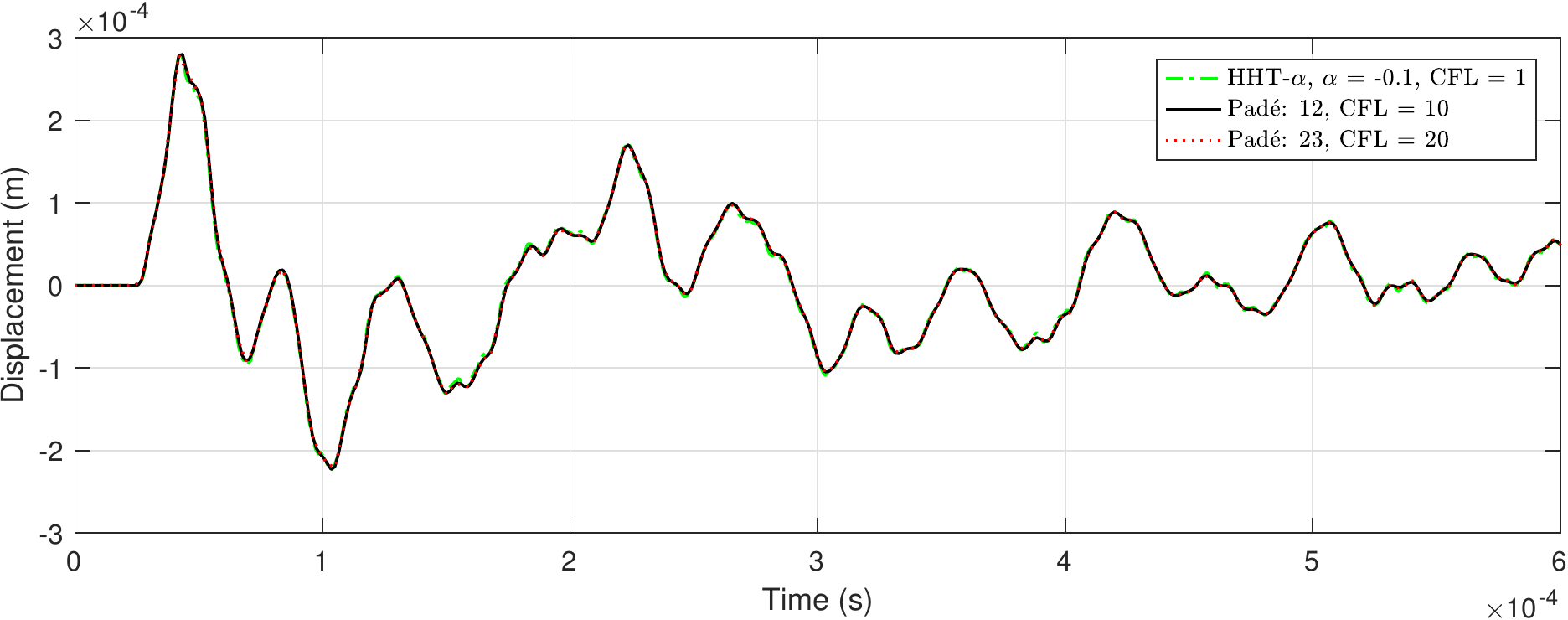}
\par\end{centering}
}

\subfloat[Velocity in $x-$direction\label{fig:Velocity-responses}]{\begin{centering}
\includegraphics[width=0.95\textwidth]{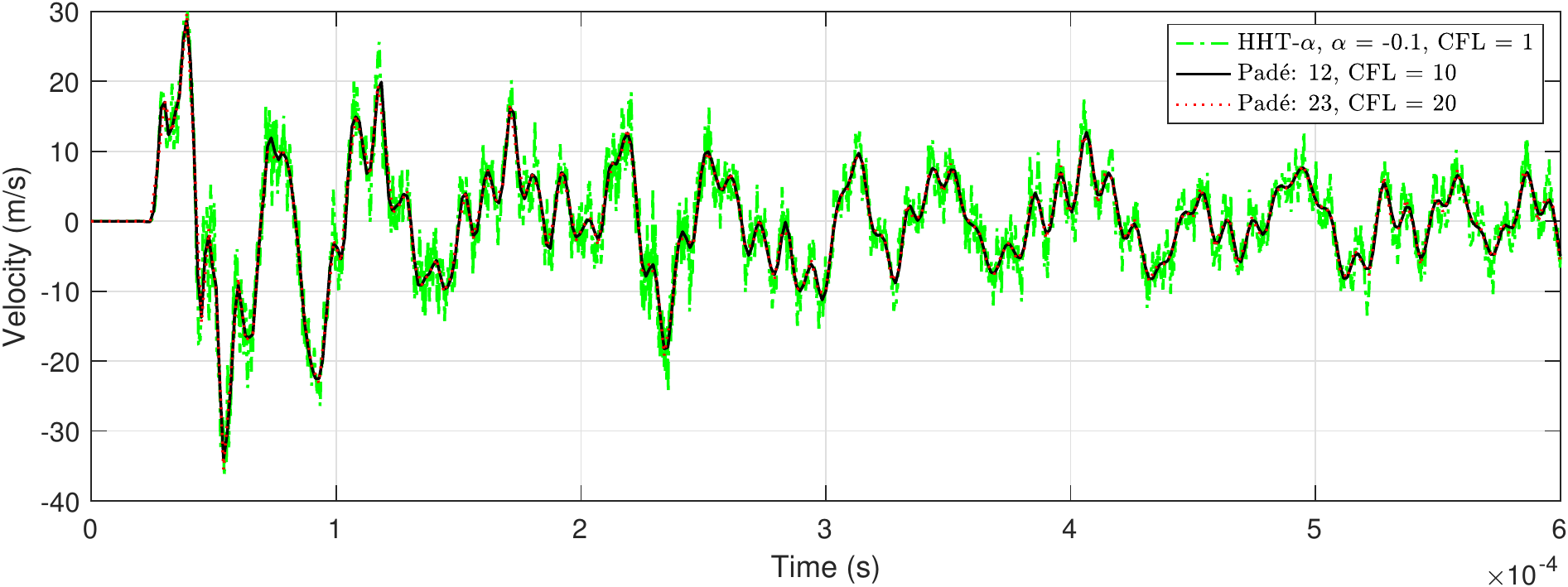}
\par\end{centering}
}

\caption{Responses of sandwich panel at point $\mathrm{A}$ obtained by the
proposed high-order scheme with $\rho_{\infty}=0.8$ and HHT-$\alpha$
method with $\alpha=-0.1$.}
\end{figure}

The displacement and velocity responses along the $x-$direction at
the middle point of the top surface $\mathrm{A}\!\left(144,36,57.6\right)$
(unit: mm), indicated by the red dot in Fig.~\ref{fig:geometry_3d},
are plotted in Fig.~\ref{fig:Displacement-responses} and Fig.~\ref{fig:Velocity-responses},
respectively. The black solid line and red dotted line represent the
responses obtained using the proposed scheme at orders (1,2) and (2,3),
respectively. The results are nearly identical to each other. An analysis
using HHT-$\alpha$ method with $\alpha=-0.1$ and CFL = 1 (with $\Delta t=0.15\times10^{-6}\,{\textstyle \mathrm{s}}$
and $4,000$ time steps) is also performed. The result is indicated
by the green dashed line. The displacement response in Fig.~\ref{fig:Displacement-responses}
is in very good agreement with those of the proposed scheme, but the
velocity in Fig.~\ref{fig:Velocity-responses} shows strong spurious
oscillations.

\begin{figure}
\subfloat[$t=2.4\times10^{-5}\,{\textstyle \mathrm{s}}$]{\includegraphics[width=0.45\textwidth]{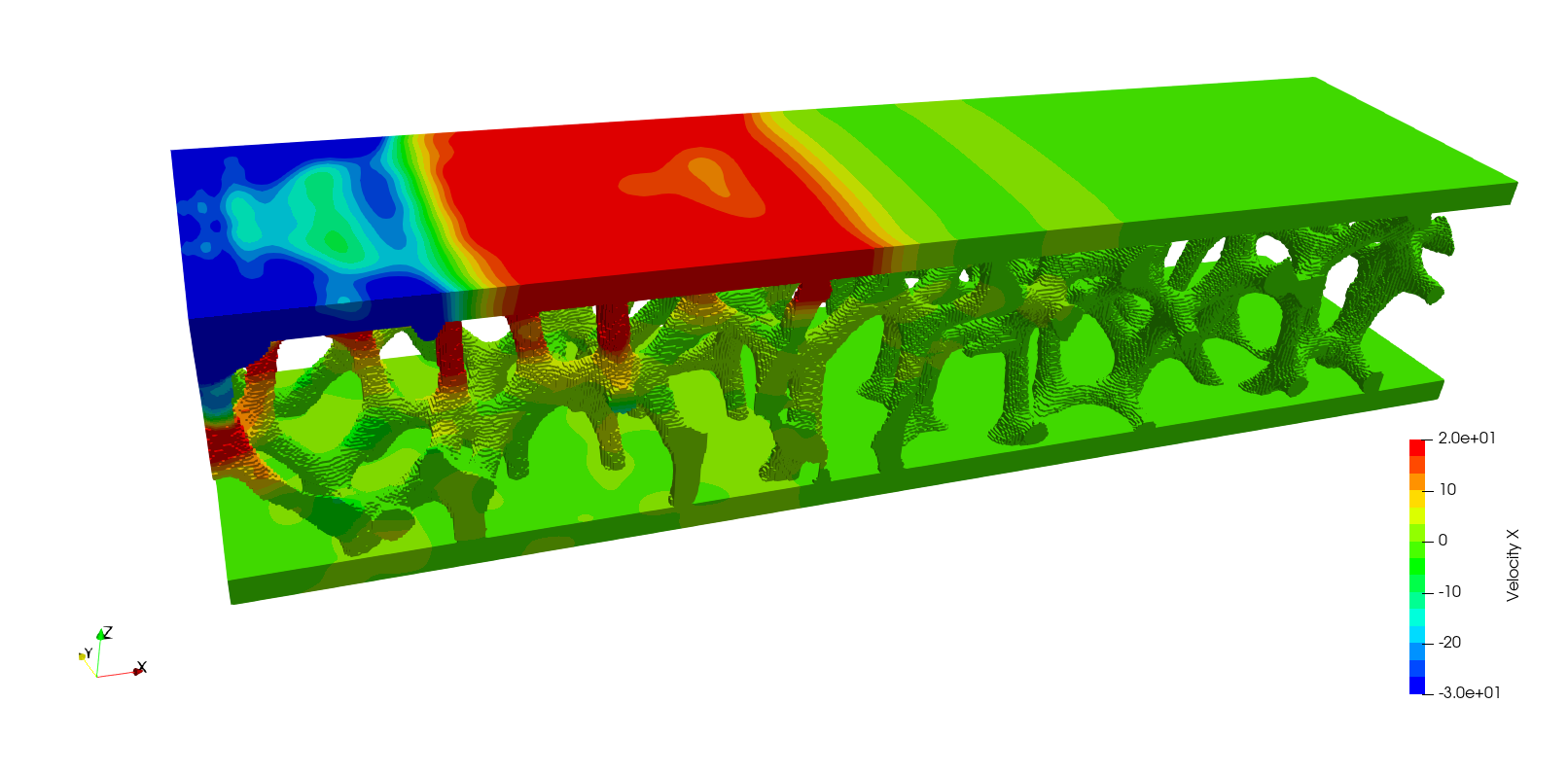}}\hfill{}\subfloat[$t=4.8\times10^{-5}\,{\textstyle \mathrm{s}}$]{\includegraphics[width=0.45\textwidth]{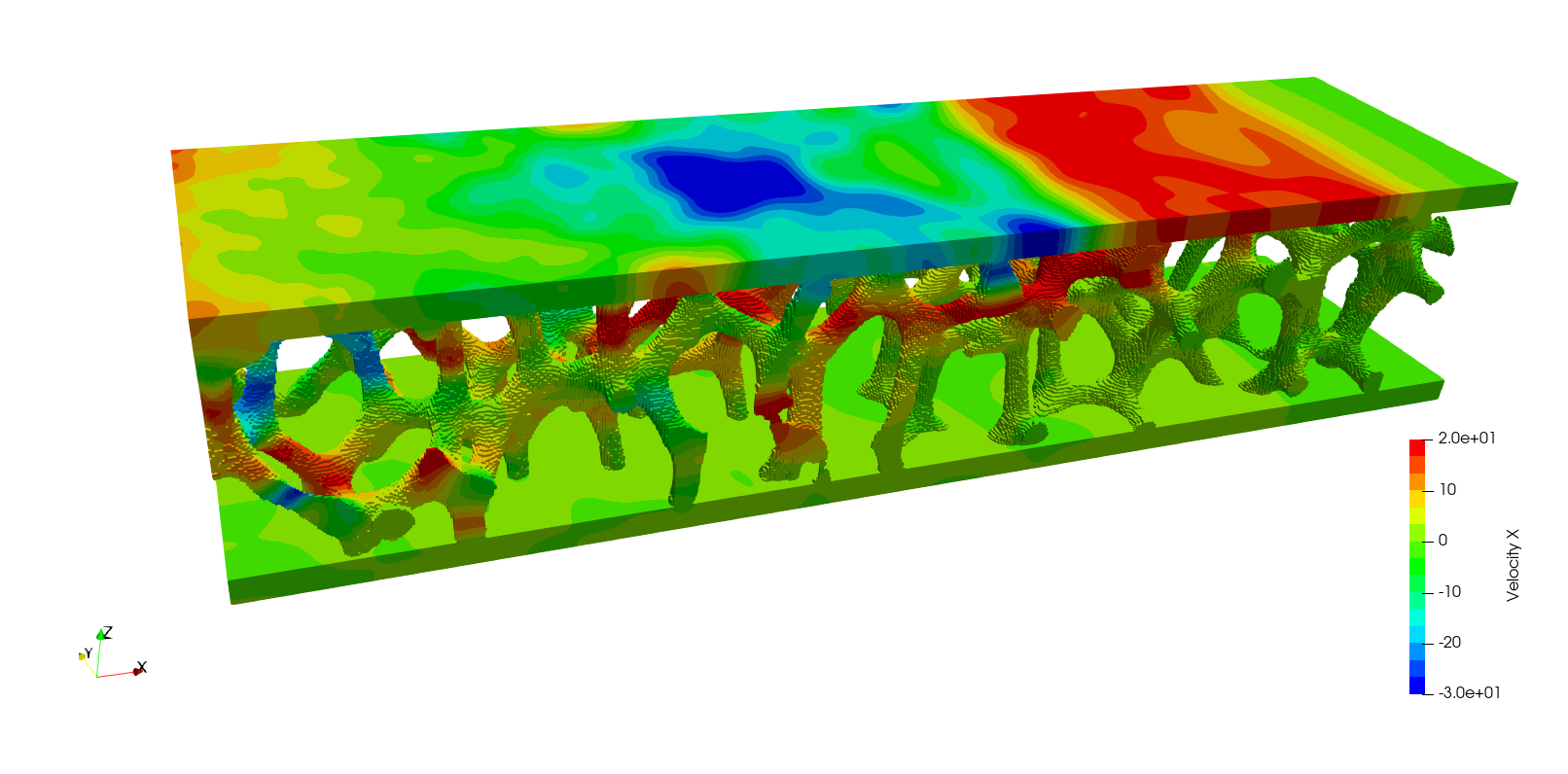}}

\subfloat[$t=7.2\times10^{-5}\,{\textstyle \mathrm{s}}$]{\includegraphics[width=0.45\textwidth]{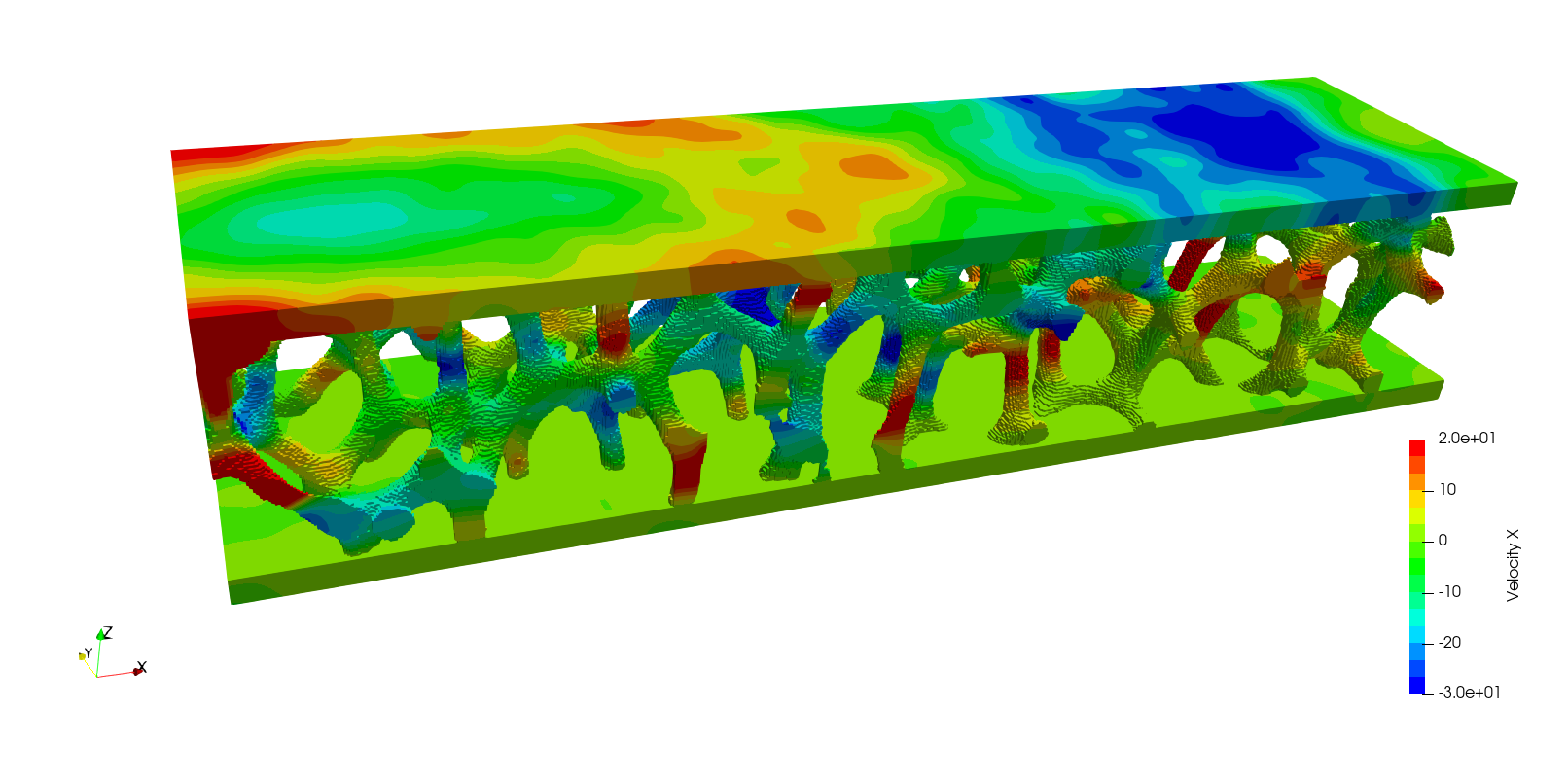}}\hfill{}\subfloat[$t=9.6\times10^{-5}\,{\textstyle \mathrm{s}}$]{\includegraphics[width=0.45\textwidth]{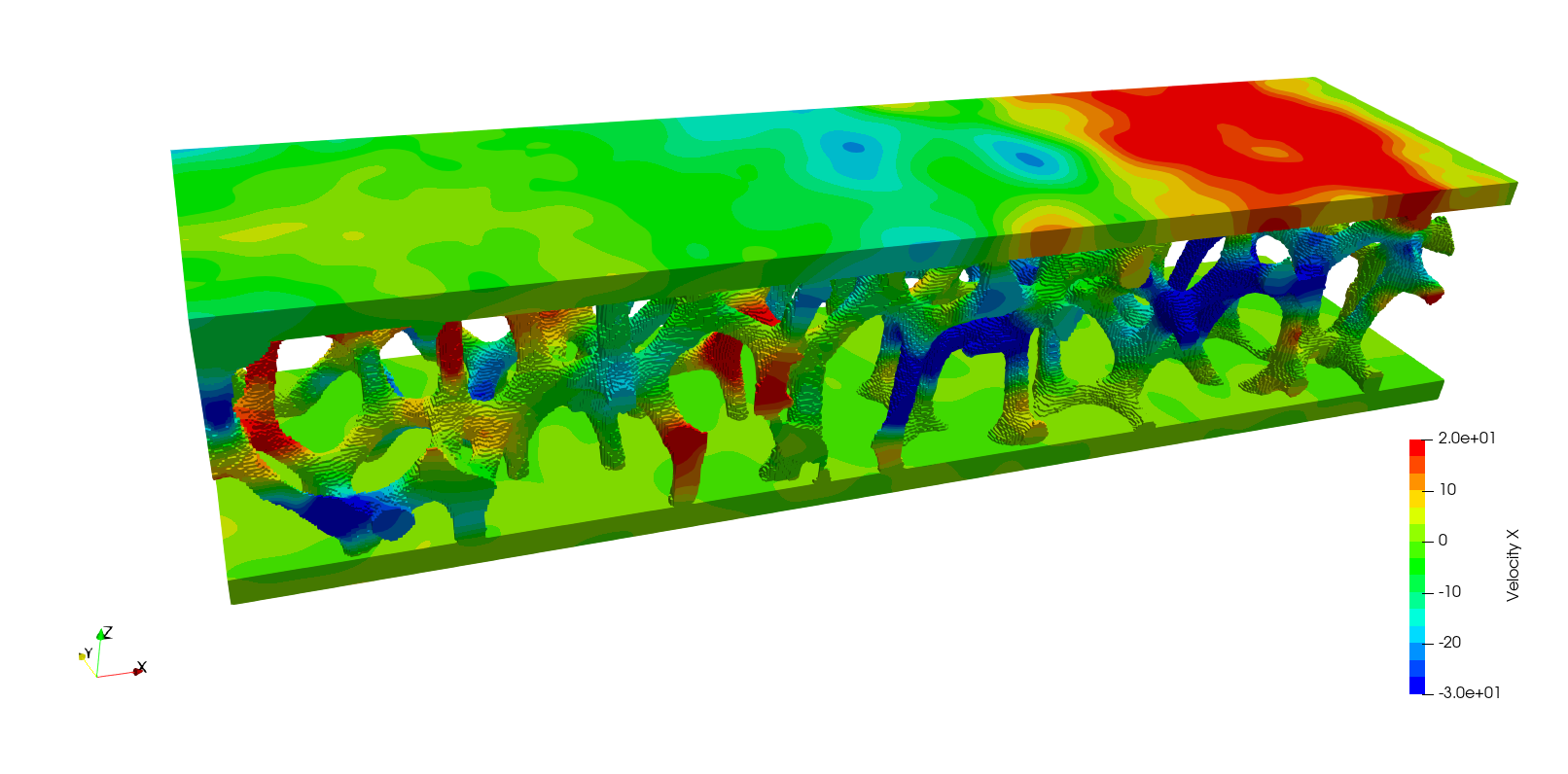}}

\subfloat[$t=3.14\times10^{-4}\,{\textstyle \mathrm{s}}$]{\includegraphics[width=0.45\textwidth]{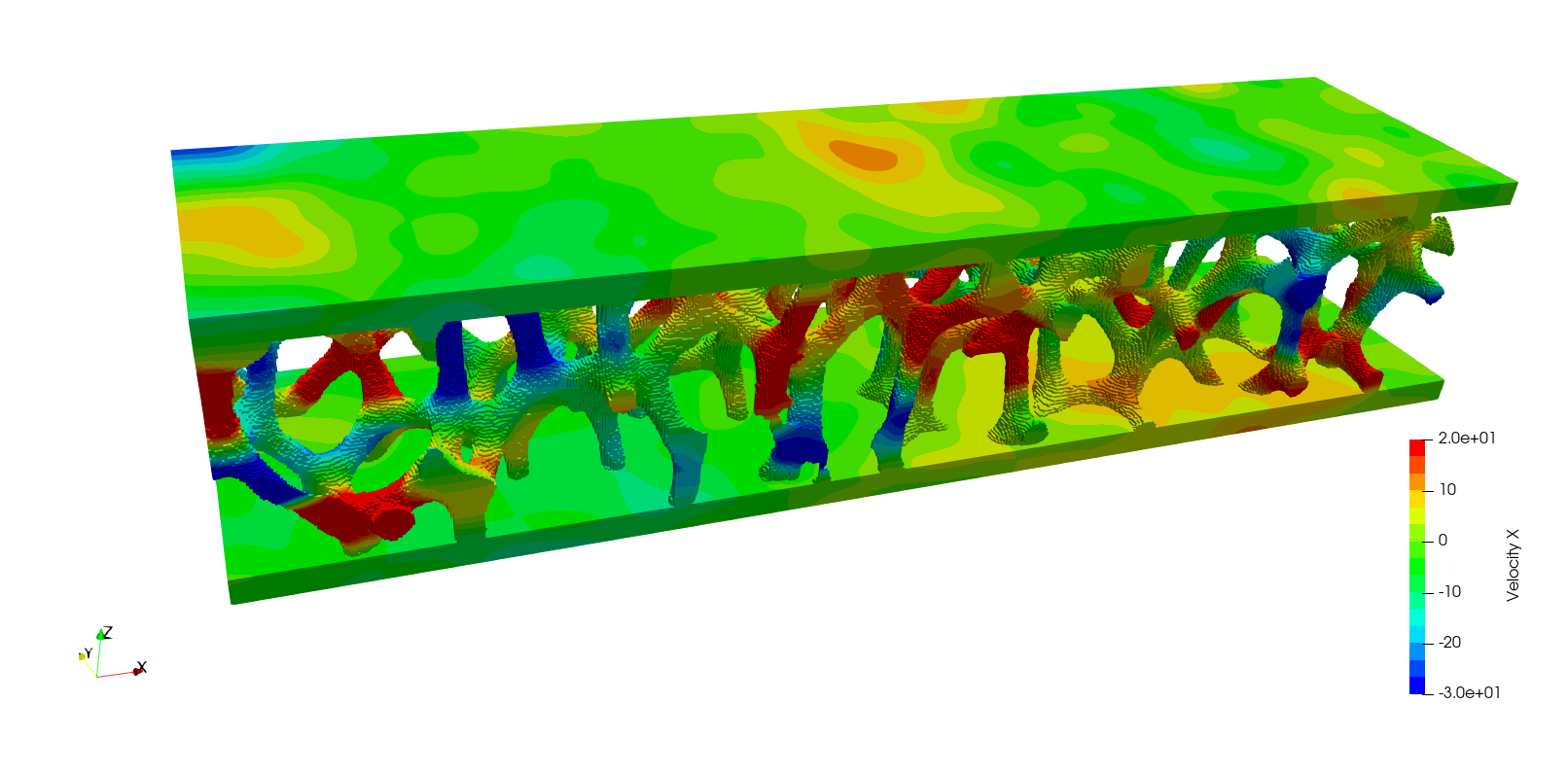}}\hfill{}\subfloat[$t=5.84\times10^{-4}\,{\textstyle \mathrm{s}}$]{\includegraphics[width=0.45\textwidth]{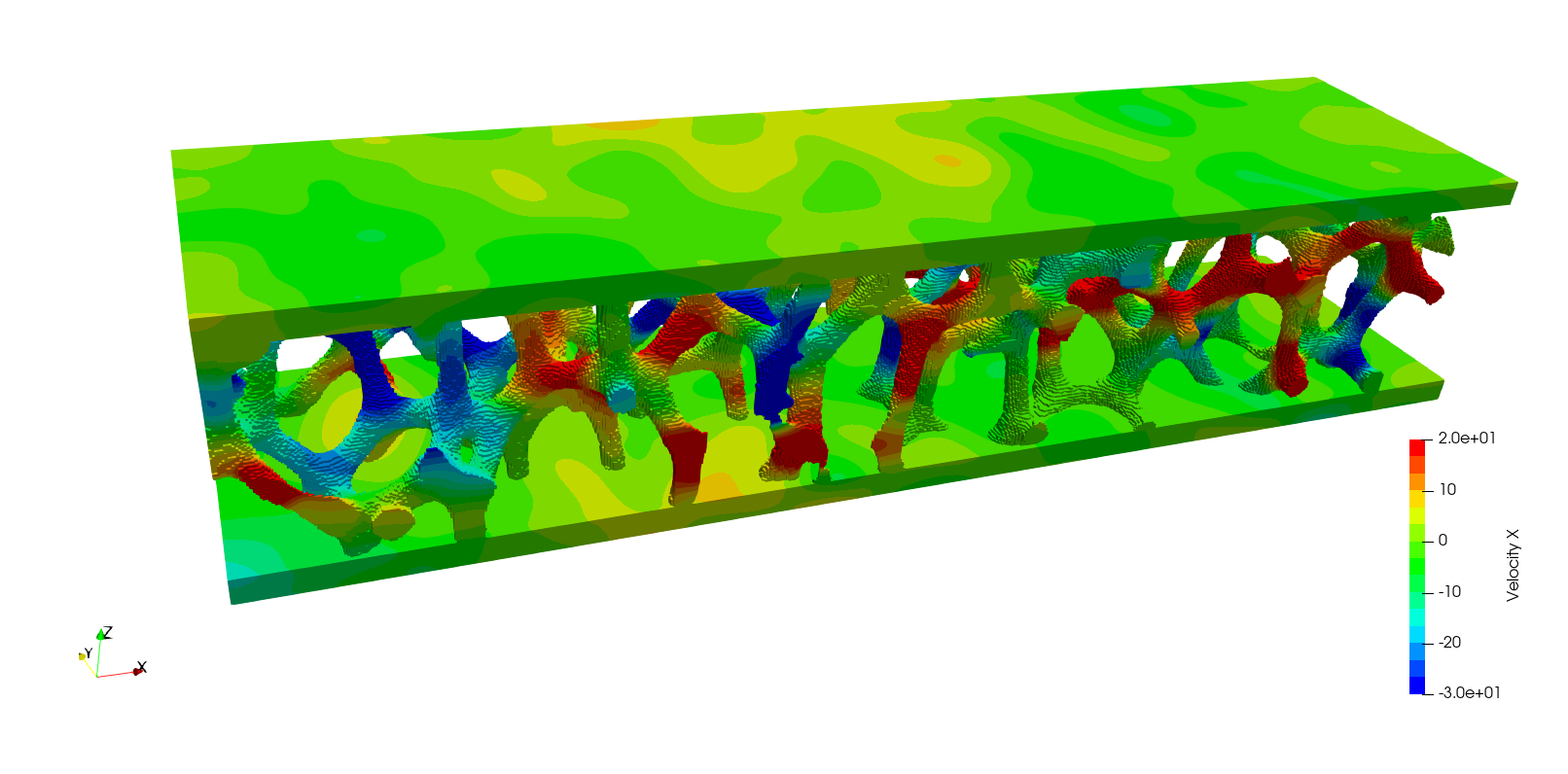}}

\caption{Contours of velocity along the $x-$direction of sandwich panel at
various time instances\label{fig:Velocity-contours}}
\end{figure}

The contours of velocity along the $x-$direction at six selected
time instances are presented in Fig.~\ref{fig:Velocity-contours}.
The waves initially concentrate in the upper cover-plate, and gradually
excite the lower cover-plate through the foam core.

The computational costs of the proposed scheme are evaluated for this
example. The computer program for time integration is written in \texttt{FORTRAN}.
The \texttt{PARDISO} direct solver of Intel’s Math Kernel Library
(MKL) is employed for the solution of simultaneous linear algebraic
equations (Eq.~\eqref{eq:realRootEq5} for a real root and Eq.~\eqref{eq:ComplexRoot_sln_y1}
for a complex root). This operation takes the majority of the running
time. For this example of elasto-dynamics, the factorization of matrices
is performed once at the beginning of the analysis. Only back-substitutions
are performed during the time stepping. The computer running times
are measured on a Dell Precision 5820 Tower Workstation with an Intel(R)
Xeon(R) W-2275 CPU and 256~GB RAM. The HHT-$\alpha$ method takes
4,007\,s for time stepping. The proposed scheme takes 960\,s at
order (1,2) and 857\,s at order (2,3), which corresponds to speedup
factors of about 4.18 and 4.69, respectively, in comparison with the
HHT-$\alpha$ method.

\section{Conclusions\label{sec:Conclusions}}

A high-order implicit time integration scheme is proposed. The amount
of numerical dissipation is controlled by using the spectral radius
$\rho_{\infty}$ at the high-frequency limit as a user-specified parameter.
The scheme varies with the specified parameter from \textit{A}-stable
(without numerical dissipation) to \textit{L}-stable (with the maximum
amount of numerical dissipation). The numerical dissipation is minimal
at the low-frequency range and rapidly approaches the maximum value
at the high-frequency range, showing better characteristics than the
second-order HHT-$\alpha$ method. Moreover, the period error is smaller
in the low-frequency range than that of the HHT-$\alpha$ method by
orders of magnitude.

From the viewpoint of application, the only user-specified parameters
are the value of the spectral radius in the high-frequency limit,
$\rho_{\infty}$, and the time step size, $\Delta t$, expressed in
terms of the CFL number. Effective dissipation of spurious high-frequency
oscillations can be achieved in a wide range of the parameters: $0\leq\rho_{\infty}\leq0.9$
and $5L\leq\mathrm{\textrm{CFL}}\leq20L$ for a scheme of order $(L=M-1,M)$
when finite elements of linear shape functions are used for spatial
discretization. The values $\rho_{\infty}=0.8$ and $\mathrm{\textrm{CFL}}=10L$
are used for wave propagation problems in this paper. 

An efficient numerical algorithm is designed, where the systems of
equations to be solved are similar in complexity to those in the standard
Newmark method. Existing computer codes of the Newmark, HHT-$\alpha$,
and generalize-$\alpha$ methods can be extended straightforwardly
to include the proposed high-order scheme. When compared with the
HHT-$\alpha$ method for the same finite element model, the proposed
scheme is not only more effective in dissipating spurious high-frequency
oscillations but also reduces computer running time. A speedup factor
of more than 4 is observed when solving the sandwich panel problem
employing our \texttt{FORTRAN} code with the Intel MKL \texttt{PARDISO}
direct solver.

\section*{Acknowledgments}

The work presented in this paper is partially supported by the Australian
Research Council through Grant Number DP200103577. The authors would
also like to thank Dr. Meysam Joulaian and Professor Alexander Düster
from Hamburg University of Technology for providing the X-ray CT scan
data, which was used in Section~\ref{subsec:Three-dimensional-wave-propagati}.

\bibliographystyle{elsarticle-num}
\bibliography{literature}

\end{document}